\documentclass{amsart}

\usepackage{caption}
\usepackage{graphicx}
\usepackage{amssymb}
\usepackage{amsmath}
\usepackage{amsthm}                   
\usepackage{amsfonts}                    
\usepackage[normalem]{ulem}              
\usepackage[noadjust,verbose,sort]{cite} 
\usepackage{setspace}                    
\usepackage{tikz}
\usepackage{array}
\newtheorem{theorem}{Theorem}[section]
\newtheorem{definition}[theorem]{Definition}
\newtheorem{lemma}[theorem]{Lemma}

\newtheorem{proposition}[theorem]{Proposition}
\newtheorem{corollary}[theorem]{Corollary}

\title[Geometric \& Combinatorial Properties of the ASM Polytope]{Geometric and Combinatorial Properties of the Alternating Sign Matrix Polytope}
\date{October 2026}
\author{Elizabeth A. Dinkelman}
\address{Department of Mathematical Sciences,
George Mason University, Fairfax, Virginia, USA and
Northern Virginia Community College, Fairfax, Virginia, USA}
\email{edinkelm@gmu.edu}

\author{Walter D. Morris, Jr.}
\address{Department of Mathematical Sciences,
George Mason University, Fairfax, Virginia, USA}
\email{wmorris@gmu.edu}

\begin{document}

\begin{abstract}
The polytope $ASM_n$, the convex hull of the $n\times n$ alternating sign matrices, was introduced by Striker and by Behrend and Knight.   A face of $ASM_n$ corresponds to an elementary flow grid defined by Striker, and each elementary flow grid determines a doubly directed graph defined by Brualdi and Dahl.  We show that a face of $ASM_n$ is centrally symmetric if and only if its doubly directed graph has all vertices of even degree.   We show that every face of $ASM_n$ is a 2-level polytope.  We show that a $d$-dimensional face of $ASM_n$ has at most $2^d$ vertices and $4(d-1)$ facets, for $d\ge 2$.  We show that a $d$-dimensional face of $ASM_n$ satisfies $vf\le d2^{d+1}$, where $v$ and $f$ are the numbers of vertices and facets of the face.   If the doubly directed graph of a $d$-dimensional face is 2-connected, then $v\le 2^{d-1}+2$.  We describe the facets of a face and a basis for the subspace parallel to a face in terms of the elementary flow grid of the face. We prove that no face of $ASM_n$ has the combinatorial type of the Birkhoff polytope $B_3$.  We list the combinatorial types of faces of $ASM_n$ that have dimension 4 or less.
\end{abstract}
\maketitle

\section{Introduction}    


Mills, Robbins, and Rumsey \cite{Mills1} introduced alternating sign matrices and posed the problem of counting the $n\times n$ alternating sign matrices. 

\begin{definition}
Alternating sign matrices ($ASM$s) are square matrices that have entries in $\{ 0,1,-1 \}$, where each row and column sum to 1, and where the nonzero entries in each row and column alternate in sign.
\end{definition}

If the matrix 
\[
D_3=\left[ { \begin{array}{ccc}
    0 & 1 & 0 \\
    1 & -1 & 1 \\
    0 & 1 & 0 \\
  \end{array} } \right]
  \]
is included with the six $3\times 3$ permutation matrices, we obtain the seven $3\times 3$ alternating sign matrices.

 An interesting account of what led to the $ASM$ conjecture of the enumeration of $n \times n \ ASMs$ proposed by Mills, Robbins and Rumsey \cite{Mills1}, and proved by Zeilberger \cite{Zeilberger1} and Kuperberg \cite{Kuperberg1} can be found in Bressoud's \textit{Proofs and Confirmations,  The Story of the Alternating Sign Matrix Conjecture} \cite{Bressoud1}. An independent proof was given by Fischer \cite{Fischer1}. The number of $n \times n$ alternating sign matrices is 
$$ \prod_{j=0}^{n-1} \frac{(3j+1)!}{(n+j)!}.$$

Striker \cite{Striker1} defined the polytope $ASM_n$ as the convex hull of the $n\times n$ alternating sign matrices.     
Striker  \cite[Theorem 2.1]{Striker1} and independently Behrend and Knight \cite[Theorem 1]{Behrend1}, proved that the $ASM$s are the vertices of the solution set of a certain system of inequalities and equations:
\begin{theorem}\label{inequalitydescription}
The convex hull of $n \times n$ alternating sign matrices consists of all $n 
\times n$ real matrices $X = \{ x_{ij} \}$ such that:
\newline
$$0 \leq \sum_{i=1}^{i'} x_{ij} \leq 1 \ \ \forall \ \ 1 \leq i' \leq n, 1 \leq j \leq n$$
$$0 \leq \sum_{j=1}^{j'} x_{ij} \leq 1 \ \ \forall \ \ 1 \leq j' \leq n, 1 \leq i \leq n$$
$$\sum_{i=1}^nx_{ij}=1 \ \ \forall \ \ 1\leq j \leq n$$
$$\sum_{j=1}^nx_{ij}=1 \ \ \forall \ \ 1\leq i \leq n$$
\end{theorem} 

  A short proof of the fact that the vertices of $ASM_n$ are the $n\times n$ $ASM$s was later discovered by Brualdi and Dahl \cite[Theorem 1]{Brualdi4}.




 Striker \cite[Theorem 3.3]{Striker1} demonstrated that $ASM_n$ has $4[(n-2)^2+1]$ facets for $n \geq 3$.  Its dimension is $(n-1)^2$ \cite[Page 3]{Striker1} and \cite[Page 22]{Behrend1}.


Analogously, the combinatorial geometric properties of the Birkhoff polytope $B_n$, the convex hull of the set of the $n\times n$ permutation matrices, have been the subject of many investigations \cite{Balinski1,BruGibI,BruGibII,BruGibIII,Billera2,Paffenholz1}, and many striking facts are known.  For example, the $n\times n$ permutation matrices are the vertices of the set of $n\times n$ matrices satisfying $0\le x_{ij}\le 1$ for $1\le i\le n, 1\le j\le n$, $\sum_{i=1}^nx_{ij}=1 \ \ \forall \ \ 1\leq j \leq n$, and 
$\sum_{j=1}^nx_{ij}=1 \ \ \forall \ \ 1\leq i \leq n$. \cite{Birkhoff1}.  
$B_n$ has $n^2$ facets, $n!$ vertices and dimension $(n-1)^2$ \cite[Page 1]{Paffenholz1}.
A prominent tool for studying $B_n$ is the combinatorial correspondence between the faces of $B_n$ and the {\it elementary subgraphs} of $K_{n,n}$ whose edge sets are the unions of the edge sets of perfect matchings \cite{Balinski1,Billera2,Paffenholz1}.

 \begin{figure}[h!] 
\centering 
\begin{tikzpicture}
\draw [stealth-](0.5,1) -- (1,1);
\draw [stealth-](0.5,2) -- (1,2);
\draw [stealth-](0.5,3) -- (1,3);

\draw [-stealth](3,1) -- (3.5,1);
\draw [-stealth](3,2) -- (3.5,2);
\draw [-stealth](3,3) -- (3.5,3);

\draw [-stealth](1,3) -- (1,3.5);
\draw [-stealth](2,3) -- (2,3.5);
\draw [-stealth](3,3) -- (3,3.5);

\draw [-stealth](1,1) -- (1,0.5);
\draw [-stealth](2,1) -- (2,0.5);
\draw [-stealth](3,1) -- (3,0.5);

\draw [stealth-](1,1) -- (2,1);
\draw [-stealth](1,2) -- (2,2);
\draw [stealth-](1,3) -- (2,3);

\draw [-stealth](2,1) -- (3,1);
\draw [stealth-](2,2) -- (3,2);
\draw [-stealth](2,3) -- (3,3);

\draw [stealth-](1,1) -- (1,2);
\draw [-stealth](1,2) -- (1,3);

\draw [-stealth](2,1) -- (2,2);
\draw [stealth-](2,2) -- (2,3);

\draw [stealth-](3,1) -- (3,2);
\draw [-stealth](3,2) -- (3,3);

\end{tikzpicture}
\caption{\label{SFGA}The simple flow grid of matrix $D_3$.} 
\end{figure}

Combinatorial objects introduced in \cite{Striker1} that correspond to faces of $ASM_n$ are the {\it elementary flow grids}, which are directed graphs whose edge sets are the unions of the edge sets of {\it simple flow grids}.  Simple flow grids are modifications of digraphs corresponding to states of a 6-vertex model from statistical physics, which Kuperberg \cite{Kuperberg1} showed to correspond to alternating sign matrices.  (See Figure \ref{SFGA}.)  Striker showed how to calculate the dimension of a face from its corresponding elementary flow grid \cite[Theorem 4.5]{Striker1}.  

Brualdi and Dahl \cite{Brualdi4} showed that the $n \times n$ ASMs form a Hilbert basis for the cone generated by the $n\times n$ ASMs.  They also defined doubly directed graphs (see the next section) and made some progress towards determining the distance in the graph of $ASM_n$ from an ASM to the nearest permutation matrix.  M\'esz\'aros, Morales and Striker \cite{Striker4} showed that certain faces of $ASM_n$ are integrally equivalent to flow polytopes.

Recently, Borsik, Frank, Madarasi and Tak\'{a}cs \cite{borsik2025prefixboundedmatrices}, in response to a challenge from Edmonds \cite{Edmonds2}, showed among other things that the ASM polytope is an intersection of generalized polymatroids.

There do not appear to be papers on the alternating sign matrix polytope other than those we have cited.  In particular, there does not seem to be an analogue of the sequence of papers by Brualdi and Gibson \cite{BruGibI}, \cite{BruGibII}, \cite{BruGibIII} which derived many properties of the Birkhoff polytope $B_n$.  We aim to present new results, primarily using the combinatorial model of elementary flow grids.

In Section 2, we establish basic properties of elementary flow grids. 
We introduce Brualdi and Dahl's doubly directed graph, which contains much of the critical information contained in the elementary flow grid of a face.  

Section 3 gives some nontrivial restrictive conditions on the face lattices of faces of the ASM polytope.  These show that central symmetry of a face can be immediately read off of its doubly directed graph. For any two vertices of $ASM_n$ the smallest face containing them has all vertices of the doubly directed graph of even degree.

Section 4 shows that a basis for the subspace parallel to the affine span of the ASMs in a face consists of cycle matrices for cycles bounding the bounded regions of the 2-connected components of the doubly directed graph corresponding to the face, analogous to a result of MacLane \cite{MacLane1} for general planar graphs. 

Section 5 is devoted to establishing that a face of an ASM is affinely equivalent to the product of some of its faces if and only if its doubly directed graph is the edge-disjoint union of 2-connected components of the doubly directed graph of the faces.  This is more difficult to show than the analogous property of the Birkhoff polytope. 


In Section 6, we describe facets of a face of the ASM polytope in terms of the elementary flow grid.  An important consequence of this is the result that faces of the ASM polytope are 2-level polytopes.  Recent work on 2-level polytopes \cite{enum2level}, \cite{2levelcombsettings} will help us determine all the combinatorial types of low-dimensional faces of ASM polytopes.  $B_n$ is also a 2-level polytope \cite{2levelcombsettings}.

In Section 7 we show that the number of facets of a $d$-dimensional face of $ASM_n$ is at most $4(d-1)$.  We also show that the number of vertices of a $d$-dimensional face that is not affinely equivalent to the product of faces of dimension of at least one is at most $2^{d-1}+2$.  


Section 8 contains a proof that the Birkhoff polytope $B_3$ is not a face of $ASM_n$ for any $n$.  We hope that the obstruction encountered may find use in excluding larger classes of combinatorial types for faces of $ASM_n$.

Finally, Section 9 has a list of combinatorial types of 2, 3 and 4-dimensional faces of $ASM_n$. We compare our list with lists from the literature of combinatorial types of faces of Birkhoff polytopes and combinatorial types of 2-level polytopes.

We now summarize the similarities and differences between the polytopes $ASM_n$ and $B_n$ that we will discuss. The papers \cite {BruGibI, BruGibII, BruGibIII,Paffenholz1} contain many results that we have not yet addressed. 

\begin{itemize}
    \item The vertices of $ASM_n$ and $B_n$ are feasible flows in certain networks.  For $B_n$, the network is $K_{n,n}$ with the arcs oriented from one side of the bipartition to the other \cite{Balinski1}.  For $ASM_n$, the network is $C_n$, the complete flow grid with horizontal edges reversed. \cite[Section 4]{Behrend1}, \cite{borsik2025prefixboundedmatrices}.  
    \item $ASM_n$ (Proposition \ref{Two level}) and $B_n$ (\cite{2levelcombsettings}) are 2-level polytopes.
    \item A face of $ASM_n$ is a product of two or more of its faces if and only if its doubly directed graph is an edge-disjoint union of its 2-connected components, which are the doubly directed graphs of its faces (Theorem \ref{productmaintheorem} and Proposition \ref{productimpliesedgedisjoint}).  
    A face of $B_n$ is a product of two or more of its faces if and only if the connected components of its elementary graphs are the vertex-disjoint unions of the elementary graphs of the faces \cite[Proposition 2.5]{Paffenholz1}.  
    \item A $d$-dimensional face of $ASM_n$ has at most $4(d-1)$ facets (Theorem \ref{maxnumfacets}), while a $d$-dimensional face of $B_n$ has at most $3(d-1)$ facets \cite[Corollary 3.4]{BruGibIII}. Proposition \ref{2 connected product} shows that a face is not affinely equivalent to a product of faces of dimension at least one if and only if its doubly directed graph is 2-connected. Using this Proposition, we show that every such face of $ASM_n$ has at most $2^{d-1}+2$ vertices (Theorem \ref{2connectedub}) while a $d$-dimensional face of $B_n$ that is not such a product has at most $2^{d-1}+1$ vertices \cite[Theorem 5.1]{BruGibI}.
    \item For two given vertices of $ASM_n$, the smallest face containing both of them is centrally symmetric and has a doubly directed graph with all vertices of even degree (Theorem \ref{centsymevendegree}).  There are such faces of $ASM_n$ that are not affinely equivalent to cubes, starting with the standard octahedron in dimension 3.  For two given vertices of $B_n$, the smallest face containing them is combinatorially equivalent to a cube of some dimension \cite[Page 16]{Paffenholz1}.  
    \item $B_3$ is not a face of $ASM_n$ for any $n$ (Theorem \ref{nofaceB3}), and $ASM_3$ is not a face of $B_n$ for any $n$ \cite[Theorem 6.1]{Paffenholz1}.  
\end{itemize}

Much of this article is based on the first author's PhD thesis \cite{Dinkelman1}. The authors thank Jim Lawrence for helpful discussions.

\section {Elementary Flow Grids and Doubly Directed Graphs}

In this section, we first review the constructions from \cite{Striker1} and \cite{Brualdi4} of graphs that correspond to faces of the alternating sign matrix polytope as well as relevant results from those papers.  We also prove Lemmas \ref{degree4} and \ref{degree2vertices} which give some new properties of these graphs that will be useful in later sections. 

Striker \cite{Striker1} modified square ice configurations to form simple flow grids.  Striker considered a directed graph on $n^2+4n$ vertices with $n^2$ vertices $(i,j)$ for $i,j=1 \dots n$ which are internal and $4n$ which are boundary vertices. The boundary vertices are labelled $(i,0),(0,j),(i,n+1)$, and $(n+1,j)$ where $i,j=1, \dots ,n.$  The complete flow grid $C_n$ has directed edges ${((i,j),(i,j \pm 1))},{((i,j),(i \pm 1,j))}$ where $i,j=1, \dots ,n.$ 

See Figure \ref{C3}, where we follow the convention of drawing a pair of oppositely directed edges as a doubly directed red edge. The internal vertices of Figure \ref{C3} are labelled $(i,j)$ for $1 \leq i,j \leq 3$ with $(1,1)$ in the top left corner consistent with matrix labelling.

\begin{figure}[!h!] 
\centering 
\begin{tikzpicture}
\draw [stealth-](0.5,1) -- (1,1);
\draw [stealth-](0.5,2) -- (1,2);
\draw [stealth-](0.5,3) -- (1,3);

\draw [-stealth](3,1) -- (3.5,1);
\draw [-stealth](3,2) -- (3.5,2);
\draw [-stealth](3,3) -- (3.5,3);

\draw [-stealth](1,3) -- (1,3.5);
\draw [-stealth](2,3) -- (2,3.5);
\draw [-stealth](3,3) -- (3,3.5);

\draw [-stealth](1,1) -- (1,0.5);
\draw [-stealth](2,1) -- (2,0.5);
\draw [-stealth](3,1) -- (3,0.5);

\draw [red,thick,{stealth-stealth}](1,1) -- (2,1);
\draw [red,thick,{stealth-stealth}](1,2) -- (2,2);
\draw [red,thick,{stealth-stealth}](1,3) -- (2,3);

\draw [red,thick,{stealth-stealth}](2,1) -- (3,1);
\draw [red,thick,{stealth-stealth}](2,2) -- (3,2);
\draw [red,thick,{stealth-stealth}](2,3) -- (3,3);

\draw [red,thick,{stealth-stealth}](1,1) -- (1,2);
\draw [red,thick,{stealth-stealth}](1,2) -- (1,3);

\draw [red,thick,{stealth-stealth}](2,1) -- (2,2);
\draw [red,thick,{stealth-stealth}](2,2) -- (2,3);

\draw [red,thick,{stealth-stealth}](3,1) -- (3,2);
\draw [red,thick,{stealth-stealth}](3,2) -- (3,3);

\end{tikzpicture}
\caption{\label{C3} The complete flow grid $C_3$.}
\end{figure}

Striker \cite{Striker1} described a simple flow grid of order $n$ as a subgraph of $C_n$ consisting of all the vertices of $C_n$, with four edges incident to each internal vertex $(i,j)$: either four edges directed inward (corresponding to $a_{ij}= -1)$, four edges directed outward (corresponding to $a_{ij}=1)$, or two horizontal edges pointing in the same direction and two vertical edges pointing in the same direction (corresponding to $a_{ij}=0)$. 
See Figure \ref{simple flow grid model}.

\begin{proposition}\label{simpleflowgridsASMs}
\cite[Proposition 3.2]{Striker1} There exists an explicit bijection between simple flow grids of order $n$ and $n \times n$ alternating sign matrices ($ASMs$).  
\end{proposition}


\begin{figure}[!h!] 
\centering 
\begin{tikzpicture}

\draw[-stealth](1,0)--(1,1);
\draw[stealth-](1,-1)--(1,0);
\draw[stealth-](0,0)--(1,0);
\draw[-stealth](1,0)--(2,0);
\node at (1.4,0.4){1};
\node at (-.3,0){d)};

\draw[stealth-](4,0)--(4,1);
\draw[-stealth](4,-1)--(4,0);
\draw[-stealth](3,0)--(4,0);
\draw[stealth-](4,0)--(5,0);
\node at (4.4,0.4){-1};
\node at (2.7,0){e)};

\draw[-stealth](7,0)--(7,1);
\draw[-stealth](7,-1)--(7,0);
\draw[-stealth](6,0)--(7,0);
\draw[-stealth](7,0)--(8,0);
\node at (7.4,0.4){0};
\node at (5.7,0){f)};

\draw[stealth-](1,3)--(1,4);
\draw[stealth-](1,2)--(1,3);
\draw[stealth-](0,3)--(1,3);
\draw[stealth-](1,3)--(2,3);
\node at (1.4,3.4){0};
\node at (-.3,3){a)};

\draw[stealth-](4,3)--(4,4);
\draw[stealth-](4,2)--(4,3);
\draw[-stealth](3,3)--(4,3);
\draw[-stealth](4,3)--(5,3);
\node at (4.4,3.4){0};
\node at (2.7,3){b)};

\draw[-stealth](7,3)--(7,4);
\draw[-stealth](7,2)--(7,3);
\draw[stealth-](6,3)--(7,3);
\draw[stealth-](7,3)--(8,3);
\node at (7.4,3.4){0};
\node at (5.7,3){c)};

\end{tikzpicture}
\caption{\label{simple flow grid model}Vertex neighbourhoods in a simple flow grid.}
\end{figure}

Striker then defined an \emph{elementary flow grid} $G$ as a subgraph of $C_n$ where the set of edges of $G$ is the union of the edge sets of simple flow grids. Then any face $F$ of $ASM_n$ corresponds to an elementary flow grid $g(F)$ whose edge set is the union of the edge sets of simple flow grids which correspond to all the vertices of that face. If $F$ consists of a single vertex, i.e., $F=\{A\}$, then $g(\{A\})$ will be denoted simply as $g(A)$. A fundamental theorem \cite [Theorem 4.3]{Striker1} is that the face lattice of $ASM_n$ is isomorphic to the lattice of all $n \times n$ elementary flow grids ordered by graph inclusion. 

In many papers, for example \cite{Behrend1,Kuperberg1}, the horizontal edges of the simple and elementary flow grids are oriented oppositely to the description above. Simple flow grids in this alternate description can be thought of as feasible flows in a network, with a flow of size 1 on edges pointing left and edges pointing down. The six cases in Figure \ref{simple flow grid model}, once the horizontal edges are reversed, describe all the ways to have flow into a vertex equal to flow out of a vertex.  On the other hand, the approach that we take is more symmetric.

A prominent feature of elementary flow grids is that they are {\it planar}.  The definition of alternating sign matrices does not give any hint as to why planarity should feature in their analysis, but it turns out to be very helpful.  We can consider bounded regions in the standard drawing of the graph.


\begin{definition}\cite{Striker1} Two grid edges in an elementary flow grid $G$ form a {\it double directed edge} if they are opposite orientations of the same undirected edge.  An edge is said to be {\it fixed} in $G$ if the oppositely oriented edge is not in $G$.   
\end{definition}

We will draw a double directed edge of a flow grid as a single red edge with two arrowheads, while fixed edges remain directed and black. Then we can refer to regions bounded by red edges.

\begin{definition}
\cite{Striker1} A doubly directed region in an elementary flow grid $G$ is a collection of cells in $G$ completely bounded by double directed edges but containing no double directed edges in its interior. 
\end{definition}

 Striker \cite[Theorem 4.5]{Striker1} proved that the dimension of a face $F$ of $ASM_n$ is the number of doubly directed regions in the corresponding elementary flow grid $g(F)$.  In particular, the edges of $ASM_n$ are represented by elementary flow grids containing exactly one cycle of double directed edges.


In Figure \ref{simple and elementary flow grids}, the simple flow grids for two $3\times 3$ ASMs and the elementary flow grid of the union of the simple flow grids are shown. The doubly directed edges are colored red and the fixed edges are black.  Note that there is one set of double directed edges forming a cycle bounding a single doubly directed region, thus the two vertices are adjacent in $ASM_3$.

     \[
\left[ { \begin{array}{ccc}
    0 & 1 & 0 \\
    1 & -1 & 1 \\
    0 & 1 & 0 \\
  \end{array} } \right],
\left[ { \begin{array}{ccc}
    1 & 0 & 0 \\
    0 & 0 & 1 \\
    0 & 1 & 0 \\
  \end{array} } \right]
  \]
  \begin{tiny}
 \begin{figure}[h!] 
\centering 
\begin{tikzpicture}
\draw [stealth-](0.5,1) -- (1,1);
\draw [stealth-](0.5,2) -- (1,2);
\draw [stealth-](0.5,3) -- (1,3);

\draw [-stealth](3,1) -- (3.5,1);
\draw [-stealth](3,2) -- (3.5,2);
\draw [-stealth](3,3) -- (3.5,3);

\draw [-stealth](1,3) -- (1,3.5);
\draw [-stealth](2,3) -- (2,3.5);
\draw [-stealth](3,3) -- (3,3.5);

\draw [-stealth](1,1) -- (1,0.5);
\draw [-stealth](2,1) -- (2,0.5);
\draw [-stealth](3,1) -- (3,0.5);

\draw [stealth-](1,1) -- (2,1);
\draw [-stealth](1,2) -- (2,2);
\draw [stealth-](1,3) -- (2,3);

\draw [-stealth](2,1) -- (3,1);
\draw [stealth-](2,2) -- (3,2);
\draw [-stealth](2,3) -- (3,3);

\draw [stealth-](1,1) -- (1,2);
\draw [-stealth](1,2) -- (1,3);

\draw [-stealth](2,1) -- (2,2);
\draw [stealth-](2,2) -- (2,3);

\draw [stealth-](3,1) -- (3,2);
\draw [-stealth](3,2) -- (3,3);

\draw (4,2) node{$\cup$};

\draw [stealth-](4.5,1) -- (5,1);
\draw [stealth-](4.5,2) -- (5,2);
\draw [stealth-](4.5,3) -- (5,3);

\draw [-stealth](7,1) -- (7.5,1);
\draw [-stealth](7,2) -- (7.5,2);
\draw [-stealth](7,3) -- (7.5,3);

\draw [-stealth](5,3) -- (5,3.5);
\draw [-stealth](6,3) -- (6,3.5);
\draw [-stealth](7,3) -- (7,3.5);

\draw [-stealth](5,1) -- (5,0.5);
\draw [-stealth](6,1) -- (6,0.5);
\draw [-stealth](7,1) -- (7,0.5);

\draw [stealth-](5,1) -- (6,1);
\draw [stealth-](5,2) -- (6,2);
\draw [-stealth](5,3) -- (6,3);

\draw [-stealth](6,1) -- (7,1);
\draw [stealth-](6,2) -- (7,2);
\draw [-stealth](6,3) -- (7,3);

\draw [stealth-](5,1) -- (5,2);
\draw [stealth-](5,2) -- (5,3);

\draw [-stealth](6,1) -- (6,2);
\draw [-stealth](6,2) -- (6,3);

\draw [stealth-](7,1) -- (7,2);
\draw [-stealth](7,2) -- (7,3);

\draw (8,2) node{$=$};
\draw [stealth-](8.5,1) -- (9,1);
\draw [stealth-](8.5,2) -- (9,2);
\draw [stealth-](8.5,3) -- (9,3);

\draw [-stealth](11,1) -- (11.5,1);
\draw [-stealth](11,2) -- (11.5,2);
\draw [-stealth](11,3) -- (11.5,3);

\draw [-stealth](9,3) -- (9,3.5);
\draw [-stealth](10,3) -- (10,3.5);
\draw [-stealth](11,3) -- (11,3.5);

\draw [-stealth](9,1) -- (9,0.5);
\draw [-stealth](10,1) -- (10,0.5);
\draw [-stealth](11,1) -- (11,0.5);

\draw [stealth-](9,1) -- (10,1);
\draw [red,thick,{stealth-stealth}](9,2) -- (10,2);
\draw [red,thick,{stealth-stealth}](9,3) -- (10,3);

\draw [-stealth](10,1) -- (11,1);
\draw [stealth-](10,2) -- (11,2);
\draw [-stealth](10,3) -- (11,3);

\draw [stealth-](9,1) -- (9,2);
\draw [red,thick,{stealth-stealth}](9,2) -- (9,3);

\draw [-stealth](10,1) -- (10,2);
\draw [red,thick,{stealth-stealth}](10,2) -- (10,3);

\draw [stealth-](11,1) -- (11,2);
\draw [-stealth](11,2) -- (11,3);

\end{tikzpicture}
  










\caption{\label{simple and elementary flow grids}Example of simple and elementary flow grids.}

\end{figure}
\end{tiny}

Given a set of $n \times n$ $ASMs$ $\mathcal{X}$, let $F(\mathcal{X})$ be the smallest face of ASM$_n$ containing  $\mathcal{X}$. Brualdi and Dahl define the \emph{doubly directed graph} of $\mathcal{X}$ as an undirected graph, which we will denote $\mathcal{G}(F(\mathcal{X}))$, obtained from the elementary flow grid $g(F(\mathcal{X}))$ associated with $F(\mathcal{X})$ by replacing each pair of double directed edges by an undirected edge, and removing the remaining edges \cite[Page 38]{Brualdi4}.  That is, the orientations of the black edges of the elementary flow grid are ignored. We call the edges of $\mathcal{G}(F(\mathcal{X}))$ doubly directed edges.  

We now prove a proposition and two lemmas containing useful properties of doubly directed graphs.  Lemma \ref{degree4} will be used in Section 3 on central symmetry, while Lemma \ref{degree2vertices} will be used in Section 3 and in Section 4 on cycles.  

\begin{proposition} \label{propositiondeg3}
If $\mathcal{X}$ is a set of $ASMs$ and $\mathcal{G}(F(\mathcal{X}))$ has a
vertex of degree 3, then $|\mathcal{X}|>2$.
\end{proposition}
\begin{proof}  The union of any two of the vertex neighbourhoods of the simple flow grids depicted in Figure \ref{simple flow grid model}  has an even number of pairs of double directed edges.  
\end{proof}

The example in Figure \ref{vertex incident to 3 dd edges} is an elementary flow grid which is the union of three simple flow grids, producing a vertex of degree 3 in the $(2,2)$ position in $\mathcal{G}(F(\mathcal{X}))$:
     \[
g( F\left( \left\{ \left[ { \begin{array}{ccc}
    0 & 1 & 0 \\
    1 & -1 & 1 \\
    0 & 1 & 0 \\
  \end{array} } \right],
\left[ { \begin{array}{ccc}
    1 & 0 & 0 \\
    0 & 0 & 1 \\
    0 & 1 & 0 \\
  \end{array} } \right],
 \left[ { \begin{array}{ccc}
    0 & 1 & 0 \\
    0 & 0 & 1 \\
    1 & 0 & 0 \\
  \end{array} }  \right] \right\} \right) )
  \]

\begin{figure}[!h!]
\centering 
\begin{tikzpicture}

\node at (-1,2) {$=$};

\draw [stealth-](0.5,1) -- (1,1);
\draw [stealth-](0.5,2) -- (1,2);
\draw [stealth-](0.5,3) -- (1,3);

\draw [-stealth](3,1) -- (3.5,1);
\draw [-stealth](3,2) -- (3.5,2);
\draw [-stealth](3,3) -- (3.5,3);

\draw [-stealth](1,3) -- (1,3.5);
\draw [-stealth](2,3) -- (2,3.5);
\draw [-stealth](3,3) -- (3,3.5);

\draw [-stealth](1,1) -- (1,0.5);
\draw [-stealth](2,1) -- (2,0.5);
\draw [-stealth](3,1) -- (3,0.5);

\draw [red,thick,{stealth-stealth}](1,1) -- (2,1);
\draw [red,thick,{stealth-stealth}](1,2) -- (2,2);
\draw [red,thick,{stealth-stealth}](1,3) -- (2,3);

\draw [-stealth](2,1) -- (3,1);
\draw [stealth-](2,2) -- (3,2);
\draw [-stealth](2,3) -- (3,3);

\draw [red,thick,{stealth-stealth}](1,1) -- (1,2);
\draw [red,thick,{stealth-stealth}](1,2) -- (1,3);

\draw [red,thick,{stealth-stealth}](2,1) -- (2,2);
\draw [red,thick,{stealth-stealth}](2,2) -- (2,3);

\draw [stealth-](3,1) -- (3,2);
\draw [-stealth](3,2) -- (3,3);

\end{tikzpicture}
\caption{\label{vertex incident to 3 dd edges}A union of three simple flow grids which gives an elementary flow grid with a vertex adjacent to 3 doubly directed edges.}
\end{figure}

\begin{lemma} \label{degree4}
If a grid vertex $(i,j)$ of $ \mathcal{G} (F( \{ A,B \} ))$ is such that $a_{ij}=1$ and $b_{ij}=-1$ then $(i,j)$ is of degree 4 in $ \mathcal{G} (F( \{ A,B \} ))$. If $ \mathcal{G} (F( \{ A,B \} ))$ has a vertex $(i,j)$ of degree 4, then $a_{ij}b_{ij}=-1$ or $a_{ij}=b_{ij}=0$.
\end{lemma}
\begin{proof}
Suppose $(i,j)$ is such that $a_{ij}=1$ and $b_{ij}=-1$. Then there are 4 edges in $ \mathcal{G} (F( \{ A,B \} ))$ pointing out from $(i,j)$ and 4 edges pointing into $(i,j)$.  Hence there are 4 doubly directed grid edges incident to $(i,j)$.



   Assume that $(i,j)$ is of degree 4 in $ \mathcal{G} (F( \{ A,B \} ))$.  Suppose 
   that $a_{ij}=\pm 1$. The arrows of $g(B)$ at $(i,j)$ are pointing in the opposite direction from those of $g(A)$, so $b_{ij}=\mp 1$ by Figure \ref{simple flow grid model}, thus $a_{ij}b_{ij}=-1$.
   If $a_{ij}=0$ then the arrows of $g(B)$ at $(i,j)$ are pointing in the opposite direction from those of $g(A)$ and by Figure \ref{simple flow grid model} $b_{i,j}=0$.
\end{proof}

For matrices $A$ and $B$ below, $(2,2)$ is of degree 4 in $ \mathcal{G} (F( \{ A,B \} ))$.  See Figure \ref{2edgenot2vertex}.

     \[
A=\left[ { \begin{array}{ccc}
    1 & 0 & 0 \\
    0 & 0 & 1 \\
    0 & 1 & 0 \\
  \end{array} } \right]
B= \left[ { \begin{array}{ccc}
    0 & 1 & 0 \\
    1 & 0 & 0 \\
    0 & 0 & 1 \\
  \end{array} } \right]
  \]

\begin{lemma} \label{degree2vertices}
Suppose a grid vertex $(i,j)$ is of degree $2$ in $ \mathcal{G} (F(\mathcal{X} ))$.  Then the following hold. 
\begin{enumerate}
    \item Suppose the doubly directed edges are incident to $(i,j)$ on the same row or column. Then $a_{ij}=0$ for each $A$ in $\mathcal{X}$.  Also, one of the two fixed edges incident to $(i,j)$ enters $(i,j)$ and the other leaves $(i,j)$.
    \item Suppose the two edges incident to $(i,j)$ are at right angles. Then there exist matrices $A$ and $B$ in $\mathcal{X}$ with $a_{ij}=0$, $b_{ij}\neq 0$ and every matrix in $\mathcal{X}$ has $(i,j)$ entry equal to $a_{ij}$ or $b_{ij}$.  Also, the fixed edges incident to vertex $(i,j)$ both enter or both leave vertex $(i,j)$.
\end{enumerate}
  
\end{lemma}
\begin{proof}
(1) If vertex $(i,j)$ is incident to two horizontal doubly directed edges, then the vertex neighbourhoods of $(i,j)$ in the simple flow grids for the matrices in $\mathcal{X}$ must be of the form a and b or of the form c and f in Figure \ref{simple flow grid model}.  If vertex $(i,j)$ is incident to two vertical doubly directed edges, then the vertex neighbourhoods of $(i,j)$ in the simple flow grids for the matrices in $\mathcal{X}$ must be of the form a and c or of the form b and f in Figure \ref{simple flow grid model}.  




(2) If vertex $(i,j)$ of ${\mathcal G}(\mathcal{X})$ has exactly two edges incident to it at right angles, then the remaining two edges incident to vertex $(i,j)$ in $g(A)$ for $A\in \mathcal{X}$ must be fixed and also at right angles.  A look at Figure  \ref{simple flow grid model} shows that there must be one type $x$ in $\{a,b,c,f\}$ and one type $y$ in $\{d,e\}$ such that the neighborhood of $(i,j)$ in every $g(A)$ for $A \in \mathcal{X}$ is either of type $x$ or $y$.
\end{proof}  

\begin{figure}[h!] 
\centering 
\begin{tikzpicture}

\draw [-stealth](1,1) -- (1,2);
\draw [-stealth](1,2) -- (1,3);
\draw [-stealth](0,2) -- (1,2);
\draw [-stealth](1,2) -- (2,2);

\node at (2.5,2) {$\cup$};

\draw [stealth-](4,1) -- (4,2);
\draw [-stealth](4,2) -- (4,3);
\draw [stealth-](3,2) -- (4,2);
\draw [-stealth](4,2) -- (5,2);

\node at (5.5,2) {$=$};

\draw [-stealth](7,2) -- (7,3);
\draw [red,thick,{stealth-stealth}](7,1) -- (7,2);
\draw [red,thick,{stealth-stealth}](6,2) -- (7,2);
\draw [-stealth](7,2) -- (8,2);

\end{tikzpicture}
\caption{\label{pk}Example showing neighbourhoods of vertex $(i,j)$ of $g(A)$ and $g(B)$ with $(i,j)$ entries of $A$ and $B$ respectively 0 and 1. The elementary flow grid would produce $(i,j)$ in $\mathcal{G}(F (\{ A,B \} ))$ of degree 2.}
\end{figure}


An important theorem of Brualdi and Dahl \cite[Theorem 19]{Brualdi4} states, 
\begin{theorem} \label{theoremB}
Let $\mathcal{X}$ be a set of at least two $ n \times n$ ASMs.  Then each connected component of $\mathcal{G}(F(\mathcal{X}))$ is 2-edge-connected, or equivalently, every grid edge is contained in a cycle.  Therefore, $\mathcal{G}(F(\mathcal{X}))$ determines connected plane regions bounded by closed curves (cycles).
\end{theorem}

We will use the term 2-connected (without the word edge), to mean a graph with at least 3 vertices that is 2-vertex-connected, that is, removing less than 2 vertices doesn't disconnect the graph.

\begin{figure}[!h!] 
\centering 
\begin{tikzpicture}
\draw [stealth-](0.5,1) -- (1,1);
\draw [stealth-](0.5,2) -- (1,2);
\draw [stealth-](0.5,3) -- (1,3);
\draw [-stealth](3,1) -- (3.5,1);
\draw [-stealth](3,2) -- (3.5,2);
\draw [-stealth](3,3) -- (3.5,3);

\draw [-stealth](1,3) -- (1,3.5);
\draw [-stealth](2,3) -- (2,3.5);
\draw [-stealth](3,3) -- (3,3.5);

\draw [-stealth](1,1) -- (1,0.5);
\draw [-stealth](2,1) -- (2,0.5);
\draw [-stealth](3,1) -- (3,0.5);

\draw [stealth-](1,1) -- (2,1);
\draw [red,thick,{stealth-stealth}](1,2) -- (2,2);
\draw [red,thick,{stealth-stealth}](1,3) -- (2,3);

\draw [red,thick,{stealth-stealth}](2,1) -- (3,1);
\draw [red,thick,{stealth-stealth}](2,2) -- (3,2);
\draw [-stealth](3,2) -- (3,3);

\draw [stealth-](1,1) -- (1,2);
\draw [red,thick,{stealth-stealth}](1,2) -- (1,3);

\draw [red,thick,{stealth-stealth}](2,1) -- (2,2);
\draw [red,thick,{stealth-stealth}](2,2) -- (2,3);

\draw [red,thick,{stealth-stealth}](3,1) -- (3,2);
\draw [-stealth](2,3) -- (3,3);

\end{tikzpicture}
\caption{\label{2edgenot2vertex}$\mathcal{G}(F(\mathcal{X}))$ is 2-edge-connected but not 2-connected. It has two 2-connected components.}
\end{figure}

 While the proof in \cite{Brualdi4} did not construct a cycle of doubly-directed edges containing the edge, we give an alternative proof that does.

\begin{proof}
If every vertex in a graph is of even degree, then the graph can be partitioned into cycles \cite[Theorem 5.7]{Agnarsson1}. Vertices in $\mathcal{G}(F(\mathcal{X}))$ are of degree 1, 2, 3, or 4.  Brualdi and Dahl showed that every doubly directed edge has an incident doubly directed edge, so no vertex of $\mathcal{G}(F(\mathcal{X}))$ is of degree 1 \cite{Brualdi4}. Let $\{(i,j),(i,j+1)\}$ be an edge in $\mathcal{G}(\mathcal{X})$. Then there are $ASMs$ $A, B \in \mathcal{X} $ such that the edge $\{(i,j),(i,j+1)\}$ is in $\mathcal{G} (F(\{ A,B \} ))$. By Proposition \ref{propositiondeg3} $\mathcal{G}(F( \{ A,B \})) $ does not contain a vertex of degree 3. Thus all vertices of $\mathcal{G}(F( \{ A,B \})) $ are of even degree, hence the edge $\{(i,j),(i,j+1)\}$ lies in a cycle of $\mathcal{G}(F( \{ A,B \}))$. We also have  $\mathcal{G}(F( \{ A,B \})) \subseteq \mathcal{G} (F(\mathcal{X} ))$. Therefore $\{(i,j),(i,j+1)\}$ is contained in a cycle of $\mathcal{G}(F(\mathcal{X}))$. The same applies to edges $\{(i,j),(i+1,j)\}$ by an analogous argument.
\end{proof}  

We conclude this section by introducing some notation which shortens the statements of some results we have proved.

 \begin{definition}\label{signfixed}
     Suppose that $e$ is a fixed grid edge in an elementary flow grid $g(F)$ and $v$ is a vertex of the grid incident to $e$.  
 
 $sgn_F(v,e) = \begin{cases}
                    +1 & \mbox{ if }e \text{ is horizontal and pointing toward }v   \\
                    +1 & \mbox{ if } e \text{ is vertical and pointing away from } v \\
                    -1 & \mbox{ otherwise}
                 \end{cases} $
                 
If $X$ is a set of grid vertices, define  $$Sgn_F(X)=\sum_{v\in X}\sum_{e: e \text{ fixed and }v \in e} \mbox{sgn}_F(v,e).$$

\end{definition}

In order to demystify this notation, note that sgn$_F(v,e)$ can be interpreted as an entry of the node-arc incidence matrix of the digraph consisting of the fixed edges of $g(F)$, if we first reverse all of the horizontal fixed arcs.  

With this notation, we can restate Proposition \ref{simpleflowgridsASMs} and Lemma \ref{degree2vertices} more concisely: For every vertex $v$ of degree 0,2,or 4 in ${\mathcal G}(F)$ we have $Sgn_F(\{v\})=0$.

For any fixed edge $e=(v,w)$ of $g(F)$, sgn$_F(v,e)+\mbox{sgn}_F(w,e)=0$.  It follows that $Sgn_F(V({\mathcal G}(F)))=0$ 
for any face $F$ of $ASM_n$, where $V({\mathcal G}(F))$ denotes the set of vertices of  ${\mathcal G}(F)$.  This implies that for any face $F$ of $ASM_n$, there are just as many vertices $v$ with degree 3 in ${\mathcal G}(F)$ and Sgn$_F(\{v\})=-1$ as there are such $v$ with Sgn$_F(\{v\})=+1$.

\section{Estranged Vertices}

In this section we show that one can very quickly derive non-trivial properties of a face $F$ from its doubly directed graph.

\begin{definition}
Vertices $A$ and $B$ of a face $F$ are estranged in $F$ if there is no proper face of $F$ that contains $A$ and $B$.
\end{definition}
\begin{theorem}
\label{Allornone}
 If every vertex of $\mathcal{G} ( F)$ has even degree and $A$ is any vertex of the face $F$, then there is a unique vertex $B \in F$ such that $B$ is estranged from $A$ in $F$. If a face $F$ has a vertex of odd degree in ${\mathcal G}(F)$ then no vertex of $F$ has a vertex that is estranged from it in $F$.
\end{theorem}
\begin{proof}
Let $A \in F $ where every vertex of $\mathcal{G} (F )$ has even degree.  Let the matrix $B$ correspond to the directed graph produced by reversing the arrows of the simple flow grid of $A$ which correspond to the edges of  $\mathcal{G} (F )$.   For every subgraph in Figure \ref{simple flow grid model} we see that reversing any even number of the arrows in the subgraph produces another subgraph in the list. It follows that $g(B)$ is a simple flow grid and $B$ is an $ASM$.   $B$ is also in $F$.  If $C$ is  a vertex of $F$ other than $B$, then the simple flow grids for $C$ and $A$ would differ on a proper subset of the edges of $\mathcal{G}$, so $F$ would not be the smallest face containing $A$ and $C$.  Therefore, $B$ is unique.

Suppose that $A$ is a vertex of $F$ and that $(i,j)$ is a grid vertex of odd degree in the doubly directed graph of $F$.  Let $B$ be the simple flow grid obtained from the simple flow grid of $A$ by reversing all the edges that are in the doubly directed graph of $F$.  If one reverses any three of the edges of one of the subgraphs in Figure  \ref{simple flow grid model}, the resulting graph is not one of the allowed neighbourhoods of a vertex in a simple flow grid.  Therefore, there is no vertex $B$ for which $F$ is the smallest face containing $A$ and $B$.  
\end{proof}

 Note that the triangle contains no estranged vertices and its doubly directed graph contains vertices of degree 3. See Figure \ref{vertex incident to 3 dd edges}.
 
 We have obtained a very restrictive condition on the combinatorial types of faces of $ASM_n$. That is, if $F$ is a face of $ASM_n$, then either every vertex of $F$ has a unique vertex that is estranged from it in $F$ or no vertex of $F$ has a vertex that is estranged from it in $F$.  This implies the following corollary.

\begin{corollary} The only two-dimensional faces of $ASM_n$ are 4-gons and 3-gons.  
\end{corollary} 
\begin{proof}  If $F$ is a $k$-gon face of $ASM_n$ for $k>4$ and $A$ is a vertex of $F$, then there is more than one vertex in $F$ that is estranged from $A$, contradicting Theorem \ref{Allornone}.
\end{proof}


\begin{lemma} \label{centrally}
Given a grid vertex $(i,j)$ and a face $F$ for which all vertices of ${\mathcal G}(F)$ have even degree, $\frac{A+B}{2}$ is the same no matter which pair $\{A,B\}$ of estranged vertices of $F$ are chosen.  
\end{lemma}
\begin{proof}  

All vertices in ${\mathcal G}(F)$ are of even degree. Suppose $(i,j)$ is incident to exactly two edges of ${\mathcal G}(F)$ which are at right angles. If both fixed edges leave $(i,j)$, without loss of generality, $a_{ij} =1$ and $b_{ij}=0$ by Lemma \ref{degree2vertices}  so $\frac{a_{ij}+b_{ij}}{2}=\frac{1}{2}$.  If both fixed edges enter $(i,j)$, without loss of generality, $a_{ij} =-1$ and $b_{ij}=0$ by Lemma \ref{degree2vertices}  so $\frac{a_{ij}+b_{ij}}{2}=-\frac{1}{2}$.
At a vertex of degree 2 with incident edges in $\mathcal{G}(F)$ not at right angles, $a_{ij} =0$ and $b_{ij}=0$ by Lemma \ref{degree2vertices} so $\frac{a_{ij}+b_{ij}}{2}=0$.  
At a vertex of degree 4 in $\mathcal{G}(F)$ either $a_{ij} =0$ and $b_{ij}=0$ or $a_{ij} =1$ and $b_{ij}=-1$ by Lemma \ref{degree4} so in either case $\frac{a_{ij}+b_{ij}}{2}=0$.  
For $(i,j)$ which are not vertices of $\mathcal{G}(F)$, $a_{ij}$ is constant, so $\frac{a_{ij}+b_{ij}}{2}$ is the same for every pair of estranged vertices.
\end{proof}

For example, consider the face (an octahedron) corresponding to the elementary flow grid of Figure \ref{octahedron}. $A$ and $B$ are two estranged vertices. There are two other pairs of estranged vertices. All the segments connecting estranged pairs have the same midpoint.

  \[
A=\left[ { \begin{array}{cccccc}
    0 & 0 & 1 & 0 & 0 & 0 \\
    0 & 1 & 0 & 0 & 0 & 0 \\
    1 & 0 & 0 & 0 & 0 & 0 \\
    0 & 0 & 0 & 0 & 0 & 1 \\
    0 & 0 & 0 & 0 & 1 & 0 \\
    0 & 0 & 0 & 1 & 0 & 0 \\
  \end{array} } \right]  
B= \left[ { \begin{array}{cccccc}
    0 & 0 & 1 & 0 & 0 & 0 \\
    0 & 0 & 0 & 1 & 0 & 0 \\
    1 & 0 & -1 & 0 & 1 & 0 \\
    0 & 1 & 0 & -1 & 0 & 1 \\
    0 & 0 & 1 & 0 & 0 & 0 \\
    0 & 0 & 0 & 1 & 0 & 0 \\
  \end{array} } \right] 
  \]

    \[ 
\frac{A+B}{2} = \left[ { \begin{array}{cccccc}
    0 & 0 & 1 & 0 & 0 & 0 \\
    0 & \frac{1}{2} & 0 & \frac{1}{2} & 0 & 0 \\
    1 & 0 & -\frac{1}{2} & 0 & \frac{1}{2} & 0 \\
    0 & \frac{1}{2} & 0 & -\frac{1}{2} & 0 & 1 \\
    0 & 0 & \frac{1}{2} & 0 & \frac{1}{2} & 0 \\
    0 & 0 & 0 & 1 & 0 & 0 \\
  \end{array} } \right] 
  \]

\begin{figure}[h!]   
\centering
\begin{tikzpicture}
\draw [stealth-](0.5,1) -- (1,1);
\draw [stealth-](0.5,2) -- (1,2);
\draw [stealth-](0.5,3) -- (1,3);
\draw [stealth-](0.5,4) -- (1,4);
\draw [stealth-](0.5,5) -- (1,5);
\draw [stealth-](0.5,6) -- (1,6);

\draw [-stealth](6,1) -- (6.5,1);
\draw [-stealth](6,2) -- (6.5,2);
\draw [-stealth](6,3) -- (6.5,3);
\draw [-stealth](6,4) -- (6.5,4);
\draw [-stealth](6,5) -- (6.5,5);
\draw [-stealth](6,6) -- (6.5,6);

\draw [-stealth](1,6) -- (1,6.5);
\draw [-stealth](2,6) -- (2,6.5);
\draw [-stealth](3,6) -- (3,6.5);
\draw [-stealth](4,6) -- (4,6.5);
\draw [-stealth](5,6) -- (5,6.5);
\draw [-stealth](6,6) -- (6,6.5);

\draw [-stealth](1,1) -- (1,0.5);
\draw [-stealth](2,1) -- (2,0.5);
\draw [-stealth](3,1) -- (3,0.5);
\draw [-stealth](4,1) -- (4,0.5);
\draw [-stealth](5,1) -- (5,0.5);
\draw [-stealth](6,1) -- (6,0.5);

\draw [stealth-](1,1) -- (2,1);
\draw [stealth-](1,2) -- (2,2);
\draw [stealth-](1,3) -- (2,3);
\draw [-stealth](1,4) -- (2,4);
\draw [stealth-](1,5) -- (2,5);
\draw [stealth-](1,6) -- (2,6);

\draw [stealth-](2,1) -- (3,1);
\draw [stealth-](2,2) -- (3,2);
\draw [red,thick,{stealth-stealth}](2,3) -- (3,3);
\draw [-stealth](2,4) -- (3,4);
\draw [red,thick,{stealth-stealth}](2,5) -- (3,5);
\draw [stealth-](2,6) -- (3,6);

\draw [stealth-](3,1) -- (4,1);
\draw [red,thick,{stealth-stealth}](3,2) -- (4,2);
\draw [red,thick,{stealth-stealth}](3,3) -- (4,3);
\draw [red,thick,{stealth-stealth}](3,4) -- (4,4);
\draw [red,thick,{stealth-stealth}](3,5) -- (4,5);
\draw [-stealth](3,6) -- (4,6);

\draw [-stealth](4,1) -- (5,1);
\draw [red,thick,{stealth-stealth}](4,2) -- (5,2);
\draw [stealth-](4,3) -- (5,3);
\draw [red,thick,{stealth-stealth}](4,4) -- (5,4);
\draw [-stealth](4,5) -- (5,5);
\draw [-stealth](4,6) -- (5,6);

\draw [-stealth](5,1) -- (6,1);
\draw [-stealth](5,2) -- (6,2);
\draw [stealth-](5,3) -- (6,3);
\draw [-stealth](5,4) -- (6,4);
\draw [-stealth](5,5) -- (6,5);
\draw [-stealth](5,6) -- (6,6);

\draw [stealth-](1,1) -- (1,2);
\draw [stealth-](1,2) -- (1,3);
\draw [stealth-](1,3) -- (1,4);
\draw [-stealth](1,4) -- (1,5);
\draw [-stealth](1,5) -- (1,6);

\draw [stealth-](2,1) -- (2,2);
\draw [stealth-](2,2) -- (2,3);
\draw [red,thick,{stealth-stealth}](2,3) -- (2,4);
\draw [red,thick,{stealth-stealth}](2,4) -- (2,5);
\draw [-stealth](2,5) -- (2,6);

\draw [stealth-](3,1) -- (3,2);
\draw [red,thick,{stealth-stealth}](3,2) -- (3,3);
\draw [red,thick,{stealth-stealth}](3,3) -- (3,4);
\draw [stealth-](3,4) -- (3,5);
\draw [stealth-](3,5) -- (3,6);

\draw [-stealth](4,1) -- (4,2);
\draw [-stealth](4,2) -- (4,3);
\draw [red,thick,{stealth-stealth}](4,3) -- (4,4);
\draw [red,thick,{stealth-stealth}](4,4) -- (4,5);
\draw [-stealth](4,5) -- (4,6);

\draw [stealth-](5,1) -- (5,2);
\draw [red,thick,{stealth-stealth}](5,2) -- (5,3);
\draw [red,thick,{stealth-stealth}](5,3) -- (5,4);
\draw [-stealth](5,4) -- (5,5);
\draw [-stealth](5,5) -- (5,6);

\draw [stealth-](6,1) -- (6,2);
\draw [stealth-](6,2) -- (6,3);
\draw [-stealth](6,3) -- (6,4);
\draw [-stealth](6,4) -- (6,5);
\draw [-stealth](6,5) -- (6,6);

\end{tikzpicture}
\caption{\label{octahedron}Elementary flow grid produced by a pair of $6 \times 6$ estranged ASM vertices $A$ and $B$.}
\end{figure}
\vskip0.1in
Lemma \ref{centrally} immediately yields the following corollary.
   
\begin{corollary}
      If a face $F$ of $ASM_n$ has a pair of estranged vertices, then $F$ is centrally symmetric.  The center of $F$ is the matrix $\frac{A+B}{2}$ for some pair of estranged vertices $A$ and $B$.
\end{corollary}

We are now in a position to state the main result of this section.

\begin{theorem}\label{centsymevendegree}
    A face $F$ of $ASM_n$  is centrally symmetric if and only if its doubly directed graph has all vertices of even degree.
\end{theorem}

\begin{proof}
    Assume that every vertex of ${\mathcal G}(F)$ is of even degree. Then by Theorem \ref{Allornone}, there exists a pair of estranged vertices in $F$.  By Lemma \ref{centrally}, $F$ is centrally symmetric.
    Assume that there exists a vertex of ${\mathcal G}(F)$ of odd degree. Then by Theorem \ref{Allornone} there does not exist a pair of estranged vertices in $F$.  If $A$ and $B$ are vertices of $F$ then there is a facet $F'$ of $F$ containing $A$ and $B$. This implies that $\frac{A+B}{2}$ is on $F'$ which is on the boundary of $F$.  Thus $F$ is not centrally symmetric.
\end{proof}

Paffenholz \cite[Page 16]{Paffenholz1} points out that for any two vertices of $B_n$, the minimal face containing both of them is a cube of some dimension.  This is not true for vertices of $ASM_n$.  Centrally symmetric faces of $ASM_n$ with 2-connected doubly directed graphs include the octahedron in dimension 3, 
the bipyramid over the 3-cube in dimension 4 (Figure \ref{bipyramidover3cube}), and the faces corresponding to the three elementary flow grids of Figure \ref{bipyramidover4cube} in dimension 5.  

\begin{figure}[h!] 
\centering 
\begin{tikzpicture}[scale=0.65]

\draw [stealth-](0.5,1) -- (1,1);
\draw [stealth-](0.5,2) -- (1,2);
\draw [stealth-](0.5,3) -- (1,3);
\draw [stealth-](0.5,4) -- (1,4);

\draw [-stealth](4,1) -- (4.5,1);
\draw [-stealth](4,2) -- (4.5,2);
\draw [-stealth](4,3) -- (4.5,3);
\draw [-stealth](4,4) -- (4.5,4);

\draw [-stealth](1,4) -- (1,4.5);
\draw [-stealth](2,4) -- (2,4.5);
\draw [-stealth](3,4) -- (3,4.5);
\draw [-stealth](4,4) -- (4,4.5);

\draw [-stealth](1,1) -- (1,0.5);
\draw [-stealth](2,1) -- (2,0.5);
\draw [-stealth](3,1) -- (3,0.5);
\draw [-stealth](4,1) -- (4,0.5);

\draw [stealth-](1,1) -- (2,1);
\draw [red,thick,{stealth-stealth}](1,2) -- (2,2);
\draw [stealth-](1,4) -- (2,4);
\draw [red,thick,{stealth-stealth}](1,3) -- (2,3);

\draw [red,thick,{stealth-stealth}](2,1) -- (3,1);
\draw [red,thick,{stealth-stealth}](2,2) -- (3,2);
\draw [red,thick,{stealth-stealth}](2,3) -- (3,3);
\draw [red,thick,{stealth-stealth}](2,4) -- (3,4);

\draw [-stealth](3,1) -- (4,1);
\draw [red,thick,{stealth-stealth}](3,2) -- (4,2);
\draw [red,thick,{stealth-stealth}](3,3) -- (4,3);
\draw [-stealth](3,4) -- (4,4);

\draw [stealth-](1,1) -- (1,2);
\draw [red,thick,{stealth-stealth}](1,2) -- (1,3);
\draw [red,thick,{stealth-stealth}](2,3) -- (2,4);

\draw [red,thick,{stealth-stealth}](2,1) -- (2,2);
\draw [red,thick,{stealth-stealth}](2,2) -- (2,3);
\draw [-stealth](1,3) -- (1,4);

\draw [red,thick,{stealth-stealth}](3,1) -- (3,2);
\draw [red,thick,{stealth-stealth}](3,2) -- (3,3);
\draw [red,thick,{stealth-stealth}](3,3) -- (3,4);

\draw [stealth-](4,1) -- (4,2);
\draw [red,thick,{stealth-stealth}](4,2) -- (4,3);
\draw [-stealth](4,3) -- (4,4);

\draw [stealth-](6.5,1) -- (7,1);
\draw [stealth-](6.5,2) -- (7,2);
\draw [stealth-](6.5,3) -- (7,3);
\draw [stealth-](6.5,4) -- (7,4);
\draw [stealth-](6.5,5) -- (7,5);

\draw [-stealth](11,1) -- (11.5,1);
\draw [-stealth](11,2) -- (11.5,2);
\draw [-stealth](11,3) -- (11.5,3);
\draw [-stealth](11,4) -- (11.5,4);
\draw [-stealth](11,5) -- (11.5,5);

\draw [-stealth](7,5) -- (7,5.5);
\draw [-stealth](8,5) -- (8,5.5);
\draw [-stealth](9,5) -- (9,5.5);
\draw [-stealth](10,5) -- (10,5.5);
\draw [-stealth](11,5) -- (11,5.5);

\draw [-stealth](7,1) -- (7,0.5);
\draw [-stealth](8,1) -- (8,0.5);
\draw [-stealth](9,1) -- (9,0.5);
\draw [-stealth](10,1) -- (10,0.5);
\draw [-stealth](11,1) -- (11,0.5);

\draw [red,thick,{stealth-stealth}](7,1) -- (8,1);
\draw [stealth-](7,2) -- (8,2);
\draw [red,thick,{stealth-stealth}](7,3) -- (8,3);
\draw [stealth-](7,4) -- (8,4);
\draw [stealth-](7,5) -- (8,5);

\draw [red,thick,{stealth-stealth}](8,1) -- (9,1);
\draw [red,thick,{stealth-stealth}](8,2) -- (9,2);
\draw [red,thick,{stealth-stealth}](8,3) -- (9,3);
\draw [stealth-](8,4) -- (9,4);
\draw [red,thick,{stealth-stealth}](8,5) -- (9,5);

\draw [-stealth](9,1) -- (10,1);
\draw [red,thick,{stealth-stealth}](9,2) -- (10,2);
\draw [red,thick,{stealth-stealth}](9,3) -- (10,3);
\draw [red,thick,{stealth-stealth}](9,4) -- (10,4);
\draw [red,thick,{stealth-stealth}](9,5) -- (10,5);

\draw [-stealth](10,1) -- (11,1);
\draw [red,thick,{stealth-stealth}](10,2) -- (11,2);
\draw [-stealth](10,3) -- (11,3);
\draw [red,thick,{stealth-stealth}](10,4) -- (11,4);
\draw [-stealth](10,5) -- (11,5);

\draw [red,thick,{stealth-stealth}](7,1) -- (7,2);
\draw [red,thick,{stealth-stealth}](7,2) -- (7,3);
\draw [-stealth](7,3) -- (7,4);
\draw [-stealth](7,4) -- (7,5);

\draw [stealth-](8,1) -- (8,2);
\draw [red,thick,{stealth-stealth}](8,2) -- (8,3);
\draw [red,thick,{stealth-stealth}](8,3) -- (8,4);
\draw [red,thick,{stealth-stealth}](8,4) -- (8,5);

\draw [red,thick,{stealth-stealth}](9,1) -- (9,2);
\draw [red,thick,{stealth-stealth}](9,2) -- (9,3);
\draw [red,thick,{stealth-stealth}](9,3) -- (9,4);
\draw [-stealth](9,4) -- (9,5);

\draw [stealth-](10,1) -- (10,2);
\draw [stealth-](10,2) -- (10,3);
\draw [red,thick,{stealth-stealth}](10,3) -- (10,4);
\draw [red,thick,{stealth-stealth}](10,4) -- (10,5);

\draw [stealth-](11,1) -- (11,2);
\draw [red,thick,{stealth-stealth}](11,2) -- (11,3);
\draw [red,thick,{stealth-stealth}](11,3) -- (11,4);
\draw [-stealth](11,4) -- (11,5);

\draw [stealth-](13.5,1) -- (14,1);
\draw [stealth-](13.5,2) -- (14,2);
\draw [stealth-](13.5,3) -- (14,3);
\draw [stealth-](13.5,4) -- (14,4);
\draw [stealth-](13.5,5) -- (14,5);
\draw [stealth-](13.5,6) -- (14,6);

\draw [-stealth](19,1) -- (19.5,1);
\draw [-stealth](19,2) -- (19.5,2);
\draw [-stealth](19,3) -- (19.5,3);
\draw [-stealth](19,4) -- (19.5,4);
\draw [-stealth](19,5) -- (19.5,5);
\draw [-stealth](19,6) -- (19.5,6);

\draw [-stealth](14,6) -- (14,6.5);
\draw [-stealth](15,6) -- (15,6.5);
\draw [-stealth](16,6) -- (16,6.5);
\draw [-stealth](17,6) -- (17,6.5);
\draw [-stealth](18,6) -- (18,6.5);
\draw [-stealth](19,6) -- (19,6.5);

\draw [-stealth](14,1) -- (14,0.5);
\draw [-stealth](15,1) -- (15,0.5);
\draw [-stealth](16,1) -- (16,0.5);
\draw [-stealth](17,1) -- (17,0.5);
\draw [-stealth](18,1) -- (18,0.5);
\draw [-stealth](19,1) -- (19,0.5);

\draw [stealth-](14,1) -- (15,1);
\draw [stealth-](14,2) -- (15,2);
\draw [-stealth](14,3) -- (15,3);
\draw [stealth-](14,4) -- (15,4);
\draw [stealth-](14,5) -- (15,5);
\draw [stealth-](14,6) -- (15,6);

\draw [stealth-](15,1) -- (16,1);
\draw [red,thick,{stealth-stealth}](15,2) -- (16,2);
\draw [-stealth](15,3) -- (16,3);
\draw [stealth-](15,4) -- (16,4);
\draw [red,thick,{stealth-stealth}](15,5) -- (16,5);
\draw [stealth-](15,6) -- (16,6);

\draw [red,thick,{stealth-stealth}](16,1) -- (17,1);
\draw [red,thick,{stealth-stealth}](16,2) -- (17,2);
\draw [red,thick,{stealth-stealth}](16,3) -- (17,3);
\draw [red,thick,{stealth-stealth}](16,4) -- (17,4);
\draw [red,thick,{stealth-stealth}](16,5) -- (17,5);
\draw [red,thick,{stealth-stealth}](16,6) -- (17,6);

\draw [-stealth](17,1) -- (18,1);
\draw [red,thick,{stealth-stealth}](17,2) -- (18,2);
\draw [-stealth](17,3) -- (18,3);
\draw [stealth-](17,4) -- (18,4);
\draw [red,thick,{stealth-stealth}](17,5) -- (18,5);
\draw [-stealth](17,6) -- (18,6);

\draw [-stealth](18,1) -- (19,1);
\draw [-stealth](18,2) -- (19,2);
\draw [-stealth](18,3) -- (19,3);
\draw [stealth-](18,4) -- (19,4);
\draw [-stealth](18,5) -- (19,5);
\draw [-stealth](18,6) -- (19,6);

\draw [stealth-](14,1) -- (14,2);
\draw [stealth-](14,2) -- (14,3);
\draw [-stealth](14,3) -- (14,4);
\draw [-stealth](14,4) -- (14,5);
\draw [-stealth](14,5) -- (14,6);

\draw [stealth-](15,1) -- (15,2);
\draw [red,thick,{stealth-stealth}](15,2) -- (15,3);
\draw [red,thick,{stealth-stealth}](15,3) -- (15,4);
\draw [red,thick,{stealth-stealth}](15,4) -- (15,5);
\draw [-stealth](15,5) -- (15,6);

\draw [red,thick,{stealth-stealth}](16,1) -- (16,2);
\draw [red,thick,{stealth-stealth}](16,2) -- (16,3);
\draw [stealth-](16,3) -- (16,4);
\draw [red,thick,{stealth-stealth}](16,4) -- (16,5);
\draw [red,thick,{stealth-stealth}](16,5) -- (16,6);

\draw [red,thick,{stealth-stealth}](17,1) -- (17,2);
\draw [red,thick,{stealth-stealth}](17,2) -- (17,3);
\draw [-stealth](17,3) -- (17,4);
\draw [red,thick,{stealth-stealth}](17,4) -- (17,5);
\draw [red,thick,{stealth-stealth}](17,5) -- (17,6);

\draw [stealth-](18,1) -- (18,2);
\draw [red,thick,{stealth-stealth}](18,2) -- (18,3);
\draw [red,thick,{stealth-stealth}](18,3) -- (18,4);
\draw [red,thick,{stealth-stealth}](18,4) -- (18,5);
\draw [-stealth](18,5) -- (18,6);

\draw [stealth-](19,1) -- (19,2);
\draw [stealth-](19,2) -- (19,3);
\draw [stealth-](19,3) -- (19,4);
\draw [-stealth](19,4) -- (19,5);
\draw [-stealth](19,5) -- (19,6);

\end{tikzpicture}
\caption{\label{bipyramidover4cube}Elementary flow grids for centrally symmetric 5-dimensional faces.}
\end{figure}

\section{Vector Space of Cycle Matrices}  

If $F$ is a face of $ASM_n$ we will show that the edges of $F$ correspond to cycles contained in the doubly directed graph ${\mathcal G}(F)$.  These cycles in turn correspond to vectors in $\mathbb{R}^{n^2}$.  This leads us to ask if there is a characterization of the vector space spanned by the edges in terms of the doubly directed graph. 

\begin{definition}
A {\it corner} of a cycle subgraph $K$ of a grid graph is a vertex that is incident to one horizontal and one vertical edge of the cycle. 
\end{definition}

\begin {definition}
A cycle matrix is an $n\times n$ matrix whose non-zero entries are in positions corresponding to the corners of a cycle contained in the $n\times n$ grid graph and alternate between 1 and -1 as one goes around the cycle.
\end{definition}

\begin{proposition} \label{simple planar cycle matrix}
     If an $ASM$ $A$ is adjacent to an $ASM$ $B$ in the graph of the $ASM$ polytope then  $C=B-A$ is a cycle matrix. 
\end{proposition}
\begin{proof}
 From Striker's characterization (\cite{Striker1}, Theorem 4.5) of the dimension of a face as the number of doubly directed regions of the elementary flow grid, we know that the elementary flow grid of $A$ and $B$, $g(F( \{ A,B \}) )$, contains a single cycle of doubly directed edges. Thus $ \mathcal{G}(F( \{ A,B \}) )$ is a cycle.  For a grid vertex $(s,t)$ not in $ \mathcal{G}(F( \{ A,B \}) )$,  $a_{st}=b_{st}$. Therefore, the matrix $C=B-A$ has zeros in those entries which are not a part of the cycle. 

 Let $(i,j)$ be a vertex of $ \mathcal{G}(F( \{ A,B \}) )$ at which the edges of $ \mathcal{G}(F( \{ A,B \}) )$ incident to $(i,j)$ are both horizontal or both vertical.  From Lemma \ref{degree2vertices}, we have $a_{ij}=b_{ij}=0$. 

 Let $(i,j)$ be a vertex of $ \mathcal{G}(F( \{ A,B \}) )$ at which the edges of $ \mathcal{G}(F( \{ A,B \}) )$ incident to $(i,j)$ are 
 at right angles.  The neighbourhoods of $(i,j)$ in the simple flow grids of $A$ and $B$ disagree in orientation on two arcs at right angles and agree on the other two.  From Figure \ref{simple flow grid model} we see that one of $a_{ij}$ and $b_{ij}$ must be 0 and the other nonzero.  Thus, $c_{ij}=a_{ij}-b_{ij}$ is $\pm 1$.  

 Suppose that $c_{ij}=1$.  Then either $a_{ij}=1$ or $b_{ij}=-1$.  This implies that the arcs of the simple flow grid for $A$ corresponding to the edges of $ \mathcal{G}(F( \{ A,B \}) )$ incident to $(i,j)$ are directed away from $(i,j)$ and those in the simple flow grid for $B$ are directed toward $(i,j)$.  

 Suppose that there is a horizontal path in $ \mathcal{G}(F( \{ A,B \}) )$ from $(i,j)$ to $(i,j+k)$, another corner of $ \mathcal{G}(F( \{ A,B \}) )$. The arcs in the simple flow grids for $A$ and $B$ form horizontal directed paths between $(i,j)$ and $(i,j+k)$, as one can see from Lemma \ref{degree2vertices}.  Thus, the arcs for the simple flow grid for $A$ corresponding to the edges of $ \mathcal{G}(F( \{ A,B \}) )$ incident to $(i,j+k)$ are directed toward $(i,j+k)$ and those in the simple flow grid for $B$ are directed away from $(i,j+k)$.  This implies that $c_{i,j+k}=a_{i,j+k}-b_{i,j+k}=-1$.  

 Similarly, the corner vertex $(i+\ell,j)$ of $ \mathcal{G}(F( \{ A,B \}) )$ at the end of a vertical path from $(i,j)$ will also have $c_{i+\ell,j}=-1$.  If we start with $c_{ij}=-1$, then we obtain $c_{i,j+k}=c_{i+\ell,j}=1$.
\end{proof}

The converse of Proposition \ref{simple planar cycle matrix} is not true.  In $C_4$, we find the cycle $((1,1),(1,2),\newline (2,2),(2,3),(1,3),(1,4),(3,4),(3,1),(1,1))$. See Figure \ref{Not polyomino}, which shows an attempt to complete the elementary flow grid from the given doubly directed cycle. Fixing the edge between $(1,2)$ and $(1,3)$ violates Lemma \ref{degree2vertices} at one of the endpoints.  The cycle matrices for this cycle are not parallel to any edge of $ASM_4$.


   \begin{figure}[h!]
\centering 
\begin{tikzpicture}
\draw [stealth-](0.5,1) -- (1,1);
\draw [stealth-](0.5,2) -- (1,2);
\draw [stealth-](0.5,3) -- (1,3);
\draw [stealth-](0.5,4) -- (1,4);

\draw [-stealth](4,1) -- (4.5,1);
\draw [-stealth](4,2) -- (4.5,2);
\draw [-stealth](4,3) -- (4.5,3);
\draw [-stealth](4,4) -- (4.5,4);

\draw [-stealth](1,4) -- (1,4.5);
\draw [-stealth](2,4) -- (2,4.5);
\draw [-stealth](3,4) -- (3,4.5);
\draw [-stealth](4,4) -- (4,4.5);

\draw [-stealth](1,1) -- (1,0.5);
\draw [-stealth](2,1) -- (2,0.5);
\draw [-stealth](3,1) -- (3,0.5);
\draw [-stealth](4,1) -- (4,0.5);

\draw [stealth-](1,1) -- (2,1);
\draw [red,thick,{stealth-stealth}](1,2) -- (2,2);
\draw [-stealth](1,3) -- (2,3);
\draw [red,thick,{stealth-stealth}](1,4) -- (2,4);

\draw [stealth-](2,1) -- (3,1);
\draw [red,thick,{stealth-stealth}](2,2) -- (3,2);
\draw [red,thick,{stealth-stealth}](2,3) -- (3,3);
\draw [-stealth](2,4) -- (3,4);

\draw [-stealth](3,1) -- (4,1);
\draw [red,thick,{stealth-stealth}](3,2) -- (4,2);
\draw [-stealth](3,3) -- (4,3);
\draw [red,thick,{stealth-stealth}](3,4) -- (4,4);

\draw [stealth-](1,1) -- (1,2);
\draw [red,thick,{stealth-stealth}](1,2) -- (1,3);
\draw [red,thick,{stealth-stealth}](1,3) -- (1,4);

\draw [stealth-](2,1) -- (2,2);
\draw [-stealth](2,2) -- (2,3);
\draw [red,thick,{stealth-stealth}](2,3) -- (2,4);

\draw [-stealth](3,1) -- (3,2);
\draw [-stealth](3,2) -- (3,3);
\draw [red,thick,{stealth-stealth}](3,3) -- (3,4);

\draw [stealth-](4,1) -- (4,2);
\draw [red,thick,{stealth-stealth}](4,2) -- (4,3);
\draw [red,thick,{stealth-stealth}](4,3) -- (4,4);

\end{tikzpicture}
\caption{\label{Not polyomino}Not an elementary flow grid.}
\end{figure}

In order to generalize this example, we introduce some useful notation which will be used often in Sections 6 - 8.  For the rest of this section, we assume that $F$ is a face of $ASM_n$.   

\begin{definition}\label{eardefinition}
  An ear of ${\mathcal G}(F)$ is a path between vertices of degree 3 or 4 in which every interior vertex has degree 2, or a cycle containing at most one vertex of degree 3 or 4.   
\end{definition}

Note that every edge of ${\mathcal G}(F)$ belongs to a unique ear of ${\mathcal G}(F)$.

 The elementary flow grid for a face $F'$ of $F$ can be obtained by replacing some doubly directed edges of $g(F)$ by single directed edges (fixing doubly directed edges.)  
 The following lemma shows that the fixed edges that are in a given ear cannot be a proper subset of the edges in the ear.

\begin{lemma}
    Suppose $F' \subseteq F$ are faces of $ASM_n$ and $e$ is a grid edge in an ear $S$ of ${\mathcal G}(F)$.  If $e$ is not in ${\mathcal G}(F')$, then the orientations in $g(F')$ of all the edges of $S$ are determined by the orientation of $e$. 
\end{lemma}


\begin{proof} Suppose $e$ is in $S\backslash {\mathcal G}(F')$ and let $e'$ be an edge of $S$ adjacent to $e$ such that $e\cap e'=(i,j)$. Because $S$ is a path or a cycle, two edges of $g(F)$ incident to $(i,j)$ are fixed.  These two edges together with $e$ are fixed in $g(F').$  One sees from Figure \ref{simple flow grid model} that the orientations of three arcs incident to a vertex $(i,j)$ determine the orientation of the fourth. 
\end{proof}

\begin{lemma}\label{twowaystofix}
    If $F$ is a face of $ASM_n$ and $S$ is an ear of ${\mathcal G}(F)$, then there are two orientations of the edges in $S$ that may appear in $g(F')$ for a face $F'$ of $F$. These two orientations are opposite to each other.
\end{lemma}

\begin{proof} 
As in the previous proof, suppose that $e$ and $e'$ are adjacent edges of $S$ incident to $(i,j)$.  The two other edges of $g(F)$ incident to $(i,j)$ are fixed.  From Figure \ref{simple flow grid model} we see that the two ways to orient $e$ and $e'$ consistent with the fixed arcs incident to $(i,j)$ are opposite to each other.
\end{proof}

\begin{figure}[!h!]
\centering 
\begin{tikzpicture}

\draw [stealth-](0.5,1) -- (1,1);
\draw [stealth-](0.5,2) -- (1,2);
\draw [stealth-](0.5,3) -- (1,3);

\draw [-stealth](3,1) -- (3.5,1);
\draw [-stealth](3,2) -- (3.5,2);
\draw [-stealth](3,3) -- (3.5,3);

\draw [-stealth](1,3) -- (1,3.5);
\draw [-stealth](2,3) -- (2,3.5);
\draw [-stealth](3,3) -- (3,3.5);

\draw [-stealth](1,1) -- (1,0.5);
\draw [-stealth](2,1) -- (2,0.5);
\draw [-stealth](3,1) -- (3,0.5);

\draw [red,thick,{stealth-stealth}](1,1) -- (2,1);
\draw [red,thick,{stealth-stealth}](1,2) -- (2,2);
\draw [red,thick,{stealth-stealth}](1,3) -- (2,3);

\draw [-stealth](2,1) -- (3,1);
\draw [stealth-](2,2) -- (3,2);
\draw [-stealth](2,3) -- (3,3);

\draw [red,thick,{stealth-stealth}](1,1) -- (1,2);
\draw [red,thick,{stealth-stealth}](1,2) -- (1,3);

\draw [red,thick,{stealth-stealth}](2,1) -- (2,2);
\draw [red,thick,{stealth-stealth}](2,2) -- (2,3);

\draw [stealth-](3,1) -- (3,2);
\draw [-stealth](3,2) -- (3,3);

\draw [stealth-](5,1) -- (5.5,1);
\draw [stealth-](5,2) -- (5.5,2);
\draw [stealth-](5,3) -- (5.5,3);

\draw [-stealth](7.5,1) -- (8,1);
\draw [-stealth](7.5,2) -- (8,2);
\draw [-stealth](7.5,3) -- (8,3);

\draw [-stealth](5.5,3) -- (5.5,3.5);
\draw [-stealth](6.5,3) -- (6.5,3.5);
\draw [-stealth](7.5,3) -- (7.5,3.5);

\draw [-stealth](5.5,1) -- (5.5,0.5);
\draw [-stealth](6.5,1) -- (6.5,0.5);
\draw [-stealth](7.5,1) -- (7.5,0.5);

\draw [red,thick,{stealth-stealth}](5.5,1) -- (6.5,1);
\draw [stealth-](5.5,2) -- (6.5,2);
\draw [red,thick,{stealth-stealth}](5.5,3) -- (6.5,3);

\draw [-stealth](6.5,1) -- (7.5,1);
\draw [stealth-](6.5,2) -- (7.5,2);
\draw [-stealth](6.5,3) -- (7.5,3);

\draw [red,thick,{stealth-stealth}](5.5,1) -- (5.5,2);
\draw [red,thick,{stealth-stealth}](5.5,2) -- (5.5,3);

\draw [red,thick,{stealth-stealth}](6.5,1) -- (6.5,2);
\draw [red,thick,{stealth-stealth}](6.5,2) -- (6.5,3);

\draw [stealth-](7.5,1) -- (7.5,2);
\draw [-stealth](7.5,2) -- (7.5,3);

\draw [stealth-](9.5,1) -- (10,1);
\draw [stealth-](9.5,2) -- (10,2);
\draw [stealth-](9.5,3) -- (10,3);

\draw [-stealth](12,1) -- (12.5,1);
\draw [-stealth](12,2) -- (12.5,2);
\draw [-stealth](12,3) -- (12.5,3);
\draw [-stealth](10,3) -- (10,3.5);
\draw [-stealth](11,3) -- (11,3.5);
\draw [-stealth](12,3) -- (12,3.5);

\draw [-stealth](10,1) -- (10,0.5);
\draw [-stealth](11,1) -- (11,0.5);
\draw [-stealth](12,1) -- (12,0.5);

\draw [stealth-](10,1) -- (11,1);
\draw [-stealth](10,2) -- (11,2);
\draw [stealth-](10,3) -- (11,3);

\draw [-stealth](11,1) -- (12,1);
\draw [stealth-](11,2) -- (12,2);
\draw [-stealth](11,3) -- (12,3);

\draw [stealth-](10,1) -- (10,2);
\draw [-stealth](10,2) -- (10,3);

\draw [-stealth](11,1) -- (11,2);
\draw [stealth-](11,2) -- (11,3);

\draw [stealth-](12,1) -- (12,2);
\draw [-stealth](12,2) -- (12,3);

\end{tikzpicture}
\caption{\label{Two ways to fix a ear} An elementary flow grid from a triangular face and two of its subfaces.}
\end{figure}

  In Figure \ref{Two ways to fix a ear} we revisit the face from Figure \ref{vertex incident to 3 dd edges}.  The doubly directed graph has three ears, each with endpoints $(2,1)$ and $(2,2)$: 
\begin{itemize}
    \item
    $(2,1), (1,1), (1,2), (2,2)$
    \item
    $(2,1), (2,2)$
    \item
    $(2,1), (3,1), (3,2), (2,2)$.
\end{itemize} 
Fixing the horizontal ear between $(2,1)$ and $(2,2)$ from right to left yields a 1-dimensional face. Fixing it from left to right forces the other ears to also be fixed, yielding a simple flow grid corresponding to a vertex. Therefore the face is a triangle.

If we fix an ear and reverse its horizontal edges, we get a directed path.  This is an advantage of the alternate description of elementary flow grids with the horizontal edges reversed as described in Section 2.

A useful property of the sgn$_F$ function is established next.

\begin{lemma}\label{signs of ear ends}  Suppose $v$ and $w$ are the endpoints of an ear $S$ of ${\mathcal G}(F)$ and $e_v$ and $e_w$ are the ear edges incident to $v$ and $w$.  Suppose that we fix $S$ in one of the two ways we can fix it.  Then  sgn$_F(v,e_v)=-\mbox{sgn}_F(w,e_w)$. 
\end{lemma}
\begin{proof}
    If $v$ and $w$ are adjacent vertices of $S$, the conclusion clearly holds.  It also holds that sgn$_F(v,e)=-\mbox{sgn}_F(v,f)$ if $e$ and $f$ are distinct edges of $S$ incident to an internal vertex $v$.   
\end{proof}

The edge between $v=(1,2)$ and $w=(1,3)$ in $C_4$ forms what we call a \emph{unisigned chord} with respect to the doubly directed cycle of Figure \ref{Not polyomino}. The following proposition generalizes the example.

\begin{definition}
    Let $F$ be a face of $ASM_n$.
    If $C$ is a cycle of ${\mathcal G}(F)$, $S$ is an ear of ${\mathcal G}(F)$ with endpoints $v$ and $w$ in $C$, but no edges of $C$, and $v,w$ have degree $3$ in ${\mathcal G}(F)$ with $\mbox{Sgn}_F(\{v\})=\mbox{Sgn}_F(\{w\})$, then $S$ is called a unisigned chord of $C$.
\end{definition}

\begin{proposition}
Let $F$ be a face of $ASM_n$. 
    If a cycle $C$ of ${\mathcal G}(F)$ has a unisigned chord, then the cycle matrix for $C$ is not parallel to any edge of $ASM_n$. 
\end{proposition}
\begin{proof}
    If the cycle matrix for $C$ is parallel to an edge $F'$ of $F$ then ${\mathcal G}(F')$ would be $C$.  In that case, the ear $S$ that is a unisigned chord would have to be fixed. After fixing $S$, Lemma \ref{signs of ear ends} shows that sgn$_F(v,e_v)=-\mbox{sgn}_F(w,e_w)$ where $e_v$ and $e_w$ are the edges of $S$ incident to $v$ and $w$.  Because $\mbox{Sgn}_F(\{v\})=\mbox{Sgn}_F(\{w\})$ then either $\mbox{sgn}_F(v,e_v)=\mbox{Sgn}_F(\{v\})$ or $\mbox{sgn}_F(w,e_w)=\mbox{Sgn}_F(\{w\})$.  If $\mbox{sgn}_F(v,e_v)=\mbox{Sgn}_F(\{v\})$ then we have $\mbox{sgn}_F(v,e_v)+\mbox{Sgn}_F(\{v\}) \neq 0.$  This implies that the edges of $C$ incident to $v$ cannot be doubly directed, meaning that ${\mathcal G}(F')$ cannot be $C$.  A similar conclusion holds if $\mbox{sgn}_F(w,e_w)=\mbox{Sgn}_F(\{w\})$. 
\end{proof}

\begin{definition}
A cycle $K$ in the doubly directed graph $\mathcal{G}(F)$ of a face $F$ is a basic cycle for $F$ if it bounds a region of the 2-connected component of $\mathcal{G}(F)$ containing $K$.
\end{definition}

\begin{definition}
    For each doubly directed region $R$ of $g(F)$, let $k(R)$ be the basic cycle that bounds $R$ and is incident to $R$.   
\end{definition}

The function $k$ is a bijection from the set of doubly directed regions of $F$ to the set of basic cycles for $F$.

As noted by Whitney \cite[Theorem 12]{Whitney1}, the decomposition of a graph into 2-connected components is unique, and every cycle of a graph is in a unique 2-connected component for the graph.  We want to show that a cycle matrix corresponding to a cycle of $\mathcal{G}(F)$ is in the span of the basic cycle matrices of $\mathcal{G}(F)$.

\begin{definition}
    A cycle matrix $C$ corresponding to a cycle $K$ in a grid graph is said to be {\it clockwise} with respect to $K$ if a maximal horizontal segment of the cycle always goes from a $+1$ entry of $C$ to a $-1$ entry if one goes clockwise around the cycle.  
\end{definition}

It follows from the definition that for a clockwise cycle $C$ a maximal vertical segment always goes from a $-1$ of $C$ to a $+1$ of $C$ as one goes clockwise around the cycle.  

\begin{proposition}\label{basiccyclesspan}
    Suppose $C$ is a clockwise cycle matrix corresponding to a cycle $K$ of ${\mathcal G}(F)$.  Let $G$ be the 2-connected component of ${\mathcal G}(F)$ containing $K$.  Let $C_1,C_2,\ldots,C_t$ be the clockwise basic cycle matrices for the regions of $G$ that are in the inside of $K$. Then $C_1+C_2+\cdots+C_t=C.$
\end{proposition}

A similar result is discussed in \cite[Theorem 5.3]{MacLane1}, but there addition is done modulo 2.  Our addition involves matrices with real entries.

\begin{proof}  In Figures \ref{Four or three cycles} and \ref{Two cycles sharing a corner} the regions $R_1,R_2,\ldots$ containing a vertex $(i,j)$ correspond to basic cycle matrices $C_1,C_2,\ldots$.  The clockwise property guarantees that the sum of the $(i,j)$ entries of the cycle matrices is zero.  The vertex $(i,j)$ is not in $K$, so the entry $(i,j)$ of $C$ is also zero.


\begin{figure}[!h!]

\centering 
\begin{tikzpicture}[scale=1.5]

\node at (0.5,2){a)};

\draw [red,thick,{stealth-stealth}](1,2) -- (2,2);
\draw [red,thick,{stealth-stealth}](2,1) -- (2,2);
\draw [red,thick,{stealth-stealth}](2,2) -- (2,3);
\draw [red,thick,{stealth-stealth}](2,2) -- (3,2);
\node at (1.8,2.2){1};
\node at (1.8,1.8){-1};
\node at (2.2,2.2){-1};
\node at (2.2,1.8){1};
\node at (1.5,2.5){$R_1$};
\node at (1.5,1.5){$R_4$};
\node at (2.5,2.5){$R_2$};
\node at (2.5,1.5){$R_3$};
\node at (4.5,2){b)};

\draw [stealth-](6,1) -- (6,2);
\draw [red,thick,{stealth-stealth}](6,2) -- (6,3);
\draw [red,thick,{stealth-stealth}](5,2) -- (6,2);
\draw [red,thick,{stealth-stealth}](6,2) -- (7,2);

\node at (6.2,2.2){-1};
\node at (5.8,1.8){0};
\node at (5.8,2.2){1};
\node at (5.5,2.5){$R_1$};
\node at (6.5,2.5){$R_2$};
\node at (6.5,1.5){$R_3$};

\end{tikzpicture}
\caption{\label{Four or three cycles} Four or three  cycles sharing a corner not in $K$.}

\end{figure}


\begin{figure}[!h!]

\centering 
\begin{tikzpicture}[scale=1.5]

\node at (0.5,2){a)};

\draw [red,thick,{stealth-stealth}](1,2) -- (2,2);
\draw [stealth-](2,1) -- (2,2);
\draw [red,thick,{stealth-stealth}](2,2) -- (2,3);
\draw [-stealth](2,2) -- (3,2);

\node at (1.8,2.2){1};
\node at (2.2,2.2){-1};
\node at (1.5,2.5){$R_1$};
\node at (2.5,1.5){$R_2$};
\node at (4.5,2){b)};
\draw [red,thick,{stealth-stealth}](5,2) -- (6,2);
\draw [stealth-](6,1) -- (6,2);
\draw [red,thick,{stealth-stealth}](6,2) -- (7,2);
\draw [stealth-](6,2) -- (6,3);

\node at (5.8,2.2){0};
\node at (6.2,1.8){0};
\node at (5.5,2.5){$R_1$};
\node at (6.5,1.5){$R_2$};

\end{tikzpicture}
\caption{\label{Two cycles sharing a corner} Two cycles sharing a vertex not in $K$.}

\end{figure}

In Figure \ref{Three or two cycles}a we assume that the edges bounding $R_4$ in Figure \ref{Four or three cycles}a are part of cycle $K$.  Thus region $R_4$ is removed. The sum of the $(i,j)$ entries of the matrices $C_1,C_2,C_3$ corresponding to the regions $R_1,R_2,R_3$ is the $(i,j)$ entry of $K$.  

In Figure \ref{Three or two cycles}b, the horizontal edges bounding $R_3$ in Figure \ref{Four or three cycles}b are part of cycle $K$ and region $R_3$ is removed.  Now the sum of the $(i,j)$ entries of $C_1,C_2$ is zero, which is also the $(i,j)$ entry of $C$. 

It could also happen that the edges bounding $R_1$ of Figure \ref{Four or three cycles} are part of $K$, and there is no region $R_1$.  In this case $K$ bounds the union of $R_2$ and $R_3$ and the sum of the $(i,j)$ entries of $C_2$ and $C_3$ equals the $(i,j)$ entry of $C$.

Finally, we consider the case when only two edges of $G$ are incident to $(i,j)$.  Then $K$ shares these edges with exactly one $R_i$ and $C$ must have the same $(i,j)$ entry as the corresponding $C_i$.
\end{proof}


\begin{figure}[!h!]

\centering 
\begin{tikzpicture}[scale=1.5]

\node at (0.5,2){a)};

\draw [red,thick,{stealth-stealth}](1,2) -- (2,2);
\draw [red,thick,{stealth-stealth}](2,1) -- (2,2);
\draw [red,thick,{stealth-stealth}](2,2) -- (2,3);
\draw [red,thick,{stealth-stealth}](2,2) -- (3,2);
\node at (1.8,2.2){1};
\node at (2.2,2.2){-1};
\node at (2.2,1.8){1};
\node at (1.5,2.5){$R_1$};
\node at (2.5,2.5){$R_2$};
\node at (2.5,1.5){$R_3$};
\node at (4.5,2){b)};

\draw [stealth-](6,1) -- (6,2);
\draw [red,thick,{stealth-stealth}](6,2) -- (6,3);
\draw [red,thick,{stealth-stealth}](5,2) -- (6,2);
\draw [red,thick,{stealth-stealth}](6,2) -- (7,2);

\node at (6.2,2.2){-1};
\node at (5.8,2.2){1};
\node at (5.5,2.5){$R_1$};
\node at (6.5,2.5){$R_2$};

\end{tikzpicture}
\caption{\label{Three or two cycles} Three or two  cycles sharing a corner in $K$.}

\end{figure}

\begin{proposition}\label{basiccyclesindependent}
For any face $F$, the collection ${\mathcal K}$ of cycle matrices corresponding to basic cycles for $F$ is linearly independent.
\end{proposition}
\begin{proof}
For every basic cycle $K$, let $(i_K,j_K)$ be the lexicographically largest vertex in $K$.  $(i_K,j_K)$ is clearly a corner vertex of $K$, so the $(i_K,j_K)$ entry of the corresponding cycle matrix is nonzero.  The vertices $(i_K,j_K)$ for $K\in {\mathcal K}$ can also be lexicographically ordered.  If $(i_K,j_K)<(i_L,j_L)$ for $K,L\in {\mathcal K},$ then $(i_L,j_L)\notin K.$  Suppose $K\neq L \in {\mathcal K},$ and $(i_K,j_K)= (i_L,j_L)$. Then $K$ and $L$ must contain the grid edge $((i_K,j_K-1),(i_K,j_K))$, so they are in the same 2-connected component $G_1$ of ${\mathcal G}(F)$. Because $(i_K,j_K)$ is the lexicographically largest corner of $K$ and $L$, the regions they bound overlap, contradicting the requirement that they bound disjoint regions of $G_1$.  

Consider a $|{\mathcal K}|\times n^2$ matrix $M$ where the columns correspond to vertices of the $n\times n$ grid graph ordered lexicographically and the rows are the cycle matrices for the cycles in ${\mathcal K}$, ordered such that the $(i_K,j_K)$ pairs increase lexicographically.  The previous discussion shows that this matrix contains a $|{\mathcal K}|\times |{\mathcal K}|$ triangular submatrix with nonzero diagonal entries.
\end{proof}

Propositions \ref{basiccyclesspan} and \ref{basiccyclesindependent} together give us the main result of this section:

\begin{theorem}
    For any face $F$, the collection ${\mathcal K}$ of cycle matrices corresponding to basic cycles for $F$ is a basis for the subspace parallel to $F$. 
\end{theorem}

With the last theorem of this section, we want to show that if $F$ is a face with each vertex of ${\mathcal G}(F)$ of even degree, that is, a centrally symmetric face by Theorem \ref{centsymevendegree}, then every cycle in ${\mathcal G}(F)$ corresponds to an edge of $F$.  In other words, a cycle such as that in Figure \ref{Not polyomino} will not occur.  

We first need to establish that one can certify that a subdigraph of the complete flow grid $C_n$ is an elementary flow grid by only looking locally at each vertex instead of explicitly writing it as a union of simple flow grids, in the special case in which each vertex is incident to an even number of doubly directed edges.

\begin{lemma}\label{even elementary flow grid}
Suppose that a subdigraph $G$ of the complete flow grid $C_n$ has all vertices of its doubly directed graph of even degree. Suppose that at every grid vertex not incident to any doubly directed edges one of the configurations of Figure \ref{simple flow grid model} occurs, and that at every grid vertex incident to two pairs of doubly directed edges, one of the configurations of Lemma \ref{degree2vertices} appears.  Then $G$ is an elementary flow grid.
\end{lemma}
\begin{proof}
 The proof is by induction on the number of doubly directed regions.  Every edge of the doubly directed graph of $G$ is contained in a cycle, so if the doubly directed graph of $G$ has no doubly directed regions then it has no edges.  
 
 Then $G$ is a simple flow grid by the assumption that at every vertex not incident to a doubly directed edge, one of the configurations of Figure \ref{simple flow grid model} occurs.  
 
 Now assume that the lemma is true for all $G{'}$ that have $k$ or fewer doubly directed regions, for some $k \geq 0$, and let $G$ have $k+1$ doubly directed regions.  Let $C$ be any cycle of the doubly directed graph of $G$.  
 
 There are two ways to fix the edges of $C$ such that at every corner of $C$ the two arcs are either both pointing into or both pointing out of the corner and at every vertex of $C$ that is not a corner one arc enters the vertex and one leaves.  Let $G_1$ and $G_2$ be the digraphs that are obtained from $G$ by fixing $C$ in one of these two ways.   
 
 Each of the subgraphs $G_1$ and $G_2$ will satisfy the same conditions that $G$ did, that is, at each vertex not incident to a doubly directed edge one of the configurations of Figure \ref{simple flow grid model} occurs and at each vertex incident to two doubly directed edges one of the configurations of Lemma \ref{degree2vertices} occurs.
 
$G_1$ and $G_2$ have fewer doubly directed regions than $G$ does, so by induction each of them is an elementary flow grid.  Let $S_1$ be a set of simple flow grids whose union is $G_1$ and let $S_2$ be a set of simple flow grids whose union is $G_2$.  Then $G$ is $S_1 \cup S_2$, so $G$ is also an elementary flow grid.   
\end{proof}

\begin{theorem} \label{ifevenallcyclesedges}
If an elementary flow grid $G$ has all vertices of its doubly directed graph of even degree then the matrix of every simple doubly directed cycle $C$ of $G$ is parallel to an edge of $f(G)$.
\end{theorem}
\begin{proof}
Let $G$ be an elementary flow grid with all vertices of even degree and let $C$ be a cycle in the doubly directed graph of $G$.  Let $G_1$ be obtained by fixing the edges of $C$ in such a way that at every corner of $C$, the two arcs are either both pointing into or both pointing out of the corner, and at every vertex of $C$ that is not a corner, one arc enters the vertex and one leaves. By Lemma \ref{even elementary flow grid}, $G_1$ is an elementary flow grid.  Let $A_1$ be one of the simple flow grids in a set whose union is $G_1$.  Let $A_2$ be the simple flow grid obtained from $A_1$ by fixing the edges of $C$ in the other one of the two ways to fix it. Then $A_1 \cup A_2$ is an elementary flow grid whose doubly directed graph is $C$.
\end{proof}

\section{Products}  

In the previous section, we showed how a basis for the subspace parallel to a face can be obtained from a union of the bases for the cycles of the 2-connected components.  This section is devoted to showing that a face is affinely equivalent to the product of the faces whose doubly directed graphs are edge-disjoint unions of 2-connected components of its doubly directed graph.

\begin{definition} \cite{Ziegler1}
 Let $P_1\subseteq \mathbb{R}^{d_1}$ and $P_2\subseteq \mathbb{R}^{d_2}$  be polytopes.  The {\it product} of $P_1$ and $P_2$ is $P_1\times P_2=\left\{\left(\begin{array}{c}w_1\\w_2\end{array}\right):   w_1\in P_1, w_2 \in P_2\right\}$.
\end{definition}

Let $F$ be a face of $ASM_n$ such that ${\mathcal G}(F)$  has two subgraphs ${\mathcal G}_1$ and ${\mathcal G}_2$ which are unions of 2-connected components of ${\mathcal G}(F)$, contain all the edges and vertices of ${\mathcal G}(F)$, and have no edges in common.  Let $A\in F$ be an ASM.  

Given an ASM $B \in F$,  define subgraphs $g_1(B)$ and $g_2(B)$ of $g(F)$ as follows.  

For $i=1,2$ and each edge $e$ of the grid graph, let $g_i(B)$ have a directed edge that agrees in orientation with $g(B)$ if $e$ is an edge of ${\mathcal G}_i$ and agrees in orientation with $g(A)$ otherwise.  


\begin{lemma}\label{edgeinface} Let $F$ be a face of $ASM_n$, let $A$ be a vertex of $F$ and let $e$ be an edge of ${\mathcal G}(F)$. There exists an edge $F'$ of $F$ such that $A\in F'$ and $e$ is an edge of ${\mathcal G}(F')$.  \end{lemma}

\begin{proof}  
The set of vertices $B$ of $F$ for which $g(B)$ agrees with $g(A)$ in orientation on edge $e$ is contained in a facet $F''$ of $F$.  There must be an edge $F'$ of $F$ containing $A$ that leaves $F''$ and contains a vertex $A'$ of $F$ for which $g(A')$ differs from $g(A)$ on edge $e$.  The cycle of ${\mathcal G}(F)$ corresponding to the edge of $F$ between $A$ and $A'$ contains $e$.  
\end{proof}

\begin{proposition}
For any $B\in F$, each of the graphs $g_1(B)$ and $g_2(B)$ is a simple flow grid.
\end{proposition}


\begin{proof}  
We first show that $g_1(B)$ is a simple flow grid.  
The proof is by induction on the dimension of $F(\{A,B\})$.  If the dimension is zero, then $A=B$ and the claim is true.  Suppose that the desired simple flow grid exists if the dimension of $F(\{A,B\})$ is $k$ or less and let the dimension of $F(\{A,B\})$ be $k+1$. If ${\mathcal G}_2$ contains no grid edge of $\mathcal{G}(F(\{A,B\}))$, $g(B)$ is the desired simple flow grid.

By Lemma \ref{edgeinface} if ${\mathcal G}_2$ contains a grid edge $e$ of $\mathcal{G}(F(\{A,B\}))$ then there is a polytope edge $F'$ of $F(\{A,B\})$ incident to $B$ with $e\in {\mathcal G}(F')$.  Let $B'$ be the other vertex  of $F'$.  The face $F(\{A,B'\})$ has smaller dimension than $F(\{A,B\})$, so by induction, $g_1(B')$ is a simple flow grid that agrees with  $g(B')$ in ${\mathcal G}_1$ and agrees with $A$ otherwise.  But $g(B)$ agrees with $g(B')$ on ${\mathcal G}_1$ because the cycle corresponding to the polytope edge containing $B$ and $B'$ must be contained in ${\mathcal G}_2$.  Therefore, $g_1(B)$ is the simple flow grid $g_1(B')$.  

An analogous argument shows that $g_2(B)$ is a simple flow grid.  \end{proof}

For each vertex $B$ of $F$, let $C_1(B)$ be the ASM corresponding to $g_1(B)$ and let $C_2(B)$ be the ASM corresponding to $g_2(B)$.  
Let $F_1=\mbox{conv}\{C_1(B):B\in F\}$ and  $F_2=\mbox{conv}\{C_2(B):B\in F\}$.  $F_1$ is a face of $F$ because it corresponds to the elementary flow grid with double directed edges corresponding to edges of ${\mathcal G}_1$ and single directed edges agreeing with $A$ otherwise.   Similarly, $F_2$ is a face of $F$.  


We now state the main result of this section.

\begin{theorem}\label{productmaintheorem} The function given by $f(B)=(C_1(B),C_2(B))$ is a bijective affine transformation from $F$ to $F_1\times F_2$.
\end{theorem}

\begin{proof}
    The transformation $f$ is the composition of three transformations:  The first subtracts $A$ from $B$, the second sends $B-A$ to $(C_1(B)-A,C_2(B)-A)$, and the third adds $(A,A)$ to the result of the second transformation.  

The second transformation is a linear transformation that is bijective on the set of matrices in $\{B-A: B\in F\}$.
\end{proof}

\begin{figure}[h!] 
\centering 
\begin{tikzpicture}
\draw [stealth-](0.5,1) -- (1,1);
\draw [stealth-](0.5,2) -- (1,2);
\draw [stealth-](0.5,3) -- (1,3);

\draw [-stealth](3,1) -- (3.5,1);
\draw [-stealth](3,2) -- (3.5,2);
\draw [-stealth](3,3) -- (3.5,3);

\draw [-stealth](1,3) -- (1,3.5);
\draw [-stealth](2,3) -- (2,3.5);
\draw [-stealth](3,3) -- (3,3.5);

\draw [-stealth](1,1) -- (1,0.5);
\draw [-stealth](2,1) -- (2,0.5);
\draw [-stealth](3,1) -- (3,0.5);

\draw [stealth-](1,1) -- (2,1);
\draw [red,thick,{stealth-stealth}](1,2) -- (2,2);
\draw [red,thick,{stealth-stealth}](1,3) -- (2,3);

\draw [red,thick,{stealth-stealth}](2,1) -- (3,1);
\draw [red,thick,{stealth-stealth}](2,2) -- (3,2);
\draw [-stealth](2,3) -- (3,3);

\draw [stealth-](1,1) -- (1,2);
\draw [red,thick,{stealth-stealth}](1,2) -- (1,3);

\draw [red,thick,{stealth-stealth}](2,1) -- (2,2);
\draw [red,thick,{stealth-stealth}](2,2) -- (2,3);

\draw [red,thick,{stealth-stealth}](3,1) -- (3,2);
\draw [-stealth](3,2) -- (3,3);

\node at (1.5,2.5){${\mathcal G}_1$};
\node at (2.5,1.5){${\mathcal G}_2$};

\end{tikzpicture}
\caption{\label{product example}Previously seen in Figure \ref{2edgenot2vertex}, this elementary flow grid corresponds to a product of two edges.}
\end{figure}
For example, let $F$ be the face with elementary flow grid in Figure \ref{product example}.  Its four vertices are: 
     \[
A=\left[ { \begin{array}{ccc}
    1 & 0 & 0 \\
    0 & 1 & 0 \\
    0 & 0 & 1 \\
  \end{array} } \right],
B=\left[ { \begin{array}{ccc}
    0 & 1 & 0 \\
    1 & -1 & 1 \\
    0 & 1 & 0 \\
  \end{array} } \right],
 \left[ { \begin{array}{ccc}
    1 & 0 & 0 \\
    0 & 0 & 1 \\
    0 & 1 & 0 \\
  \end{array} } \right],
\left[ { \begin{array}{ccc}
    0 & 1 & 0 \\
    1 & 0 & 0 \\
    0 & 0 & 1 \\
  \end{array} } \right].
  \]


 $F_1= \mbox{conv}\left\{\left[ { \begin{array}{ccc}
    1 & 0 & 0 \\
    0 & 1 & 0 \\
    0 & 0 & 1 \\
  \end{array} } \right],
\left[ { \begin{array}{ccc}
    0 & 1 & 0 \\
    1 & 0 & 0 \\
    0 & 0 & 1 \\
  \end{array} } \right]\right\}$, which agree with $A$ on entries corresponding to grid vertices not in $V({\mathcal G}_1)$, and 
  
 $F_2=\mbox{conv}\left\{ \left[ { \begin{array}{ccc}
    1 & 0 & 0 \\
    0 & 0 & 1 \\
    0 & 1 & 0 \\
  \end{array} } \right],
\left[ { \begin{array}{ccc}
    1 & 0 & 0 \\
    0 & 1 & 0 \\
    0 & 0 & 1 \\
  \end{array} } \right]\right\}$, which agree with $A$ on entries corresponding to vertices not in $V({\mathcal G}_2)$.
  
  We trace the mapping of matrix $B$:
     \[ B=
\left[ { \begin{array}{ccc}
    0 & 1 & 0 \\
    1 & -1 & 1 \\
    0 & 1 & 0 \\ 
 \end{array} } \right]
  \rightarrow
 \left[ { \begin{array}{ccc}
    -1 & 1 & 0 \\
    1 & -2 & 1 \\
    0 & 1 & -1 \\
  \end{array} } \right]
  \rightarrow
\left( \left[ { \begin{array}{ccc}
    -1 & 1 & 0 \\
    1 & -1 & 0 \\
    0 & 0 & 0 \\
  \end{array} } \right]
 , \left[ { \begin{array}{ccc}
    0 & 0 & 0 \\
    0 & -1 & 1 \\
    0 & 1 & -1 \\
  \end{array} } \right]
  \right)
  \]
  \[
  \rightarrow
  \left(
\left[ { \begin{array}{ccc}
    0 & 1 & 0 \\
    1 & 0 & 0 \\
    0 & 0 & 1 \\
  \end{array} } \right]
  ,
  \left[ { \begin{array}{ccc}
    1 & 0 & 0 \\
    0 & 0 & 1 \\
    0 & 1 & 0 \\
  \end{array} } \right]
  \right) =(C_1(B),C_2(B))
\]

\begin{proposition}\label{productimpliesedgedisjoint}
    Suppose that a face $F$ of $ASM_n$ is combinatorially equivalent to the product of two faces $F_1$ and $F_2$ of $F$.  Then the edge sets of ${\mathcal G}(F_1)$ and ${\mathcal G}(F_2)$ partition the edge set of ${\mathcal G}(F)$.
\end{proposition}

\begin{proof}
    Suppose that $e$ is an edge of ${\mathcal G}(F)$.  Lemma \ref{edgeinface} says that there is an edge $F'$ of $F$ such that ${\mathcal G}(F')$ contains $e$.  Every edge of $F_1\times F_2$ is of the form $A\times F''$ where $A$ is a vertex of $F_1$ and $F''$ is an edge of $F_2$ or of the form $F''\times A$ for an edge $F''$ of $F_1$ and a vertex $A$ of $F_2$ (\cite [Page 10]{Ziegler1}).  
    Thus $e\in {\mathcal G}(F_1)$ or $e \in {\mathcal G}(F_2)$.  Suppose that $F_1'$ is an edge of $F_1$ and $F_2'$ is an edge of $F_2$.  Then $F_1'\times F_2'$ is a face of $F_1\times F_2$ combinatorially equivalent to a square.  Thus ${\mathcal G}(F_1'\times F_2')$ has no vertices of degree 3, by the argument after Figure \ref{Two ways to fix a ear}, so ${\mathcal G}(F_1')$ and ${\mathcal G}(F_2')$ are edge-disjoint.  By Lemma \ref{edgeinface}, ${\mathcal G(F_1)}$ is the union of ${\mathcal G}(F_1')$ for edges $F_1'$ of $F_1$ and ${\mathcal G(F_2)}$ is the union of ${\mathcal G}(F_2')$ for edges $F_2'$ of $F_2$, so ${\mathcal G(F_1)}$ and ${\mathcal G(F_2)}$ are edge-disjoint.

\end{proof}

\begin{proposition}\label{2 connected product}
 A face of $ASM_n$ is affinely equivalent to a product of two or more of its faces, each of dimension at least 1, if and only if its doubly directed graph is not 2-connected.   
\end{proposition}

\begin{proof} 
Suppose the doubly directed graph ${\mathcal G}(F)$ of a face $F$ of $ASM_n$ is not 2-connected.  Then there exist subgraphs ${\mathcal G}_1$ and ${\mathcal G}_2$ of ${\mathcal G}(F)$ that each contain cycles, such that the faces $F_1$ and $F_2$ from Theorem  \ref{productmaintheorem} each have dimension  at least 1.  Theorem \ref{productmaintheorem} then says that $F$ is affinely equivalent to $F_1\times F_2$.

Conversely, suppose that $F$ is affinely equivalent to $F_1\times F_2$, where $F_1$ and $F_2$ are faces of $F$, each of dimension at least 1.  Proposition \ref{productimpliesedgedisjoint} implies that ${\mathcal G(F_1)}$ and ${\mathcal G(F_2)}$ partition the edge set of ${\mathcal G}(F)$.  Each of ${\mathcal G(F_1)}$ and ${\mathcal G(F_2)}$ contains a cycle. 

We proceed by induction on the dimension of $F_1\times F_2$ to show that no cycle of ${\mathcal G}(F_1\times F_2)$ contains edges of both 
${\mathcal G(F_1)}$ and ${\mathcal G(F_2)}$.  If dim$(F_1\times F_2)=2$, then ${\mathcal G}(F_1\times F_2)$ is the edge-disjoint union of ${\mathcal G}(F_1)$ and ${\mathcal G}(F_2)$, which are both cycles.  Now suppose that the assertion holds if the dimension of the product is at least 2 and less than $k$, and let dim$(F_1\times F_2) = k$. Suppose $C$ is a cycle of ${\mathcal G}(F_1\times F_2)$ containing edges of both ${\mathcal G}(F_1)$ and ${\mathcal G}(F_2)$. 

By \cite [Page 10]{Ziegler1}, every facet of $F_1\times F_2$ is the product $F_1$ with of a facet of $F_2$ or the product of $F_2$ with a facet of $F_1$. Let $F_1'\times F_2$ be a facet of $F_1\times F_2$, where $F_1'$ is a facet of $F_1$.    ${\mathcal G}(F_1'\times F_2)$ must be obtained from fixing an ear $S$ of ${\mathcal G}(F_1\times F_2)$. 
If $C$ has no edges of $S$, then ${\mathcal G}(F_1'\times F_2)$ contains $C$ with edges in both ${\mathcal G}(F_1)$ and ${\mathcal G}(F_2)$, contradicting the inductive hypothesis. 
Thus $C$ must contain $S$.  $S$ is a path in ${\mathcal G}(F_1)$, and there must be a cycle $C'$ of ${\mathcal G}(F_1)$ containing $S$.  The union of $C'\backslash{S}$ and $C\backslash{S}$ contains a cycle $C''$ containing edges of both ${\mathcal G}(F_1)$ and ${\mathcal G}(F_2)$ but no edges of $S$.  This cycle $C''$ is contained in ${\mathcal G}(F_1'\times F_2)$, contradicting the inductive hypothesis.  Because ${\mathcal G}(F_1\times F_2)$ has no cycle containing edges of both ${\mathcal G}(F_1)$ and ${\mathcal G}(F_2)$, the edges of ${\mathcal G}(F_1)$ cannot be in the same 2-connected component of ${\mathcal G}(F_1\times F_2)$ as ${\mathcal G}(F_2)$ \cite[Lemma 3.1.3]{Diestel}. \end{proof}

\section{Facets of Faces}

To find the facets of faces we look for minimal nonempty sets of edges of the doubly directed graph that could be replaced by singly directed edges to form other elementary flow grids.  These minimal sets are unions of what we call ears.
Recall from Definition \ref{eardefinition} that an ear of ${\mathcal G}(F)$ is a path between vertices of degree 3 or 4 in which every interior vertex has degree 2, or a cycle containing at most one vertex of degree 3 or 4.   Recall from Lemma \ref{twowaystofix} that if $F$ is a face of $ASM_n$ and $S$ is an ear of ${\mathcal G}(F)$, then there are two orientations of the edges in $S$ that may appear in $g(F')$ for a face $F'$ of $F$, and these two orientations are opposite to each other.

\begin{proposition}\label{Two level}
    If $F$ is a face of $ASM_n$ and $S$ is an ear of ${\mathcal G}(F)$, let $X_1$ be the set of vertices of $F$ with simple flow grid having one of the two orientations of the edges in $S$ and let $X_2$ be the set of vertices of $F$ with simple flow grid having the other orientation of the edges in $S$. 
    Then $\mbox{conv}(X_1)$ and $\mbox{conv}(X_2)$ are disjoint faces of $F$ containing all the vertices of $F$, where $\mbox{conv}(X_1)$ is the convex hull of the set of vertices in $X_1$.
\end{proposition} 

\begin{proof}
    By Lemma \ref{twowaystofix}, the union of the disjoint sets $X_1$ and $X_2$ is the set of all of the vertices of $F$.  Let $e\in S$ be an edge.  Without loss of generality, suppose $e$ has endpoints $(i,j)$ and $(i+1,j).$  The two ways to fix $e$ correspond to adding the equations $\sum_{i'=1}^ia_{i'j}=1$ and $\sum_{i'=1}^ia_{i'j}=0$ to the inequalities and equations defining $F$.  The two hyperplanes defined by these equations are parallel. Because every matrix $A\in ASM_n$ satisfies $0\le \sum_{i'=1}^ia_{i'j}\le 1$, the two hyperplanes determine faces of $F$.  
\end{proof}

 The two faces $\mbox{conv}(X_1)$ and $\mbox{conv}(X_2)$ may be facets of $F$, but it is possible for neither to be a facet of $F$.

A facet of a face $F$ corresponds to a maximal set of ears in ${\mathcal G}(F)$ that can be fixed to reduce the number of doubly directed regions by no more than one.  Sometimes a single ear is such a set, and sometimes fixing one ear creates a doubly directed graph that has cut edges that can be fixed without further lowering the number of doubly directed regions.  

\begin{proposition}\label{facetcase(1)and(2)}
    Let $F$ and $F'$ be faces of $ASM_n$, let $F'$ be a face of $F$ and let ${\mathcal G}(F)$ be 2-connected.  Then $F'$ is a facet of $F$ if and only if one of the following holds.  
    \begin{enumerate}
        \item ${\mathcal G}(F')$ is 2-connected and $g(F')$ is obtained from $g(F)$ by fixing one ear.
        \item ${\mathcal G}(F')$ is not 2-connected and ${\mathcal G}(F)$can be expressed as a chain $(P_1,Q_1,P_2,\newline Q_2,\ldots,P_k,Q_k)$, where $k\ge 2$ and $Q_1,Q_2,\ldots,Q_k$ edge-disjoint 2-connected components of ${\mathcal G}(F)$; for each $i=1,2,\ldots,k$, $P_i$ is an ear (possibly a single vertex) connecting $Q_{i-1}$ to $Q_i$, with indices modulo $k$, and ${\mathcal G}(F')$ is obtained from ${\mathcal G}(F)$ by removing the interiors of the $P_i$'s.  
    \end{enumerate}
\end{proposition}

These possibilities are analogous to those described in \cite[Corollary 2.11]{BruGibI} and \cite[Theorem 4.1]{Paffenholz1}.  The term {\it fully indecomposable} in \cite{BruGibI} corresponds to 2-connected and the term {\it has total support} in \cite{BruGibI} means that the matrix corresponds to a face of $B_n$. 

\begin{proof}
    Suppose $F'$ is a facet of $F$. There must be an ear $P$ of ${\mathcal G}(F)$ that is not in ${\mathcal G}(F')$.   
    Let ${\mathcal G}$ be obtained from ${\mathcal G}(F)$ by removing the edges and interior vertices of $P$.  If ${\mathcal G}$ is 2-connected, then fixing any other ears would reduce the number of regions further, so ${\mathcal G}={\mathcal G}(F')$ and we have case (1). 
    If ${\mathcal G}$ is not 2-connected, then one can create the graph $B(G)$ which has a vertex $v_Q$ for each 2-connected component $Q$, a vertex $v_P$ for each path (possibly a single vertex) joining two connected components, and an edge between $v_P$ and $v_Q$ whenever $P$ has an endpoint in $Q$. 
    The graph $B(G)$ must be a tree, (See, e.g. \cite[Section 6.3]{Agnarsson1}).  Because the original graph ${\mathcal G}(F)$ is 2-connected, $B(G)$ must actually be a path.  Thus we can define $P_1=P$ and $(Q_1,P_2,Q_2,\ldots,P_k,Q_k)$ be the sequence of 2-connected components $Q_i$ and paths corresponding to the path $B(G)$.  For ${\mathcal G}(F)$ to be 2-connected, $P_1$ must connect $Q_1$ to $Q_k$. Theorem \ref{theoremB} implies that none of the edges of $P_2,P_3,\ldots,P_k$ can be in ${\mathcal G}(F')$, and removing any edges of $Q_1,\ldots,Q_k$ from ${\mathcal G}$ would reduce the number of regions further.  Thus ${\mathcal G}(F')$ is obtained from ${\mathcal G}(F)$ by removing the interiors of the $P_i$'s.

    If either (1) or (2) occurs, then the number of regions of ${\mathcal G}(F)$ is one more than the number of regions of ${\mathcal G}(F')$, so $F'$ is a facet of $F$.  
\end{proof}

\begin{lemma}
    Suppose that $F'$ is a facet of $F$ and case (2) of Proposition \ref{facetcase(1)and(2)} holds.  Then there are only two ways to fix the collection of ears $P_1,P_2,\ldots,P_k$ and these two ways are opposite to each other.  
\end{lemma}  

\begin{proof}
    Suppose that there are $i$ and $j$ such that $P_i$ and $P_j$ each have at least one edge.  If the ear $P_i$ has been fixed, then $P_j$ cannot be in the resulting doubly directed graph, by Theorem \ref{theoremB}.  Thus fixing $P_i$ determines a unique way to fix $P_j$.  Thus fixing one ear determines how each of the other ears is to be fixed.  There are two ways to fix $P_i$, by Lemma \ref{twowaystofix}.
\end{proof}  

 2-level polytopes have been the subject of several recent papers \cite{enum2level, 2levelcombsettings}.

\begin{definition}
    \cite{enum2level} A (convex) polytope $P \subseteq \mathbb{R}^d$ is said to be a 2-level polytope if every hyperplane $H$ that is facet-defining for $P$ has a parallel hyperplane $H'$ that contains all the vertices of $P$ which are not contained in $H$.
\end{definition}

The decomposition in Proposition \ref{Two level} gives us the following Corollary. 

\begin{corollary}
    Every face of $ASM_n$ is a 2-level polytope.
\end{corollary}

Among several classes of polytopes that are known to be 2-level are the faces of the Birkhoff polytope \cite{2levelcombsettings}.
\begin{proposition} \cite[Lemma 5]{enum2level} If two 2-level polytopes are combinatorially equivalent, then they are affinely equivalent.  
\end{proposition}

In particular, this implies that the 4-gon faces of $ASM_n$ are all parallelograms and the octahedron faces of $ASM_n$ are all affine images of the standard octahedron.  See \cite[p. 88]{Dinkelman1}, for an affine map from the standard octahedron (convex hull of $\pm e_i$ in $\mathbb{R}^3$) to a face of $ASM_4$.  

Another useful property of 2-level polytopes is the following.

\begin{proposition}\label{2level01polytope}\cite[Lemma 9]{enum2level} If $F$ is a $d$-dimensional 2-level polytope, then $F$ is affinely equivalent to the convex hull of a set of vertices of the standard $d$-cube $[0,1]^d$.  
\end{proposition}

The following corollary can be proved for faces of $ASM_n$ using what we have proved, but we state it with a more general hypothesis, since it follows from Proposition \ref{2level01polytope}.

\begin{corollary}\label{maxvertices}
    The number of vertices of a $d$-dimensional 2-level polytope is at most $2^d$.  If a $d$-dimensional 2-level polytope has $2^d$ vertices, then it is affinely equivalent to a $d$-cube.
\end{corollary}

In preparation for the next section we examine what happens when an ear incident to a degree three vertex is fixed.

\begin{lemma}\label{fixdegree3oneway}
Suppose $(i,j)$ is a degree 3 vertex in ${\mathcal G}({F})$ and $S$ is an ear of ${\mathcal G}(F)$ with endpoint $(i,j)$.  Then, at most one of the ways to fix the edges of $S$ defines a facet of $F$.   
\end{lemma}

\begin{proof}
Let $e$ be the fixed edge incident to vertex $(i,j)$ and let $f$ be a doubly directed edge incident to vertex $(i,j)$.  For one way to fix edge $f$, we will have sgn$_F((i,j),e)=\mbox{sgn}_F((i,j),f)$.  It follows that the two other edges incident to $(i,j)$ must be fixed in order to satisfy $Sgn_F(\{v\})=0.$ 
   If the face $F'$ of $F$ has all the edges of $g(F')$ incident to $(i,j)$ fixed, then the number of doubly directed regions of $g(F')$ will not be one less than the number of doubly directed regions of $g(F)$, so $F'$ will not be a facet of $F$. 
\end{proof}






\section{Inequalities}

In the following, assume that $s$ is the number of ears of ${\mathcal G}(F)$ and $f$ is the number of facets of a face $F$ of $ASM_n$.  We first derive an upper bound on the number of facets of a $d$-dimensional face of $ASM_n$.

Recall (\cite{Agnarsson1}) Euler's formula $V-E+R=2$ for a connected planar graph, where $V$ is the number of vertices, $E$ is the number of edges and $R$ is the number of 2-dimensional regions. Let $V_3$ and $V_4$ be the number of degree 3 and degree 4 vertices of ${\mathcal G}({F})$ and let $d$ be the dimension of $F$.  Define a planar graph $G$ with vertex set equal to the set of degree 3 and 4 vertices of ${\mathcal G}({F})$ and edge set equal to the set of ears of ${\mathcal G}({F})$.  Then the number of regions of $G$ is $d+1$ and Euler's formula applied to $G$ yields, for a connected doubly directed graph  ${\mathcal G}({F})$ that has at least one vertex of degree $>2$: 
\begin{equation}\label{Euler}V_3+V_4 - s + (d+1)=2.\end{equation} 

By counting the ear-vertex incidences, we get 
\begin{equation}\label{earvertexincidences} 2s=4V_4+3V_3.\end{equation}
Eliminating the number of ears by adding two times equation (\ref{Euler}) to equation (\ref{earvertexincidences}) gives us 
\begin{equation}\label{2V4+V3=2(d-1)}2V_4+V_3=2(d-1).\end{equation}

\begin{proposition}\label{connected4(d-1)} If the doubly directed graph of $F$ is connected and has at least one vertex of degree $>2$, then $F$ has at most $4(d-1)$ facets.    
\end{proposition}
\begin{proof}
Each ear incident to a degree 4 vertex can be fixed to form a facet in at most 2 ways, while each ear incident to a degree 3 vertex can be fixed to form a facet in at most 1 way.  Every facet is formed by fixing at least one ear.  
Similarly to equation (\ref{earvertexincidences}) we get \begin{equation}2f\le 8V_4 +3V_3,\end{equation}
which implies, by equation (\ref{2V4+V3=2(d-1)}), 
\begin{equation} f \le 4V_4+\frac{3}{2}V_3\le 4V_4+2V_3=4(d-1).\end{equation}
\end{proof}

We can generalize Proposition \ref{connected4(d-1)} to the case where the doubly directed graph of $F$ is not connected, using the material of Section 5.

\begin{theorem}\label{maxnumfacets}  
If a face $F$ of $ASM_n$ has dimension $d\ge 2$, then $F$ has at most $4(d-1)$ facets.   
\end{theorem}

\begin{proof}
    If ${\mathcal G}(F)$ is connected and the dimension of $F$ is at least 2, then ${\mathcal G}(F)$ has at least one vertex of degree greater than 2 and by Proposition \ref{connected4(d-1)}, $F$ has at most $4(d-1)$ facets.  Otherwise, by Theorem \ref{productmaintheorem}, $F$ is the product of faces $F_1$ and $F_2$ such that ${\mathcal{G}}(F_1)$ and ${\mathcal G}(F_2)$ have no vertices in common. From \cite[Page 10]{Ziegler1}, the dimension of $F$ is the dimension of $F_1$ plus the dimension of $F_2$ and the number of facets of $F$ is the number of facets of $F_1$ plus the number of facets of $F_2$.  If $F_1$ and $F_2$ both have dimension 1, then $F$ has dimension $1+1=2$ and $2+2=4=4((1+1)-1)$ facets. If $F_1$ has dimension $d_1=1$ and $F_2$ has dimension $d_2\ge 2$, then $F$ has dimension $d_2+1$ and at most $2+4(d_2-1)\le 4((d_2+1)-1)$ facets.  If $F_1$ has dimension $d_1\ge 2$ and $F_2$ has dimension $d_2\ge 2$, then $F$ has at most $4(d_1-1)+4(d_2-1)< 4(d_1+d_2-1)$ facets. By induction, the inequality holds for all $d\ge 2$.
\end{proof}

One can compare this with \cite[Corollary 3.4]{BruGibIII}, which states that a $d$-dimensional face of $B_n$ has at most $3(d-1)$ facets, for $d\ge 3$.  Note that in order to have equality in (5), it is necessary that $V_3=0$, so $F$ must be centrally symmetric.  
\begin{proposition}
 If the doubly directed graph of $F$ is connected and has at least one vertex of degree $>2$, then $2(d-1) \leq s \leq 3(d-1)$.    
\end{proposition}
\begin{proof}
From equation (2),
$2s = 4V_4 + 3V_3$.
From equation (1),
$s = V_4 + V_3+d-1$.
Then subtracting these two equations we get
$s = 3V_4 + 2V_3  - d+1$.
Then
$3V_4+ \frac{3}{2} V_3 -d+1 \leq s \leq 4V_4 + 2V_3 - d+1$.
But 
$2V_4+V_3=2(d-1)$ by equation (3).
Then, $3(d-1)-d+1 \leq s \leq 4(d-1) - d + 1 = 3(d-1)$. So $2(d-1) \leq s \leq 3(d-1)$.
\end{proof}

The next theorem involves faces $F$ for which ${\mathcal G}(F)$ is 2-connected.  It is therefore critical to cite the theorem of Whitney, \cite[Theorem 18]{Whitney1}, which states that a 2-connected graph that is not a cycle can always be obtained by adding a path to another 2-connected graph. We will need a related observation:

\begin{lemma}\label{removableearbydeg3vertex}
    If a graph ${\mathcal G}$ is 2-connected and $v$ is a vertex of ${\mathcal G}$ of degree 3, then we can remove one of the ears incident to $v$, except for its endpoints, without destroying 2-connectivity.
\end{lemma}

\begin{proof} By 2-connectivity, none of the ears incident to $v$ has both of its endpoints equal to $v$.  Let $P_1$, $P_2$, $P_3$ be the three ears incident to $v$ and let $v_1,v_2,v_3$ be the endpoints of the three ears that are not $v$.  If $v_i=v_j$ for some $i\neq j$, then we can remove $P_i$ without losing 2-connectivity, so assume $v_1,v_2, v_3$ are distinct.  There must be a cycle $C_1$ of ${\mathcal G}$ containing $P_1$.  We can without loss of generality assume that $C_1$ consists of $P_1$, $P_2$, and a path $P$ from $v_1$ to $v_2$.  Similarly there must be a cycle $C_2$ containing $P_3$, which can be assumed to be made up of $P_3$, $P_1$ and a path $P'$ from $v_1$ to $v_3$.  If we then remove path $P_1$, except for its endpoints, from ${\mathcal G}$ the resulting graph will be 2-connected, because we have a cycle consisting of $P_2,P,P',P_3$ which contains the endpoints of $P_1$.  
\end{proof}

For a face $F$ of $ASM_n$, let $v(F)$ be the number of vertices of $F$.
Theorem 5.1 of \cite{BruGibI} shows that if a $d$-dimensional face of the Birkhoff polytope $B_n$ corresponds to a connected subgraph of a complete bipartite graph, then the face contains at most $2^{d-1}+1$ vertices.  We will now show an analogous bound for faces of $ASM_n$, but slightly larger.

\begin{lemma}\label{chainofsimplecycles}
If the 2-connected doubly directed graph of a face $F$ is made up of a chain of cycles $(C_1,C_2,\ldots,C_t)$ and paths $P'_1,P'_2,\ldots,P'_t$, some but not all of which may be single vertices, such that $P'_i$ has one endpoint in $C_i$ and one endpoint in $C_{i+1}$ for $i=1,2,\ldots,t$, indices modulo $t$, then $F$ has at most $2^t+1$ vertices.  
\end{lemma}
\begin{proof} The proof is by induction on $d$.  If $t=1$, then $F$ is a triangle, with $2^1+1$ vertices.  Suppose the conclusion holds when $t=k\ge 1$ and let $t=k+1$. Assume $P'_{k+1}$ is not a single vertex. Consider an ear $P' \subseteq C_{k+1}$ from the intersection of the degree three vertex  $P'_{k}\cap C_{k+1}$ to the vertex $P'_{k+1}\cap C_{k+1}$.  One way to fix $P'$ leaves a polytope of dimension at most $k$ and the other leaves a polytope to which we can apply the inductive hypothesis.  By induction, the number of vertices of $F$ is at most $2^{k} + 2^{k}+1=2^{k+1}+1$.  
\end{proof}

\begin{theorem} \label{2connectedub} Suppose a face $F$ has dimension $d$ and ${\mathcal G}(F)$ is 2-connected.  Then $v(F)\le 2^{d-1}+2$ if ${\mathcal G}(F)$ has no degree 3 vertex, while $v(F)\le 2^{d-1}+1$ otherwise. 
\end{theorem}

\begin{proof}  The proof is by induction on $d$.  It is true for $d< 3$ by inspection.  

Suppose $d\ge 3$ and ${\mathcal G}(F)$ has no degree 3 vertices.   
Let $P$ be an ear of ${\mathcal G}(F)$ such that removing  $P$ (but not its endpoints) results in a 2-connected graph.  Let $F_1$ and $F_2$ be the two faces of $F$ obtained by fixing the edges of $P$.  If the endpoints of $P$ have degree 3 in ${\mathcal G}(F_1)$, then by induction $v(F_1)\le 2^{d-2}+1$.  
If an endpoint of $P$ has degree less than 3 in ${\mathcal G}(F_1)$ then $F_1$ is not a facet.  In that case $v(F_1)\le 2^{d-2}$, by Corollary \ref{maxvertices}.  The same conclusions hold for $F_2$.  Therefore, $v(F_1)+v(F_2)\le 2^{d-2}+1+2^{d-2}+1=2^{d-1}+2$.   

Next suppose that ${\mathcal G}(F)$ has a vertex $(i,j)$ of degree 3.  Let $P_1,P_2,P_3$ be the three ears emanating from $(i,j)$.  The three ears must be distinct, due to 2-connectivity.  By Lemma \ref{removableearbydeg3vertex}, we can assume that removing $P_1$ but not its endpoints leaves a 2-connected graph.  Let $F_{11}$ and $F_{12}$ be the faces obtained by the two ways of fixing $P_1$.  One of the ways to fix $P_1$ fixes all three ears incident to $(i,j)$ (See Lemma \ref{fixdegree3oneway}), and the resulting face, say $F_{11}$, has dimension at most $d-2$. Hence, by  Corollary \ref{maxvertices} $F_{11}$ has at most $2^{d-2}$ vertices.  

The other way of fixing $P_1$ leaves a 2-connected graph. If $F_{12}$ has dimension less than $d-1$, then $F$ has at most $2^{d-2}+2^{d-2}=2^{d-1}$ vertices.  Suppose $F_{12}$ is $(d-1)$-dimensional, so that $v(F_{12})\le 2^{d-2}+2$ by induction.  If $v(F_{11})<2^{d-2}$ then $v(F)\le 2^{d-2}-1+2^{d-2}+2=2^{d-1}+1$.  Therefore suppose $v(F_{11})=2^{d-2}$, so by Corollary \ref{maxvertices}, $F_{11}$ is a $(d-2)$-cube and by Proposition \ref{productimpliesedgedisjoint}, ${\mathcal G}(F_{11})$ has $d-2$ two-connected components which are all cycles.   
Let $C_1,C_2,\ldots,C_{d-2}$ be the cycles of ${\mathcal G}(F_{11})$.  These cycles are contained in ${\mathcal G}(F_{12})$, and ${\mathcal G}_(F_{12})$ contains one more cycle, which must contain ears $P_2$ and $P_3$. 
Because ${\mathcal G}(F_{12})$ is 2-connected, ${\mathcal G}(F_{12})$ is the union of $C_1,C_2,\ldots,C_{d-2}$ and paths $P'_1,P'_2,\ldots,P'_{d-2}$ of which some may be single vertices but one of which contains $P_1\cup P_2$, such that for $i=1,2,\ldots,d-2$, $P'_i$ has one endpoint in $C_i$ and one endpoint in $C_{i+1}$, indices modulo $d-2$.

By Lemma \ref{chainofsimplecycles}, $F_{12}$ has at most $2^{d-2}+1$ vertices, so $F$ has at most $2^{d-2}+2^{d-2}+1=2^{d-1}+1$ vertices.  



\end{proof}

Examples of faces that meet the upper bounds of both Theorems \ref{maxnumfacets} and \ref{2connectedub} are bipyramids over cubes.  Such examples can be obtained by taking a doubly directed graph whose doubly directed regions correspond to a cycle of cycles bounding a doubly directed region.  For example, see Figure \ref{bipyramidover3cube}.

\begin{figure}[h!] 
\centering 
\begin{tikzpicture}

\draw [stealth-](5.5,1) -- (6,1);
\draw [stealth-](5.5,2) -- (6,2);
\draw [stealth-](5.5,3) -- (6,3);
\draw [stealth-](5.5,4) -- (6,4);

\draw [-stealth](9,1) -- (9.5,1);
\draw [-stealth](9,2) -- (9.5,2);
\draw [-stealth](9,3) -- (9.5,3);
\draw [-stealth](9,4) -- (9.5,4);

\draw [-stealth](6,4) -- (6,4.5);
\draw [-stealth](7,4) -- (7,4.5);
\draw [-stealth](8,4) -- (8,4.5);
\draw [-stealth](9,4) -- (9,4.5);

\draw [-stealth](6,1) -- (6,0.5);
\draw [-stealth](7,1) -- (7,0.5);
\draw [-stealth](8,1) -- (8,0.5);
\draw [-stealth](9,1) -- (9,0.5);

\draw [stealth-](6,1) -- (7,1);
\draw [red,thick,{stealth-stealth}](6,2) -- (7,2);
\draw [stealth-](6,3) -- (7,3);
\draw [red,thick,{stealth-stealth}](6,4) -- (7,4);

\draw [red,thick,{stealth-stealth}](7,1) -- (8,1);
\draw [red,thick,{stealth-stealth}](7,2) -- (8,2);
\draw [red,thick,{stealth-stealth}](7,3) -- (8,3);
\draw [red,thick,{stealth-stealth}](7,4) -- (8,4);

\draw [-stealth](8,1) -- (9,1);
\draw [red,thick,{stealth-stealth}](8,2) -- (9,2);
\draw [red,thick,{stealth-stealth}](8,3) -- (9,3);
\draw [-stealth](8,4) -- (9,4);

\draw [stealth-](6,1) -- (6,2);
\draw [red,thick,{stealth-stealth}](6,2) -- (6,3);
\draw [red,thick,{stealth-stealth}](6,3) -- (6,4);

\draw [red,thick,{stealth-stealth}](7,1) -- (7,2);
\draw [red,thick,{stealth-stealth}](7,2) -- (7,3);
\draw [-stealth](7,3) -- (7,4);

\draw [red,thick,{stealth-stealth}](8,1) -- (8,2);
\draw [red,thick,{stealth-stealth}](8,2) -- (8,3);
\draw [red,thick,{stealth-stealth}](8,3) -- (8,4);

\draw [stealth-](9,1) -- (9,2);
\draw [red,thick,{stealth-stealth}](9,2) -- (9,3);
\draw [-stealth](9,3) -- (9,4);

\end{tikzpicture}
\caption{\label{bipyramidover3cube}Elementary flow grid for a bipyramid over a 3-cube.}
\end{figure}

Theorem \ref{2connectedub} can be used to prove an inequality that is known for many classes of 2-level polytopes \cite{2levelcombsettings}.
\begin{corollary}
    Suppose that $F$ is a $d$-dimensional face of $ASM_n$ for some $n$.  If $v$ and $f$ are the numbers of vertices and facets of $F$, then $vf\le d2^{d+1}$.  
\end{corollary}
\begin{proof}
    It was shown in \cite[Lemma 2.2]{2levelcombsettings} that if $F=F_1\times F_2$ where $F_1$ and $F_2$ satisfy the inequality, then $F$ also satisfies it.  On the other hand, if ${\mathcal G}(F)$ is 2-connected, then by Theorems \ref{maxnumfacets} and \ref{2connectedub}, $vf\le (2^{d-1}+2)(4(d-1))\le d2^{d+1}$.
\end{proof}

\section{No face of $ASM_n$ has the combinatorial type of $B_3$.}

The polytope $B_3$, the convex hull of the $3\times 3$ permutation matrices, is known to have dimension 4, 6 vertices, and 9 facets (\cite{Paffenholz1}).  
The purpose of this section is to prove Theorem \ref{nofaceB3}, showing that no face of $ASM_n$ is combinatorially equivalent to $B_3$.  The only properties of $B_3$ that we will use are its dimension, number of vertices and number of facets.  We will need the following lemmas, which derive properties that the doubly directed graph of a face combinatorially equivalent to $B_3$ would have to have.

\begin{lemma}\label{onefacetperear}
  Suppose that $F$ is a 4-dimensional face of $ASM_n$ for some $n$, with less than 8 vertices, and let $S$ be an ear of the doubly directed graph of $F$.  At most one of the two ways to fix $S$ defines a facet of $F$.  
\end{lemma}

\begin{proof} A facet of $F$ would have to contain at least 4 vertices.  The two disjoint faces of $F$ defined by fixing $S$ together have at most 7 vertices.    
\end{proof}  

Lemma \ref{onefacetperear} is next used to derive the number of ears and the number of vertices of degrees 3 and 4 in the doubly directed graph of a face combinatorially equivalent to $B_3$.

\begin{lemma}\label{9earsdistinctfacets} If $F$ is a 4-dimensional face of $ASM_n$ with 9 facets and less than 8 vertices, then the doubly directed graph of $F$ has 6 vertices of degree 3 and no vertices of degree 4.  Furthermore, ${\mathcal G}(F)$ has 9 ears which define 9 distinct facets.  \end{lemma} 

\begin{proof} By Lemma \ref{onefacetperear}, the doubly directed graph of $F$ would have to have at least 9 ears.  In that case, equation (\ref{Euler}) from the discussion before Proposition \ref{connected4(d-1)} leads to $V_3+V_4+3\ge 9$, or $2V_3+2V_4\ge 12.$  Equation (\ref{2V4+V3=2(d-1)}) tells us that $2V_4+V_3=6$, so subtracting gives us $V_3\ge 6.$ But $V_4\ge 0,$ so $V_3=6$ and $V_4=0$.    
\end{proof}

The following proposition applies in general to faces of $ASM_n$.

\begin{proposition}\label{twoearcutset}
    If, for a face $F$ of $ASM_n$,  ${\mathcal G}(F)$ is connected but removing two ears $S_1$ and $S_2$ from ${\mathcal G}(F)$ creates a disconnected graph, then fixing one of $S_1$ or $S_2$ also fixes the other ear.  
\end{proposition}

\begin{proof} If $S_1$ is fixed and $S_2$ is not, then the edges of $S_2$ are not contained in a cycle of the doubly directed graph of the face obtained by fixing $S_1$, contradicting Theorem \ref{theoremB}.
\end{proof}

Because, by Lemma \ref{9earsdistinctfacets}, each of the 9 ears of ${\mathcal G}(F)$ for a face $F$ combinatorially equivalent to $B_3$ must correspond to a different facet, we get the following corollary.

\begin{corollary}\label{twoearcutsetcor}If $F$ has the combinatorial type of $B_3$, then ${\mathcal G}(F)$ cannot have two ears whose removal disconnects ${\mathcal G}(F)$.    
\end{corollary}

We now have a characterization of ${\mathcal G}(F)$ for $F$ combinatorially equivalent to $B_3$.

\begin{proposition}\label{triangularprism}
    Let $F$ be a face of $ASM_n$ of the combinatorial type of $B_3$.  Let $G$ be the multigraph obtained from ${\mathcal G}(F)$ by retaining each vertex of degree 3 or 4 and replacing each ear by a single edge.  Then $G$ is isomorphic to the graph of the triangular prism.
\end{proposition}

\begin{proof}
    If $G$ has parallel edges joining vertices $v$ and $w$, then the remaining edges incident to $v$ and $w$ form a 2-edge cutset of $G$, contradicting Corollary \ref{twoearcutsetcor}.  Thus $G$ is a 3-regular planar graph on six vertices. It is easy to see that such a graph is isomorphic to the graph of the triangular prism. 
\end{proof}

   \begin{figure}[h!]
\centering 
\begin{tikzpicture}
\draw [stealth-](0.5,1) -- (1,1);
\draw [stealth-](0.5,2) -- (1,2);
\draw [stealth-](0.5,3) -- (1,3);
\draw [stealth-](0.5,4) -- (1,4);
\draw [stealth-](0.5,5) -- (1,5);

\draw [-stealth](5,1) -- (5.5,1);
\draw [-stealth](5,2) -- (5.5,2);
\draw [-stealth](5,3) -- (5.5,3);
\draw [-stealth](5,4) -- (5.5,4);
\draw [-stealth](5,5) -- (5.5,5);

\draw [-stealth](1,5) -- (1,5.5);
\draw [-stealth](2,5) -- (2,5.5);
\draw [-stealth](3,5) -- (3,5.5);
\draw [-stealth](4,5) -- (4,5.5);
\draw [-stealth](5,5) -- (5,5.5);

\draw [-stealth](1,1) -- (1,0.5);
\draw [-stealth](2,1) -- (2,0.5);
\draw [-stealth](3,1) -- (3,0.5);
\draw [-stealth](4,1) -- (4,0.5);
\draw [-stealth](5,1) -- (5,0.5);

\draw [stealth-](1,1) -- (2,1);
\draw [stealth-](1,2) -- (2,2);
\draw [-stealth](1,3) -- (2,3);
\draw [stealth-](1,4) -- (2,4);
\draw [stealth-](1,5) -- (2,5);

\node at (2.2,5.2){$1$};

\draw [stealth-](2,1) -- (3,1);
\draw [red,thick,{stealth-stealth}](2,2) -- (3,2);
\draw [red,thick,{stealth-stealth}](2,3) -- (3,3);
\draw [red,thick,{stealth-stealth}](2,4) -- (3,4);
\draw [red,thick,{stealth-stealth}](2,5) -- (3,5);

\node at (2.5,4.2){$5$};
\node at (3.5,4.2){$6$};
\node at (4.2,4.5){$2$};
\node at (1.8,3.5){$4$};
\node at (3.2,3.5){$8$};
\node at (4.2,3.2){$9$};
\node at (2.5,2.8){$7$};
\node at (5.2,2.5){$3$};

\draw [-stealth](3,1) -- (4,1);
\draw [red,thick,{stealth-stealth}](3,2) -- (4,2);
\draw [red,thick,{stealth-stealth}](3,3) -- (4,3);
\draw [red,thick,{stealth-stealth}](3,4) -- (4,4);
\draw [red,thick,{stealth-stealth}](3,5) -- (4,5);

\draw [-stealth](4,1) -- (5,1);
\draw [red,thick,{stealth-stealth}](4,2) -- (5,2);
\draw [-stealth](4,3) -- (5,3);
\draw [-stealth](4,4) -- (5,4);
\draw [red,thick,{stealth-stealth}](4,5) -- (5,5);

\draw [stealth-](1,1) -- (1,2);
\draw [stealth-](1,2) -- (1,3);
\draw [-stealth](1,3) -- (1,4);
\draw [-stealth](1,4) -- (1,5);

\draw [stealth-](2,1) -- (2,2);
\draw [red,thick,{stealth-stealth}](2,2) -- (2,3);
\draw [red,thick,{stealth-stealth}](2,3) -- (2,4);
\draw [red,thick,{stealth-stealth}](2,4) -- (2,5);

\draw [-stealth](3,1) -- (3,2);
\draw [-stealth](3,2) -- (3,3);
\draw [red,thick,{stealth-stealth}](3,3) -- (3,4);
\draw [-stealth](3,4) -- (3,5);

\draw [stealth-](4,1) -- (4,2);
\draw [stealth-](4,2) -- (4,3);
\draw [red,thick,{stealth-stealth}](4,3) -- (4,4);
\draw [red,thick,{stealth-stealth}](4,4) -- (4,5);

\draw [stealth-](5,1) -- (5,2);
\draw [red,thick,{stealth-stealth}](5,2) -- (5,3);
\draw [red,thick,{stealth-stealth}](5,3) -- (5,4);
\draw [red,thick,{stealth-stealth}](5,4) -- (5,5);
\end{tikzpicture}
\caption{\label{3reg9ears}3-regular doubly directed graph with 9 ears.}
\end{figure}

In the example of Figure \ref{3reg9ears}, there are six degree 3 vertices of the doubly directed graph: $(1,4),(2,2),(2,3),(2,4),(3,2),$ and $(3,3)$.  There are 9 ears, labelled 1 through 9.  There is a doubly directed region bounded by the 3 ears labelled 6, 8 and 9.  The cycle on the outside is made up of ears 1, 3 and 4.  These two cycles are connected by the ears 2, 5 and 7.  The face $F$ of $ASM_5$ corresponding to this elementary flow grid is a 2-fold iterated pyramid over a square base.  Each of the ways of fixing ears 3 and 9 gives a triangular face of $F$.  One of the ways to fix ear 5 gives the square face of $F$ and the other yields the edge containing the two vertices not in the square face.  

Recall from Proposition \ref{triangularprism} that because a triangular prism has two triangular faces, ${\mathcal G}(F)$ for a face $F$ combinatorially equivalent to $B_3$ must contain a cycle of three degree 3 vertices connected by ears.  
In the following, we will assume that $(i_1,j_1)$, $(i_2,j_2)$, and $(i_3,j_3)$ are vertices that have degree 3 in ${\mathcal G}(F)$ such that whenever $k\neq k'$ in $\{1,2,3\}$, there is an ear with endpoints $(i_k,j_k)$ and $(i_{k'},j_{k'})$.  See Figure \ref{3cycle}, where the black edges are the fixed edges.

  \begin{figure}[h!] 
\centering 
\begin{tikzpicture}

\draw [green,thick,{stealth-stealth}](0,2) -- (1,2);
\draw [green,thick,{stealth-stealth}](3,3) -- (3,4);
\draw [green,thick,{stealth-stealth}](4,0) -- (4,1);

\draw [red,thick,{stealth-stealth}](1,2) -- (2,2);
\draw [red,thick,{stealth-stealth}](1,3) -- (2,3);
\draw [red,thick,{stealth-stealth}](2,1) -- (3,1);
\draw [red,thick,{stealth-stealth}](2,3) -- (3,3);
\draw [red,thick,{stealth-stealth}](3,1) -- (4,1);
\draw [red,thick,{stealth-stealth}](3,2) -- (4,2);

\draw [thick,{-stealth}](3,3) -- (4,3);
\draw [thick,{stealth-}](4,1) -- (5,1);
\draw [thick,{-stealth}](1,1) -- (1,2);

\draw [red,thick,{stealth-stealth}](1,2) -- (1,3);
\draw [red,thick,{stealth-stealth}](2,1) -- (2,2);
\draw [red,thick,{stealth-stealth}](3,2) -- (3,3);
\draw [red,thick,{stealth-stealth}](4,1) -- (4,2);

\node at (0.4,2.3){$(i_3,j_3)$};
\node at (4.6,0.7){$(i_1,j_1)$};
\node at (3.6,3.4){$(i_2,j_2)$};

\end{tikzpicture}
\caption{\label{3cycle}
Cycle with 3 ears.} 
\end{figure}

We want to show that at most 2 of the ears connecting these 3 vertices can be fixed to form facets of $F$.  We first define variables that describe the fixed edges incident to the vertices $(i_k,j_k)$.  
For $i=1,2,3$, define $e_i$ to be the fixed edge incident to $(i_k,j_k)$.

                 


Suppose that we have fixed the ear connecting $(i_k,j_k)$ and $(i_{k'},j_{k'})$, for $k, k' \in \{1,2,3\}$.  Then we can define $f_{kk'}$ to be the edge of this ear incident to $(i_k,j_k)$.  
                 



 Recall that Lemma \ref{fixdegree3oneway} said that if vertex $(i_k,j_k)$ has degree 3 in ${\mathcal G}(F)$ then at most one of the ways to fix the ear from    $(i_k,j_k)$  to  $(i_{k'},j_{k'})$  defines a facet.  We identify one of the ways that does not define a facet, in terms of the sgn$_F$ function introduced at the end of Section 2.

\begin{lemma}\label{fixingnotfacet}
    The given way of fixing the ear from $(i_k,j_k)$ to $(i_{k'},j_{k'})$ does not define a facet of $F$ if sgn$_F((i_k,j_k),e_k)=$sgn$_F((i_k,j_k),f_{kk'}).$
    \end{lemma}

\begin{proof} If sgn$_F((i_k,j_k),e_k)=\mbox{sgn}_F((i_k,j_k),f_{kk'})$, it follows that the two other edges incident to $(i_k,j_k)$ must be fixed to satisfy $\sum_{e: e \text{ fixed, }(i_k,j_k) \in e} \mbox{sgn}_F((i_k,j_k),e) =0.$ This implies that the number of doubly connected regions decreases by more than one.
\end{proof} 

We now apply Lemma \ref{signs of ear ends} to the ear from $(i_k,j_k)$ to $(i_{k'},j_{k'})$.

\begin{lemma}\label{twoendpointsoffixedear}
    sgn$_F((i_k,j_k),f_{kk'})= -\text{sgn}_F((i_{k'},j_{k'}),f_{k'k})$
\end{lemma}


The following lemma uses Lemma \ref{twoendpointsoffixedear} and Lemma \ref{fixingnotfacet} to determine a case in which neither way to fix the ear between $(i_k,j_k)$ and $(i_{k'},j_{k'})$ defines a facet.

\begin{lemma}\label{neitherwayfacet}
    If sgn$_F((i_k,j_k),e_k)=\mbox{sgn}_F((i_{k'},j_{k'}),e_{k'})$,
    then neither way to fix the ear connecting $(i_k,j_k)$ and $(i_{k'},j_{k'})$ defines a facet.
\end{lemma}

\begin{proof}
   Suppose we fix the ear from $(i_k,j_k)$ to $(i_{k'},j_{k'})$ to define a facet. Then sgn$_F((i_k,j_k),e_k)=\mbox{sgn}_F((i_k,j_k),f_{kk'})$ 
   and sgn$_F((i_{k'},j_{k'}),e_{k'})=\mbox{sgn}_F((i_{k'},j_{k'}),f_{k'k})$
   by Lemma \ref{fixingnotfacet}. Then Lemma \ref{twoendpointsoffixedear} implies  sgn$_F((i_k,j_k),e_k)=-\mbox{sgn}_F((i_{k'},j_{k'}),e_{k'})$.
\end{proof}

We have proved that the values of sgn$_F((i_k,j_k),e_k)$ 
have to be different on the endpoints of an ear if one of the ways to fix the ear defines a facet.  

\begin{proposition}\label{Atmost2earsfacets}
    At most 2 of the ears connecting $(i_1,j_1),(i_2,j_2)$ and $(i_3,j_3)$ can be fixed to define facets.  
\end{proposition}

\begin{proof}
    The values sgn$_F((i_k,j_k),e_k)$ 
    are contained in the 2-element set $\{-1,1\}$.  Because the cycle $((i_1,j_1),(i_2,j_2),(i_3,j_3))$ has odd length, there must be a pair $(k,k')$ for which sgn$_F((i_k,j_k),e_k)=\mbox{sgn}_F((i_{k'},j_{k'}),e_{k'})$
    By Lemma \ref{neitherwayfacet}, neither way of fixing the ear connecting $(i_k,j_k)$ and $(i_{k'},j_{k'})$ defines a facet of $F$.  
\end{proof}

We finally tie the preceding sequence of lemmas together.

\begin{theorem}\label{nofaceB3}
    No face of $ASM_n$ has the combinatorial type of $B_3$.
\end{theorem}

\begin{proof}
    By Lemma \ref{9earsdistinctfacets}, the elementary flow grid for a face $F$ of the combinatorial type of $B_3$ would have 9 ears each of which defines one facet of $F$.  By Proposition \ref{triangularprism} this elementary flow grid would have a cycle made up of 3 ears.  By Proposition \ref{Atmost2earsfacets}, not all ears of this cycle can define distinct facets.
\end{proof}







One can try to extend the proof to other cases of a combinatorial type with less than $2d$ vertices and $3(d-1)$ facets, but the graph $G$ obtained from ${\mathcal G}(F)$ by replacing ears with edges may be bipartite.  Consider the face for which the elementary flow grid is in Figure \ref{derangementpoly}.  The vertices of  this face are the permutation matrices for the nine permutations $\sigma$ of $\{1,2,3,4\}$ such that $\sigma(1)\neq 2, \sigma(2)\neq 4$, $\sigma(3)\neq 1$ and $\sigma(4)\neq 3$, which are the images of the derangements of $\{1,2,3,4\}$ under a permutation.

\begin{figure}[h!] 
\centering 
\begin{tikzpicture}

\draw [stealth-](5.5,1) -- (6,1);
\draw [stealth-](5.5,2) -- (6,2);
\draw [stealth-](5.5,3) -- (6,3);
\draw [stealth-](5.5,4) -- (6,4);

\draw [-stealth](9,1) -- (9.5,1);
\draw [-stealth](9,2) -- (9.5,2);
\draw [-stealth](9,3) -- (9.5,3);
\draw [-stealth](9,4) -- (9.5,4);

\draw [-stealth](6,4) -- (6,4.5);
\draw [-stealth](7,4) -- (7,4.5);
\draw [-stealth](8,4) -- (8,4.5);
\draw [-stealth](9,4) -- (9,4.5);

\draw [-stealth](6,1) -- (6,0.5);
\draw [-stealth](7,1) -- (7,0.5);
\draw [-stealth](8,1) -- (8,0.5);
\draw [-stealth](9,1) -- (9,0.5);

\draw [red,thick,{stealth-stealth}](6,1) -- (7,1);
\draw [stealth-](6,2) -- (7,2);
\draw [red,thick,{stealth-stealth}](6,3) -- (7,3);
\draw [red,thick,{stealth-stealth}](6,4) -- (7,4);

\draw [red,thick,{stealth-stealth}](7,1) -- (8,1);
\draw [red,thick,{stealth-stealth}](7,2) -- (8,2);
\draw [red,thick,{stealth-stealth}](7,3) -- (8,3);
\draw [red,thick,{stealth-stealth}](7,4) -- (8,4);

\draw [red,thick,{stealth-stealth}](8,1) -- (9,1);
\draw [red,thick,{stealth-stealth}](8,2) -- (9,2);
\draw [-stealth](8,3) -- (9,3);
\draw [red,thick,{stealth-stealth}](8,4) -- (9,4);

\draw [red,thick,{stealth-stealth}](6,1) -- (6,2);
\draw [red,thick,{stealth-stealth}](6,2) -- (6,3);
\draw [red,thick,{stealth-stealth}](6,3) -- (6,4);

\draw [red,thick,{stealth-stealth}](7,1) -- (7,2);
\draw [red,thick,{stealth-stealth}](7,2) -- (7,3);
\draw [-stealth](7,3) -- (7,4);

\draw [stealth-](8,1) -- (8,2);
\draw [red,thick,{stealth-stealth}](8,2) -- (8,3);
\draw [red,thick,{stealth-stealth}](8,3) -- (8,4);

\draw [red,thick,{stealth-stealth}](9,1) -- (9,2);
\draw [red,thick,{stealth-stealth}](9,2) -- (9,3);
\draw [red,thick,{stealth-stealth}](9,3) -- (9,4);

\end{tikzpicture}
\caption{\label{derangementpoly}Elementary flow grid for a 5-dimensional face with 9 vertices and 12 facets.}
\end{figure}

\section{Combinatorial types of low-dimensional faces}

Paffenholz  \cite{Paffenholz1}, by computer enumeration, found eleven 4-dimensional combinatorial types of faces of $B_n$. A similar study of $ASMs$ completed by hand for small values of $n$ is summarized in Table \ref{2level}.  In Table \ref{2level}, Type is the combinatorial type of face, $d$ is the dimension, $v$ is the number of vertices, $f$ is the number of facets.  In addition, $LS$ represents line segment facets, $Sq$ represents square facets, $T$ represents triangles, $S$ represents tetrahedral facets,  $Py$ represents pyramid with a square base facets, $Pr$ represents triangular prism facets, $O$ is octahedral facets, and $C$ is 3-cube facets.  The column ``line" lists the line number of the polytope as found in the table on page 196 of \cite{semidefinite4D}.  The last column gives an example of an elementary flow grid contained in $C_n$ for the smallest possible value of $n$. Reference {\cite{semidefinite4D}} has Schlegel diagrams for many four-dimensional polytopes in a class that strictly contains the 2-level polytopes.

\begin{table}[h]
\begin{Large}
    \centering
\caption{Low dimensional combinatorial types of faces of $ASM_n$.}
\label{2level}
\begin{tabular}{ | m{3cm}| m{0.5cm} | m{0.5cm}| m{0.5cm} | m{3cm}| m{0.6cm} | m{2cm} |}  
 \hline
 Type & $d$ & $v$ & $f$ & Facets & line & EFG  \\ [0.5ex] 
 \hline\hline
 square & 2 & 4 & 4 & 4$LS$ & - &

 \begin{tikzpicture} [scale=0.4]
\draw [stealth-](0.5,1) -- (1,1);
\draw [stealth-](0.5,2) -- (1,2);
\draw [stealth-](0.5,3) -- (1,3);

\draw [-stealth](3,1) -- (3.5,1);
\draw [-stealth](3,2) -- (3.5,2);
\draw [-stealth](3,3) -- (3.5,3);

\draw [-stealth](1,3) -- (1,3.5);
\draw [-stealth](2,3) -- (2,3.5);
\draw [-stealth](3,3) -- (3,3.5);

\draw [-stealth](1,1) -- (1,0.5);
\draw [-stealth](2,1) -- (2,0.5);
\draw [-stealth](3,1) -- (3,0.5);

\draw [red,thick,{stealth-stealth}](1,1) -- (2,1);
\draw [red,thick,{stealth-stealth}](1,2) -- (2,2);
\draw [stealth-](1,3) -- (2,3);

\draw [-stealth](2,1) -- (3,1);
\draw [red,thick,{stealth-stealth}](2,2) -- (3,2);
\draw [red,thick,{stealth-stealth}](2,3) -- (3,3);

\draw [red,thick,{stealth-stealth}](1,1) -- (1,2);
\draw [-stealth](1,2) -- (1,3);

\draw [red,thick,{stealth-stealth}](2,1) -- (2,2);
\draw [red,thick,{stealth-stealth}](2,2) -- (2,3);

\draw [stealth-](3,1) -- (3,2);
\draw [red,thick,{stealth-stealth}](3,2) -- (3,3);

\end{tikzpicture} \\
 \hline
 triangle & 2 & 3 & 3 & 3$LS$ & - &

 \begin{tikzpicture} [scale=0.4]
\draw [stealth-](0.5,1) -- (1,1);
\draw [stealth-](0.5,2) -- (1,2);
\draw [stealth-](0.5,3) -- (1,3);

\draw [-stealth](3,1) -- (3.5,1);
\draw [-stealth](3,2) -- (3.5,2);
\draw [-stealth](3,3) -- (3.5,3);

\draw [-stealth](1,3) -- (1,3.5);
\draw [-stealth](2,3) -- (2,3.5);
\draw [-stealth](3,3) -- (3,3.5);

\draw [-stealth](1,1) -- (1,0.5);
\draw [-stealth](2,1) -- (2,0.5);
\draw [-stealth](3,1) -- (3,0.5);

\draw [stealth-](1,1) -- (2,1);
\draw [red,thick,{stealth-stealth}](1,2) -- (2,2);
\draw [red,thick,{stealth-stealth}](1,3) -- (2,3);

\draw [red,thick,{stealth-stealth}](2,1) -- (3,1);
\draw [red,thick,{stealth-stealth}](2,2) -- (3,2);
\draw [red,thick,{stealth-stealth}](2,3) -- (3,3);

\draw [stealth-](1,1) -- (1,2);
\draw [red,thick,{stealth-stealth}](1,2) -- (1,3);

\draw [red,thick,{stealth-stealth}](2,1) -- (2,2);
\draw [-stealth](2,2) -- (2,3);

\draw [red,thick,{stealth-stealth}](3,1) -- (3,2);
\draw [red,thick,{stealth-stealth}](3,2) -- (3,3);

\end{tikzpicture}  \\
 
 \hline
 
 pyramid with a square base & 3 & 5 & 5 & 1$Sq$ 4$T$ & - &
 \begin{tikzpicture} [scale=0.4]
\draw [stealth-](0.5,1) -- (1,1);
\draw [stealth-](0.5,2) -- (1,2);
\draw [stealth-](0.5,3) -- (1,3);

\draw [-stealth](3,1) -- (3.5,1);
\draw [-stealth](3,2) -- (3.5,2);
\draw [-stealth](3,3) -- (3.5,3);

\draw [-stealth](1,3) -- (1,3.5);
\draw [-stealth](2,3) -- (2,3.5);
\draw [-stealth](3,3) -- (3,3.5);

\draw [-stealth](1,1) -- (1,0.5);
\draw [-stealth](2,1) -- (2,0.5);
\draw [-stealth](3,1) -- (3,0.5);

\draw [red,thick,{stealth-stealth}](1,1) -- (2,1);
\draw [red,thick,{stealth-stealth}](1,2) -- (2,2);
\draw [red,thick,{stealth-stealth}](1,3) -- (2,3);

\draw [red,thick,{stealth-stealth}](2,1) -- (3,1);
\draw [red,thick,{stealth-stealth}](2,2) -- (3,2);
\draw [-stealth](2,3) -- (3,3);

\draw [red,thick,{stealth-stealth}](1,1) -- (1,2);
\draw [red,thick,{stealth-stealth}](1,2) -- (1,3);

\draw [red,thick,{stealth-stealth}](2,1) -- (2,2);
\draw [red,thick,{stealth-stealth}](2,2) -- (2,3);

\draw [red,thick,{stealth-stealth}](3,1) -- (3,2);
\draw [-stealth](3,2) -- (3,3);

\end{tikzpicture}
 \\
 \hline
 tetrahedron & 3 & 4 & 4 & 4$T$ & - &
 \begin{tikzpicture} [scale=0.4]
\draw [stealth-](0.5,1) -- (1,1);
\draw [stealth-](0.5,2) -- (1,2);
\draw [stealth-](0.5,3) -- (1,3);

\draw [-stealth](3,1) -- (3.5,1);
\draw [-stealth](3,2) -- (3.5,2);
\draw [-stealth](3,3) -- (3.5,3);

\draw [-stealth](1,3) -- (1,3.5);
\draw [-stealth](2,3) -- (2,3.5);
\draw [-stealth](3,3) -- (3,3.5);

\draw [-stealth](1,1) -- (1,0.5);
\draw [-stealth](2,1) -- (2,0.5);
\draw [-stealth](3,1) -- (3,0.5);

\draw [red,thick,{stealth-stealth}](1,1) -- (2,1);
\draw [red,thick,{stealth-stealth}](1,2) -- (2,2);
\draw [red,thick,{stealth-stealth}](1,3) -- (2,3);

\draw [red,thick,{stealth-stealth}](2,1) -- (3,1);
\draw [red,thick,{stealth-stealth}](2,2) -- (3,2);
\draw [red,thick,{stealth-stealth}](2,3) -- (3,3);

\draw [red,thick,{stealth-stealth}](1,1) -- (1,2);
\draw [red,thick,{stealth-stealth}](1,2) -- (1,3);

\draw [red,thick,{stealth-stealth}](2,1) -- (2,2);
\draw [-stealth](2,2) -- (2,3);

\draw [red,thick,{stealth-stealth}](3,1) -- (3,2);
\draw [red,thick,{stealth-stealth}](3,2) -- (3,3);

\end{tikzpicture}
 \\ 
 \hline
 3-cube & 3 & 8 & 6 & 6$Sq$ & - &
 
\begin{tikzpicture} [scale=0.4]
\draw [stealth-](0.5,1) -- (1,1);
\draw [stealth-](0.5,2) -- (1,2);
\draw [stealth-](0.5,3) -- (1,3);
\draw [stealth-](0.5,4) -- (1,4);

\draw [-stealth](4,1) -- (4.5,1);
\draw [-stealth](4,2) -- (4.5,2);
\draw [-stealth](4,3) -- (4.5,3);
\draw [-stealth](4,4) -- (4.5,4);

\draw [-stealth](1,4) -- (1,4.5);
\draw [-stealth](2,4) -- (2,4.5);
\draw [-stealth](3,4) -- (3,4.5);
\draw [-stealth](4,4) -- (4,4.5);

\draw [-stealth](1,1) -- (1,0.5);
\draw [-stealth](2,1) -- (2,0.5);
\draw [-stealth](3,1) -- (3,0.5);
\draw [-stealth](4,1) -- (4,0.5);

\draw [stealth-](1,1) -- (2,1);
\draw [red,thick,{stealth-stealth}](1,2) -- (2,2);
\draw [red,thick,{stealth-stealth}](1,3) -- (2,3);
\draw [stealth-](1,4) -- (2,4);

\draw [red,thick,{stealth-stealth}](2,1) -- (3,1);
\draw [red,thick,{stealth-stealth}](2,2) -- (3,2);
\draw [stealth-](2,3) -- (3,3);
\draw [-stealth](2,4) -- (3,4);

\draw [-stealth](3,1) -- (4,1);
\draw [red,thick,{stealth-stealth}](3,2) -- (4,2);
\draw [red,thick,{stealth-stealth}](3,3) -- (4,3);
\draw [-stealth](3,4) -- (4,4);

\draw [stealth-](1,1) -- (1,2);
\draw [red,thick,{stealth-stealth}](1,2) -- (1,3);
\draw [-stealth](1,3) -- (1,4);

\draw [red,thick,{stealth-stealth}](2,1) -- (2,2);
\draw [red,thick,{stealth-stealth}](2,2) -- (2,3);
\draw [stealth-](2,3) -- (2,4);

\draw [red,thick,{stealth-stealth}](3,1) -- (3,2);
\draw [red,thick,{stealth-stealth}](3,2) -- (3,3);
\draw [-stealth](3,3) -- (3,4);

\draw [stealth-](4,1) -- (4,2);
\draw [red,thick,{stealth-stealth}](4,2) -- (4,3);
\draw [-stealth](4,3) -- (4,4);

\end{tikzpicture} 
 \\
 \hline
 octahedron & 3 & 6 & 8 & 8$T$ & - &
 \begin{tikzpicture} [scale=0.4]
\draw [stealth-](0.5,1) -- (1,1);
\draw [stealth-](0.5,2) -- (1,2);
\draw [stealth-](0.5,3) -- (1,3);
\draw [stealth-](0.5,4) -- (1,4);

\draw [-stealth](4,1) -- (4.5,1);
\draw [-stealth](4,2) -- (4.5,2);
\draw [-stealth](4,3) -- (4.5,3);
\draw [-stealth](4,4) -- (4.5,4);

\draw [-stealth](1,4) -- (1,4.5);
\draw [-stealth](2,4) -- (2,4.5);
\draw [-stealth](3,4) -- (3,4.5);
\draw [-stealth](4,4) -- (4,4.5);

\draw [-stealth](1,1) -- (1,0.5);
\draw [-stealth](2,1) -- (2,0.5);
\draw [-stealth](3,1) -- (3,0.5);
\draw [-stealth](4,1) -- (4,0.5);

\draw [stealth-](1,1) -- (2,1);
\draw [red,thick,{stealth-stealth}](1,2) -- (2,2);
\draw [stealth-](1,3) -- (2,3);
\draw [red,thick,{stealth-stealth}](1,4) -- (2,4);

\draw [red,thick,{stealth-stealth}](2,1) -- (3,1);
\draw [red,thick,{stealth-stealth}](2,2) -- (3,2);
\draw [red,thick,{stealth-stealth}](2,3) -- (3,3);
\draw [red,thick,{stealth-stealth}](2,4) -- (3,4);

\draw [red,thick,{stealth-stealth}](3,1) -- (4,1);
\draw [-stealth](3,2) -- (4,2);
\draw [red,thick,{stealth-stealth}](3,3) -- (4,3);
\draw [-stealth](3,4) -- (4,4);

\draw [stealth-](1,1) -- (1,2);
\draw [red,thick,{stealth-stealth}](1,2) -- (1,3);
\draw [red,thick,{stealth-stealth}](1,3) -- (1,4);

\draw [red,thick,{stealth-stealth}](2,1) -- (2,2);
\draw [red,thick,{stealth-stealth}](2,2) -- (2,3);
\draw [-stealth](2,3) -- (2,4);

\draw [stealth-](3,1) -- (3,2);
\draw [red,thick,{stealth-stealth}](3,2) -- (3,3);
\draw [red,thick,{stealth-stealth}](3,3) -- (3,4);

\draw [red,thick,{stealth-stealth}](4,1) -- (4,2);
\draw [red,thick,{stealth-stealth}](4,2) -- (4,3);
\draw [-stealth](4,3) -- (4,4);

\end{tikzpicture}
 \\ 
 \hline
 triangular prism & 3 & 6 & 5 & 3$Sq$ 2$T$ & - &
 \begin{tikzpicture} [scale=0.4]

\draw [stealth-](5.5,1) -- (6,1);
\draw [stealth-](5.5,2) -- (6,2);
\draw [stealth-](5.5,3) -- (6,3);
\draw [stealth-](5.5,4) -- (6,4);

\draw [-stealth](9,1) -- (9.5,1);
\draw [-stealth](9,2) -- (9.5,2);
\draw [-stealth](9,3) -- (9.5,3);
\draw [-stealth](9,4) -- (9.5,4);

\draw [-stealth](6,4) -- (6,4.5);
\draw [-stealth](7,4) -- (7,4.5);
\draw [-stealth](8,4) -- (8,4.5);
\draw [-stealth](9,4) -- (9,4.5);

\draw [-stealth](6,1) -- (6,0.5);
\draw [-stealth](7,1) -- (7,0.5);
\draw [-stealth](8,1) -- (8,0.5);
\draw [-stealth](9,1) -- (9,0.5);

\draw [stealth-](6,1) -- (7,1);
\draw [stealth-](6,2) -- (7,2);
\draw [red,thick,{stealth-stealth}](6,3) -- (7,3);
\draw [red,thick,{stealth-stealth}](6,4) -- (7,4);

\draw [stealth-](7,1) -- (8,1);
\draw [-stealth](7,2) -- (8,2);
\draw [stealth-](7,3) -- (8,3);
\draw [-stealth](7,4) -- (8,4);

\draw [red,thick,{stealth-stealth}](8,1) -- (9,1);
\draw [red,thick,{stealth-stealth}](8,2) -- (9,2);
\draw [red,thick,{stealth-stealth}](8,3) -- (9,3);
\draw [-stealth](8,4) -- (9,4);

\draw [stealth-](6,1) -- (6,2);
\draw [stealth-](6,2) -- (6,3);
\draw [red,thick,{stealth-stealth}](6,3) -- (6,4);

\draw [stealth-](7,1) -- (7,2);
\draw [-stealth](7,2) -- (7,3);
\draw [red,thick,{stealth-stealth}](7,3) -- (7,4);

\draw [red,thick,{stealth-stealth}](8,1) -- (8,2);
\draw [red,thick,{stealth-stealth}](8,2) -- (8,3);
\draw [-stealth](8,3) -- (8,4);

\draw [red,thick,{stealth-stealth}](9,1) -- (9,2);
\draw [red,thick,{stealth-stealth}](9,2) -- (9,3);
\draw [-stealth](9,3) -- (9,4);

\end{tikzpicture}
 \\ 
 \hline
 4-simplex & 4 & 5 & 5 & 5$S$ & 1 & 
  \begin{tikzpicture} [scale=0.4]

\draw [stealth-](5.5,1) -- (6,1);
\draw [stealth-](5.5,2) -- (6,2);
\draw [stealth-](5.5,3) -- (6,3);
\draw [stealth-](5.5,4) -- (6,4);

\draw [-stealth](9,1) -- (9.5,1);
\draw [-stealth](9,2) -- (9.5,2);
\draw [-stealth](9,3) -- (9.5,3);
\draw [-stealth](9,4) -- (9.5,4);

\draw [-stealth](6,4) -- (6,4.5);
\draw [-stealth](7,4) -- (7,4.5);
\draw [-stealth](8,4) -- (8,4.5);
\draw [-stealth](9,4) -- (9,4.5);

\draw [-stealth](6,1) -- (6,0.5);
\draw [-stealth](7,1) -- (7,0.5);
\draw [-stealth](8,1) -- (8,0.5);
\draw [-stealth](9,1) -- (9,0.5);

\draw [stealth-](6,1) -- (7,1);
\draw [red,thick,{stealth-stealth}](6,2) -- (7,2);
\draw [red,thick,{stealth-stealth}](6,3) -- (7,3);
\draw [red,thick,{stealth-stealth}](6,4) -- (7,4);

\draw [-stealth](7,1) -- (8,1);
\draw [red,thick,{stealth-stealth}](7,2) -- (8,2);
\draw [red,thick,{stealth-stealth}](7,3) -- (8,3);
\draw [red,thick,{stealth-stealth}](7,4) -- (8,4);

\draw [-stealth](8,1) -- (9,1);
\draw [red,thick,{stealth-stealth}](8,2) -- (9,2);
\draw [red,thick,{stealth-stealth}](8,3) -- (9,3);
\draw [red,thick,{stealth-stealth}](8,4) -- (9,4);

\draw [stealth-](6,1) -- (6,2);
\draw [red,thick,{stealth-stealth}](6,2) -- (6,3);
\draw [red,thick,{stealth-stealth}](6,3) -- (6,4);

\draw [-stealth](7,1) -- (7,2);
\draw [-stealth](7,2) -- (7,3);
\draw [red,thick,{stealth-stealth}](7,3) -- (7,4);

\draw [stealth-](8,1) -- (8,2);
\draw [stealth-](8,2) -- (8,3);
\draw [red,thick,{stealth-stealth}](8,3) -- (8,4);

\draw [stealth-](9,1) -- (9,2);
\draw [red,thick,{stealth-stealth}](9,2) -- (9,3);
\draw [red,thick,{stealth-stealth}](9,3) -- (9,4);

\end{tikzpicture}
 \\ [1ex] 
 \hline
 square pyramid pyramid & 4 & 6 & 6 & 4$S$ 2$Py$ & 2 &
 \begin{tikzpicture} [scale=0.4]

\draw [stealth-](5.5,1) -- (6,1);
\draw [stealth-](5.5,2) -- (6,2);
\draw [stealth-](5.5,3) -- (6,3);
\draw [stealth-](5.5,4) -- (6,4);

\draw [-stealth](9,1) -- (9.5,1);
\draw [-stealth](9,2) -- (9.5,2);
\draw [-stealth](9,3) -- (9.5,3);
\draw [-stealth](9,4) -- (9.5,4);

\draw [-stealth](6,4) -- (6,4.5);
\draw [-stealth](7,4) -- (7,4.5);
\draw [-stealth](8,4) -- (8,4.5);
\draw [-stealth](9,4) -- (9,4.5);

\draw [-stealth](6,1) -- (6,0.5);
\draw [-stealth](7,1) -- (7,0.5);
\draw [-stealth](8,1) -- (8,0.5);
\draw [-stealth](9,1) -- (9,0.5);

\draw [red,thick,{stealth-stealth}](6,1) -- (7,1);
\draw [red,thick,{stealth-stealth}](6,2) -- (7,2);
\draw [red,thick,{stealth-stealth}](6,3) -- (7,3);
\draw [stealth-](6,4) -- (7,4);

\draw [red,thick,{stealth-stealth}](7,1) -- (8,1);
\draw [red,thick,{stealth-stealth}](7,2) -- (8,2);
\draw [red,thick,{stealth-stealth}](7,3) -- (8,3);
\draw [-stealth](7,4) -- (8,4);

\draw [-stealth](8,1) -- (9,1);
\draw [stealth-](8,2) -- (9,2);
\draw [-stealth](8,3) -- (9,3);
\draw [-stealth](8,4) -- (9,4);

\draw [red,thick,{stealth-stealth}](6,1) -- (6,2);
\draw [red,thick,{stealth-stealth}](6,2) -- (6,3);
\draw [-stealth](6,3) -- (6,4);

\draw [red,thick,{stealth-stealth}](7,1) -- (7,2);
\draw [red,thick,{stealth-stealth}](7,2) -- (7,3);
\draw [stealth-](7,3) -- (7,4);

\draw [red,thick,{stealth-stealth}](8,1) -- (8,2);
\draw [red,thick,{stealth-stealth}](8,2) -- (8,3);
\draw [-stealth](8,3) -- (8,4);

\draw [stealth-] (9,1) -- (9,2);
\draw [-stealth](9,2) -- (9,3);
\draw [-stealth](9,3) -- (9,4);

\end{tikzpicture}
 \\ 
 \hline
 pyramid over a triangular prism & 4 & 7 & 6 & 2$S$ 3$Py$ 1$Pr$ & 8 & 
\begin{tikzpicture} [scale=0.4]

\draw [stealth-](5.5,1) -- (6,1);
\draw [stealth-](5.5,2) -- (6,2);
\draw [stealth-](5.5,3) -- (6,3);
\draw [stealth-](5.5,4) -- (6,4);

\draw [-stealth](9,1) -- (9.5,1);
\draw [-stealth](9,2) -- (9.5,2);
\draw [-stealth](9,3) -- (9.5,3);
\draw [-stealth](9,4) -- (9.5,4);

\draw [-stealth](6,4) -- (6,4.5);
\draw [-stealth](7,4) -- (7,4.5);
\draw [-stealth](8,4) -- (8,4.5);
\draw [-stealth](9,4) -- (9,4.5);

\draw [-stealth](6,1) -- (6,0.5);
\draw [-stealth](7,1) -- (7,0.5);
\draw [-stealth](8,1) -- (8,0.5);
\draw [-stealth](9,1) -- (9,0.5);

\draw [red,thick,{stealth-stealth}](6,1) -- (7,1);
\draw [red,thick,{stealth-stealth}](6,2) -- (7,2);
\draw [red,thick,{stealth-stealth}](6,3) -- (7,3);
\draw [red,thick,{stealth-stealth}](6,4) -- (7,4);

\draw [-stealth](7,1) -- (8,1);
\draw [stealth-](7,2) -- (8,2);
\draw [red,thick,{stealth-stealth}](7,3) -- (8,3);
\draw [red,thick,{stealth-stealth}](7,4) -- (8,4);

\draw [-stealth](8,1) -- (9,1);
\draw [stealth-](8,2) -- (9,2);
\draw [-stealth](8,3) -- (9,3);
\draw [-stealth](8,4) -- (9,4);

\draw [red,thick,{stealth-stealth}](6,1) -- (6,2);
\draw [red,thick,{stealth-stealth}](6,2) -- (6,3);
\draw [red,thick,{stealth-stealth}](6,3) -- (6,4);

\draw [red,thick,{stealth-stealth}](7,1) -- (7,2);
\draw [red,thick,{stealth-stealth}](7,2) -- (7,3);
\draw [red,thick,{stealth-stealth}](7,3) -- (7,4);

\draw [stealth-](8,1) -- (8,2);
\draw [stealth-](8,2) -- (8,3);
\draw [red,thick,{stealth-stealth}](8,3) -- (8,4);

\draw [stealth-] (9,1) -- (9,2);
\draw [-stealth](9,2) -- (9,3);
\draw [-stealth](9,3) -- (9,4);

\end{tikzpicture}
 \\ 
 \hline
 $ASM_3$ & 4 & 7 & 8 & 4$S$ 4$Py$ & 11 &
  \begin{tikzpicture} [scale=0.4]
\draw [stealth-](0.5,1) -- (1,1);
\draw [stealth-](0.5,2) -- (1,2);
\draw [stealth-](0.5,3) -- (1,3);

\draw [-stealth](3,1) -- (3.5,1);
\draw [-stealth](3,2) -- (3.5,2);
\draw [-stealth](3,3) -- (3.5,3);

\draw [-stealth](1,3) -- (1,3.5);
\draw [-stealth](2,3) -- (2,3.5);
\draw [-stealth](3,3) -- (3,3.5);

\draw [-stealth](1,1) -- (1,0.5);
\draw [-stealth](2,1) -- (2,0.5);
\draw [-stealth](3,1) -- (3,0.5);

\draw [red,thick,{stealth-stealth}](1,1) -- (2,1);
\draw [red,thick,{stealth-stealth}](1,2) -- (2,2);
\draw [red,thick,{stealth-stealth}](1,3) -- (2,3);

\draw [red,thick,{stealth-stealth}](2,1) -- (3,1);
\draw [red,thick,{stealth-stealth}](2,2) -- (3,2);
\draw [red,thick,{stealth-stealth}](2,3) -- (3,3);

\draw [red,thick,{stealth-stealth}](1,1) -- (1,2);
\draw [red,thick,{stealth-stealth}](1,2) -- (1,3);

\draw [red,thick,{stealth-stealth}](2,1) -- (2,2);
\draw [red,thick,{stealth-stealth}](2,2) -- (2,3);

\draw [red,thick,{stealth-stealth}](3,1) -- (3,2);
\draw [red,thick,{stealth-stealth}](3,2) -- (3,3);

\end{tikzpicture} 
 \\ 
 \hline
 tetrahedron prism & 4 & 8 & 6 & 2$S$ 4$Pr$ & 4 & 
   \begin{tikzpicture} [scale=0.4]

\draw [stealth-](5.5,1) -- (6,1);
\draw [stealth-](5.5,2) -- (6,2);
\draw [stealth-](5.5,3) -- (6,3);
\draw [stealth-](5.5,4) -- (6,4);

\draw [-stealth](9,1) -- (9.5,1);
\draw [-stealth](9,2) -- (9.5,2);
\draw [-stealth](9,3) -- (9.5,3);
\draw [-stealth](9,4) -- (9.5,4);

\draw [-stealth](6,4) -- (6,4.5);
\draw [-stealth](7,4) -- (7,4.5);
\draw [-stealth](8,4) -- (8,4.5);
\draw [-stealth](9,4) -- (9,4.5);

\draw [-stealth](6,1) -- (6,0.5);
\draw [-stealth](7,1) -- (7,0.5);
\draw [-stealth](8,1) -- (8,0.5);
\draw [-stealth](9,1) -- (9,0.5);

\draw [red,thick,{stealth-stealth}](6,1) -- (7,1);
\draw [stealth-](6,2) -- (7,2);
\draw [red,thick,{stealth-stealth}](6,3) -- (7,3);
\draw [stealth-](6,4) -- (7,4);

\draw [red,thick,{stealth-stealth}](7,1) -- (8,1);
\draw [red,thick,{stealth-stealth}](7,2) -- (8,2);
\draw [red,thick,{stealth-stealth}](7,3) -- (8,3);
\draw [stealth-](7,4) -- (8,4);

\draw [-stealth](8,1) -- (9,1);
\draw [-stealth](8,2) -- (9,2);
\draw [red,thick,{stealth-stealth}](8,3) -- (9,3);
\draw [red,thick,{stealth-stealth}](8,4) -- (9,4);

\draw [red,thick,{stealth-stealth}](6,1) -- (6,2);
\draw [red,thick,{stealth-stealth}](6,2) -- (6,3);
\draw [-stealth](6,3) -- (6,4);

\draw [red,thick,{stealth-stealth}](7,1) -- (7,2);
\draw [red,thick,{stealth-stealth}](7,2) -- (7,3);
\draw [-stealth](7,3) -- (7,4);

\draw [red,thick,{stealth-stealth}](8,1) -- (8,2);
\draw [red,thick,{stealth-stealth}](8,2) -- (8,3);
\draw [red,thick,{stealth-stealth}](8,3) -- (8,4);

\draw [stealth-] (9,1) -- (9,2);
\draw [stealth-](9,2) -- (9,3);
\draw [red,thick,{stealth-stealth}](9,3) -- (9,4);

\end{tikzpicture}
 \\ [1ex] 
 \hline

\end{tabular}
\end{Large}
\end{table}

\begin{table}[h]
\begin{Large}

\begin{tabular}{ | m{3cm}| m{0.5cm} | m{0.5cm}| m{0.5cm} | m{3cm}| m{0.6cm} | m{2cm} |}  
 \hline
 Type & $d$ & $v$ & $f$ & Facets & line & EFG  \\ [0.5ex] 
 \hline\hline
 
 polytope Y & 4 & 8 & 7 & 1$S$ 4$Py$ 2$Pr$ & 10 & 
\begin{tikzpicture} [scale=0.4]

\draw [stealth-](5.5,1) -- (6,1);
\draw [stealth-](5.5,2) -- (6,2);
\draw [stealth-](5.5,3) -- (6,3);
\draw [stealth-](5.5,4) -- (6,4);

\draw [-stealth](9,1) -- (9.5,1);
\draw [-stealth](9,2) -- (9.5,2);
\draw [-stealth](9,3) -- (9.5,3);
\draw [-stealth](9,4) -- (9.5,4);

\draw [-stealth](6,4) -- (6,4.5);
\draw [-stealth](7,4) -- (7,4.5);
\draw [-stealth](8,4) -- (8,4.5);
\draw [-stealth](9,4) -- (9,4.5);

\draw [-stealth](6,1) -- (6,0.5);
\draw [-stealth](7,1) -- (7,0.5);
\draw [-stealth](8,1) -- (8,0.5);
\draw [-stealth](9,1) -- (9,0.5);

\draw [red,thick,{stealth-stealth}](6,1) -- (7,1);
\draw [red,thick,{stealth-stealth}](6,2) -- (7,2);
\draw [red,thick,{stealth-stealth}](6,3) -- (7,3);
\draw [stealth-](6,4) -- (7,4);

\draw [-stealth](7,1) -- (8,1);
\draw [red,thick,{stealth-stealth}](7,2) -- (8,2);
\draw [red,thick,{stealth-stealth}](7,3) -- (8,3);
\draw [red,thick,{stealth-stealth}](7,4) -- (8,4);

\draw [-stealth](8,1) -- (9,1);
\draw [-stealth](8,2) -- (9,2);
\draw [stealth-](8,3) -- (9,3);
\draw [-stealth](8,4) -- (9,4);

\draw [red,thick,{stealth-stealth}](6,1) -- (6,2);
\draw [red,thick,{stealth-stealth}](6,2) -- (6,3);
\draw [-stealth](6,3) -- (6,4);

\draw [red,thick,{stealth-stealth}](7,1) -- (7,2);
\draw [red,thick,{stealth-stealth}](7,2) -- (7,3);
\draw [red,thick,{stealth-stealth}](7,3) -- (7,4);

\draw [stealth-](8,1) -- (8,2);
\draw [red,thick,{stealth-stealth}](8,2) -- (8,3);
\draw [red,thick,{stealth-stealth}](8,3) -- (8,4);

\draw [stealth-] (9,1) -- (9,2);
\draw [stealth-](9,2) -- (9,3);
\draw [-stealth](9,3) -- (9,4);

\end{tikzpicture} 
 \\ 
 \hline
 polytope X & 4 & 8 & 9 & 4$S$ 4$Py$ 1$O$ & 19 & 
 \begin{tikzpicture} [scale=0.4]

\draw [stealth-](5.5,1) -- (6,1);
\draw [stealth-](5.5,2) -- (6,2);
\draw [stealth-](5.5,3) -- (6,3);
\draw [stealth-](5.5,4) -- (6,4);

\draw [-stealth](9,1) -- (9.5,1);
\draw [-stealth](9,2) -- (9.5,2);
\draw [-stealth](9,3) -- (9.5,3);
\draw [-stealth](9,4) -- (9.5,4);

\draw [-stealth](6,4) -- (6,4.5);
\draw [-stealth](7,4) -- (7,4.5);
\draw [-stealth](8,4) -- (8,4.5);
\draw [-stealth](9,4) -- (9,4.5);

\draw [-stealth](6,1) -- (6,0.5);
\draw [-stealth](7,1) -- (7,0.5);
\draw [-stealth](8,1) -- (8,0.5);
\draw [-stealth](9,1) -- (9,0.5);

\draw [red,thick,{stealth-stealth}](6,1) -- (7,1);
\draw [red,thick,{stealth-stealth}](6,2) -- (7,2);
\draw [stealth-](6,3) -- (7,3);
\draw [red,thick,{stealth-stealth}](6,4) -- (7,4);

\draw [red,thick,{stealth-stealth}](7,1) -- (8,1);
\draw [red,thick,{stealth-stealth}](7,2) -- (8,2);
\draw [red,thick,{stealth-stealth}](7,3) -- (8,3);
\draw [red,thick,{stealth-stealth}](7,4) -- (8,4);

\draw [red,thick,{stealth-stealth}](8,1) -- (9,1);
\draw [-stealth](8,2) -- (9,2);
\draw [red,thick,{stealth-stealth}](8,3) -- (9,3);
\draw [-stealth](8,4) -- (9,4);

\draw [red,thick,{stealth-stealth}](6,1) -- (6,2);
\draw [red,thick,{stealth-stealth}](6,2) -- (6,3);
\draw [red,thick,{stealth-stealth}](6,3) -- (6,4);

\draw [red,thick,{stealth-stealth}](7,1) -- (7,2);
\draw [red,thick,{stealth-stealth}](7,2) -- (7,3);
\draw [-stealth](7,3) -- (7,4);

\draw [stealth-](8,1) -- (8,2);
\draw [red,thick,{stealth-stealth}](8,2) -- (8,3);
\draw [red,thick,{stealth-stealth}](8,3) -- (8,4);

\draw [red,thick,{stealth-stealth}] (9,1) -- (9,2);
\draw [red,thick,{stealth-stealth}](9,2) -- (9,3);
\draw [-stealth](9,3) -- (9,4);

\end{tikzpicture}
 \\ 
 \hline
 3-3 duoprism & 4 & 9 & 6 & 6$Pr$ & 6 & 
 \begin{tikzpicture} [scale=0.4]

\draw [stealth-](5.5,1) -- (6,1);
\draw [stealth-](5.5,2) -- (6,2);
\draw [stealth-](5.5,3) -- (6,3);
\draw [stealth-](5.5,4) -- (6,4);

\draw [-stealth](9,1) -- (9.5,1);
\draw [-stealth](9,2) -- (9.5,2);
\draw [-stealth](9,3) -- (9.5,3);
\draw [-stealth](9,4) -- (9.5,4);

\draw [-stealth](6,4) -- (6,4.5);
\draw [-stealth](7,4) -- (7,4.5);
\draw [-stealth](8,4) -- (8,4.5);
\draw [-stealth](9,4) -- (9,4.5);

\draw [-stealth](6,1) -- (6,0.5);
\draw [-stealth](7,1) -- (7,0.5);
\draw [-stealth](8,1) -- (8,0.5);
\draw [-stealth](9,1) -- (9,0.5);

\draw [stealth-](6,1) -- (7,1);
\draw [stealth-](6,2) -- (7,2);
\draw [red,thick,{stealth-stealth}](6,3) -- (7,3);
\draw [red,thick,{stealth-stealth}](6,4) -- (7,4);

\draw [stealth-](7,1) -- (8,1);
\draw [-stealth](7,2) -- (8,2);
\draw [red,thick,{stealth-stealth}](7,3) -- (8,3);
\draw [red,thick,{stealth-stealth}](7,4) -- (8,4);

\draw [red,thick,{stealth-stealth}](8,1) -- (9,1);
\draw [red,thick,{stealth-stealth}](8,2) -- (9,2);
\draw [red,thick,{stealth-stealth}](8,3) -- (9,3);
\draw [-stealth](8,4) -- (9,4);

\draw [stealth-](6,1) -- (6,2);
\draw [stealth-](6,2) -- (6,3);
\draw [red,thick,{stealth-stealth}](6,3) -- (6,4);

\draw [stealth-](7,1) -- (7,2);
\draw [-stealth](7,2) -- (7,3);
\draw [red,thick,{stealth-stealth}](7,3) -- (7,4);

\draw [red,thick,{stealth-stealth}](8,1) -- (8,2);
\draw [red,thick,{stealth-stealth}](8,2) -- (8,3);
\draw [red,thick,{stealth-stealth}](8,3) -- (8,4);

\draw [red,thick,{stealth-stealth}](9,1) -- (9,2);
\draw [red,thick,{stealth-stealth}](9,2) -- (9,3);
\draw [-stealth](9,3) -- (9,4);

\end{tikzpicture}
 \\ [1ex] 
 \hline 
 cubic pyramid & 4 & 9 & 7 & 6$Py$ 1$C$ & 16 & 
 \begin{tikzpicture} [scale=0.4]

\draw [stealth-](5.5,1) -- (6,1);
\draw [stealth-](5.5,2) -- (6,2);
\draw [stealth-](5.5,3) -- (6,3);
\draw [stealth-](5.5,4) -- (6,4);

\draw [-stealth](9,1) -- (9.5,1);
\draw [-stealth](9,2) -- (9.5,2);
\draw [-stealth](9,3) -- (9.5,3);
\draw [-stealth](9,4) -- (9.5,4);

\draw [-stealth](6,4) -- (6,4.5);
\draw [-stealth](7,4) -- (7,4.5);
\draw [-stealth](8,4) -- (8,4.5);
\draw [-stealth](9,4) -- (9,4.5);

\draw [-stealth](6,1) -- (6,0.5);
\draw [-stealth](7,1) -- (7,0.5);
\draw [-stealth](8,1) -- (8,0.5);
\draw [-stealth](9,1) -- (9,0.5);

\draw [stealth-](6,1) -- (7,1);
\draw [red,thick,{stealth-stealth}](6,2) -- (7,2);
\draw [red,thick,{stealth-stealth}](6,3) -- (7,3);
\draw [stealth-](6,4) -- (7,4);

\draw [-stealth](7,1) -- (8,1);
\draw [red,thick,{stealth-stealth}](7,2) -- (8,2);
\draw [red,thick,{stealth-stealth}](7,3) -- (8,3);
\draw [red,thick,{stealth-stealth}](7,4) -- (8,4);

\draw [-stealth](8,1) -- (9,1);
\draw [red,thick,{stealth-stealth}](8,2) -- (9,2);
\draw [red,thick,{stealth-stealth}](8,3) -- (9,3);
\draw [-stealth](8,4) -- (9,4);

\draw [stealth-](6,1) -- (6,2);
\draw [red,thick,{stealth-stealth}](6,2) -- (6,3);
\draw [-stealth](6,3) -- (6,4);

\draw [-stealth](7,1) -- (7,2);
\draw [red,thick,{stealth-stealth}](7,2) -- (7,3);
\draw [red,thick,{stealth-stealth}](7,3) -- (7,4);

\draw [stealth-](8,1) -- (8,2);
\draw [red,thick,{stealth-stealth}](8,2) -- (8,3);
\draw [red,thick,{stealth-stealth}](8,3) -- (8,4);

\draw [stealth-] (9,1) -- (9,2);
\draw [red,thick,{stealth-stealth}](9,2) -- (9,3);
\draw [-stealth](9,3) -- (9,4);

\end{tikzpicture}
 \\ [1ex] 
 \hline
 polytope Z & 4 & 9 & 9 & 3$S$ 3$Py$ 1$Pr$ 2$O$ & 23 &
 \begin{tikzpicture} [scale=0.4]
\draw [stealth-](0.5,1) -- (1,1);
\draw [stealth-](0.5,2) -- (1,2);
\draw [stealth-](0.5,3) -- (1,3);
\draw [stealth-](0.5,4) -- (1,4);
\draw [stealth-](0.5,5) -- (1,5);

\draw [-stealth](5,1) -- (5.5,1);
\draw [-stealth](5,2) -- (5.5,2);
\draw [-stealth](5,3) -- (5.5,3);
\draw [-stealth](5,4) -- (5.5,4);
\draw [-stealth](5,5) -- (5.5,5);

\draw [-stealth](1,5) -- (1,5.5);
\draw [-stealth](2,5) -- (2,5.5);
\draw [-stealth](3,5) -- (3,5.5);
\draw [-stealth](4,5) -- (4,5.5);
\draw [-stealth](5,5) -- (5,5.5);

\draw [-stealth](1,1) -- (1,0.5);
\draw [-stealth](2,1) -- (2,0.5);
\draw [-stealth](3,1) -- (3,0.5);
\draw [-stealth](4,1) -- (4,0.5);
\draw [-stealth](5,1) -- (5,0.5);

\draw [stealth-](1,1) -- (2,1);
\draw [stealth-](1,2) -- (2,2);
\draw [red,thick,{stealth-stealth}](1,3) -- (2,3);
\draw [stealth-](1,4) -- (2,4);
\draw [red,thick,{stealth-stealth}](1,5) -- (2,5);

\draw [red,thick,{stealth-stealth}](2,1) -- (3,1);
\draw [stealth-](2,2) -- (3,2);
\draw [red,thick,{stealth-stealth}](2,3) -- (3,3);
\draw [red,thick,{stealth-stealth}](2,4) -- (3,4);
\draw [red,thick,{stealth-stealth}](2,5) -- (3,5);

\draw [red,thick,{stealth-stealth}](3,1) -- (4,1);
\draw [stealth-](3,2) -- (4,2);
\draw [-stealth](3,3) -- (4,3);
\draw [red,thick,{stealth-stealth}](3,4) -- (4,4);
\draw [-stealth](3,5) -- (4,5);

\draw [red,thick,{stealth-stealth}](4,1) -- (5,1);
\draw [red,thick,{stealth-stealth}](4,2) -- (5,2);
\draw [-stealth](4,3) -- (5,3);
\draw [red,thick,{stealth-stealth}](4,4) -- (5,4);
\draw [-stealth](4,5) -- (5,5);

\draw [stealth-](1,1) -- (1,2);
\draw [stealth-](1,2) -- (1,3);
\draw [red,thick,-stealth](1,3) -- (1,4);
\draw [red,thick,-stealth](1,4) -- (1,5);

\draw [red,thick,{stealth-stealth}](2,1) -- (2,2);
\draw [red,thick,{stealth-stealth}](2,2) -- (2,3);
\draw [red,thick,-stealth](2,3) -- (2,4);
\draw [-stealth](2,4) -- (2,5);

\draw [stealth-](3,1) -- (3,2);
\draw [stealth-](3,2) -- (3,3);
\draw [red,thick,{stealth-stealth}](3,3) -- (3,4);
\draw [red,thick,{stealth-stealth}](3,4) -- (3,5);

\draw [red,thick,{stealth-stealth}](4,1) -- (4,2);
\draw [-stealth](4,2) -- (4,3);
\draw [-stealth](4,3) -- (4,4);
\draw [-stealth](4,4) -- (4,5);

\draw [red,thick,{stealth-stealth}](5,1) -- (5,2);
\draw [red,thick,{stealth-stealth}](5,2) -- (5,3);
\draw [red,thick,{stealth-stealth}](5,3) -- (5,4);
\draw [-stealth](5,4) -- (5,5);
\end{tikzpicture}
 \\ 
 \hline
 square pyramid prism & 4 & 10 & 7 & 2$Py$ 4$Pr$ 1$C$ & 20 &
  \begin{tikzpicture} [scale=0.4]

\draw [stealth-](5.5,1) -- (6,1);
\draw [stealth-](5.5,2) -- (6,2);
\draw [stealth-](5.5,3) -- (6,3);
\draw [stealth-](5.5,4) -- (6,4);

\draw [-stealth](9,1) -- (9.5,1);
\draw [-stealth](9,2) -- (9.5,2);
\draw [-stealth](9,3) -- (9.5,3);
\draw [-stealth](9,4) -- (9.5,4);

\draw [-stealth](6,4) -- (6,4.5);
\draw [-stealth](7,4) -- (7,4.5);
\draw [-stealth](8,4) -- (8,4.5);
\draw [-stealth](9,4) -- (9,4.5);

\draw [-stealth](6,1) -- (6,0.5);
\draw [-stealth](7,1) -- (7,0.5);
\draw [-stealth](8,1) -- (8,0.5);
\draw [-stealth](9,1) -- (9,0.5);

\draw [stealth-](6,1) -- (7,1);
\draw [stealth-](6,2) -- (7,2);
\draw [red,thick,{stealth-stealth}](6,3) -- (7,3);
\draw [red,thick,{stealth-stealth}](6,4) -- (7,4);

\draw [red,thick,{stealth-stealth}](7,1) -- (8,1);
\draw [red,thick,{stealth-stealth}](7,2) -- (8,2);
\draw [red,thick,{stealth-stealth}](7,3) -- (8,3);
\draw [-stealth](7,4) -- (8,4);

\draw [-stealth](8,1) -- (9,1);
\draw [red,thick,{stealth-stealth}](8,2) -- (9,2);
\draw [red,thick,{stealth-stealth}](8,3) -- (9,3);
\draw [-stealth](8,4) -- (9,4);

\draw [stealth-](6,1) -- (6,2);
\draw [stealth-](6,2) -- (6,3);
\draw [red,thick,{stealth-stealth}](6,3) -- (6,4);

\draw [red,thick,{stealth-stealth}](7,1) -- (7,2);
\draw [red,thick,{stealth-stealth}](7,2) -- (7,3);
\draw [red,thick,{stealth-stealth}](7,3) -- (7,4);

\draw [red,thick,{stealth-stealth}](8,1) -- (8,2);
\draw [red,thick,{stealth-stealth}](8,2) -- (8,3);
\draw [-stealth](8,3) -- (8,4);

\draw [stealth-] (9,1) -- (9,2);
\draw [red,thick,{stealth-stealth}](9,2) -- (9,3);
\draw [-stealth](9,3) -- (9,4);

\end{tikzpicture}
 \\
 \hline 
 cubical bipyramid & 4 & 10 & 12 & 12$Py$ & 26 &
 \begin{tikzpicture} [scale=0.4]

\draw [stealth-](5.5,1) -- (6,1);
\draw [stealth-](5.5,2) -- (6,2);
\draw [stealth-](5.5,3) -- (6,3);
\draw [stealth-](5.5,4) -- (6,4);

\draw [-stealth](9,1) -- (9.5,1);
\draw [-stealth](9,2) -- (9.5,2);
\draw [-stealth](9,3) -- (9.5,3);
\draw [-stealth](9,4) -- (9.5,4);

\draw [-stealth](6,4) -- (6,4.5);
\draw [-stealth](7,4) -- (7,4.5);
\draw [-stealth](8,4) -- (8,4.5);
\draw [-stealth](9,4) -- (9,4.5);

\draw [-stealth](6,1) -- (6,0.5);
\draw [-stealth](7,1) -- (7,0.5);
\draw [-stealth](8,1) -- (8,0.5);
\draw [-stealth](9,1) -- (9,0.5);

\draw [stealth-](6,1) -- (7,1);
\draw [red,thick,{stealth-stealth}](6,2) -- (7,2);
\draw [stealth-](6,3) -- (7,3);
\draw [red,thick,{stealth-stealth}](6,4) -- (7,4);

\draw [red,thick,{stealth-stealth}](7,1) -- (8,1);
\draw [red,thick,{stealth-stealth}](7,2) -- (8,2);
\draw [red,thick,{stealth-stealth}](7,3) -- (8,3);
\draw [red,thick,{stealth-stealth}](7,4) -- (8,4);

\draw [-stealth](8,1) -- (9,1);
\draw [red,thick,{stealth-stealth}](8,2) -- (9,2);
\draw [red,thick,{stealth-stealth}](8,3) -- (9,3);
\draw [-stealth](8,4) -- (9,4);

\draw [stealth-](6,1) -- (6,2);
\draw [red,thick,{stealth-stealth}](6,2) -- (6,3);
\draw [red,thick,{stealth-stealth}](6,3) -- (6,4);

\draw [red,thick,{stealth-stealth}](7,1) -- (7,2);
\draw [red,thick,{stealth-stealth}](7,2) -- (7,3);
\draw [-stealth](7,3) -- (7,4);

\draw [red,thick,{stealth-stealth}](8,1) -- (8,2);
\draw [red,thick,{stealth-stealth}](8,2) -- (8,3);
\draw [red,thick,{stealth-stealth}](8,3) -- (8,4);

\draw [stealth-](9,1) -- (9,2);
\draw [red,thick,{stealth-stealth}](9,2) -- (9,3);
\draw [-stealth](9,3) -- (9,4);

\end{tikzpicture}
 \\ [1ex] 
 \hline
 3-4 duoprism & 4 & 12 & 7 & 4$Pr$ 3$C$ & 28 &
  \begin{tikzpicture} [scale=0.4]

\draw [stealth-](5.5,1) -- (6,1);
\draw [stealth-](5.5,2) -- (6,2);
\draw [stealth-](5.5,3) -- (6,3);
\draw [stealth-](5.5,4) -- (6,4);

\draw [-stealth](9,1) -- (9.5,1);
\draw [-stealth](9,2) -- (9.5,2);
\draw [-stealth](9,3) -- (9.5,3);
\draw [-stealth](9,4) -- (9.5,4);

\draw [-stealth](6,4) -- (6,4.5);
\draw [-stealth](7,4) -- (7,4.5);
\draw [-stealth](8,4) -- (8,4.5);
\draw [-stealth](9,4) -- (9,4.5);

\draw [-stealth](6,1) -- (6,0.5);
\draw [-stealth](7,1) -- (7,0.5);
\draw [-stealth](8,1) -- (8,0.5);
\draw [-stealth](9,1) -- (9,0.5);

\draw [stealth-](6,1) -- (7,1);
\draw [red,thick,{stealth-stealth}](6,2) -- (7,2);
\draw [red,thick,{stealth-stealth}](6,3) -- (7,3);
\draw [stealth-](6,4) -- (7,4);

\draw [red,thick,{stealth-stealth}](7,1) -- (8,1);
\draw [red,thick,{stealth-stealth}](7,2) -- (8,2);
\draw [red,thick,{stealth-stealth}](7,3) -- (8,3);
\draw [red,thick,{stealth-stealth}](7,4) -- (8,4);

\draw [red,thick,{stealth-stealth}](8,1) -- (9,1);
\draw [red,thick,{stealth-stealth}](8,2) -- (9,2);
\draw [-stealth](8,3) -- (9,3);
\draw [-stealth](8,4) -- (9,4);

\draw [stealth-](6,1) -- (6,2);
\draw [red,thick,{stealth-stealth}](6,2) -- (6,3);
\draw [-stealth](6,3) -- (6,4);

\draw [red,thick,{stealth-stealth}](7,1) -- (7,2);
\draw [red,thick,{stealth-stealth}](7,2) -- (7,3);
\draw [red,thick,{stealth-stealth}](7,3) -- (7,4);

\draw [red,thick,{stealth-stealth}](8,1) -- (8,2);
\draw [stealth-](8,2) -- (8,3);
\draw [red,thick,{stealth-stealth}](8,3) -- (8,4);

\draw [red,thick,{stealth-stealth}](9,1) -- (9,2);
\draw [-stealth](9,2) -- (9,3);
\draw [-stealth](9,3) -- (9,4);

\end{tikzpicture}
 \\
 \hline
 octahedral prism & 4 & 12 & 10 & 8$Pr$ 2$O$ & 27 &
  \begin{tikzpicture} [scale=0.4]
\draw [stealth-](0.5,1) -- (1,1);
\draw [stealth-](0.5,2) -- (1,2);
\draw [stealth-](0.5,3) -- (1,3);
\draw [stealth-](0.5,4) -- (1,4);
\draw [stealth-](0.5,5) -- (1,5);

\draw [-stealth](5,1) -- (5.5,1);
\draw [-stealth](5,2) -- (5.5,2);
\draw [-stealth](5,3) -- (5.5,3);
\draw [-stealth](5,4) -- (5.5,4);
\draw [-stealth](5,5) -- (5.5,5);

\draw [-stealth](1,5) -- (1,5.5);
\draw [-stealth](2,5) -- (2,5.5);
\draw [-stealth](3,5) -- (3,5.5);
\draw [-stealth](4,5) -- (4,5.5);
\draw [-stealth](5,5) -- (5,5.5);

\draw [-stealth](1,1) -- (1,0.5);
\draw [-stealth](2,1) -- (2,0.5);
\draw [-stealth](3,1) -- (3,0.5);
\draw [-stealth](4,1) -- (4,0.5);
\draw [-stealth](5,1) -- (5,0.5);

\draw [stealth-](1,1) -- (2,1);
\draw [stealth-](1,2) -- (2,2);
\draw [red,thick,{stealth-stealth}](1,3) -- (2,3);
\draw [stealth-](1,4) -- (2,4);
\draw [red,thick,{stealth-stealth}](1,5) -- (2,5);

\draw [red,thick,{stealth-stealth}](2,1) -- (3,1);
\draw [stealth-](2,2) -- (3,2);
\draw [red,thick,{stealth-stealth}](2,3) -- (3,3);
\draw [red,thick,{stealth-stealth}](2,4) -- (3,4);
\draw [red,thick,{stealth-stealth}](2,5) -- (3,5);

\draw [red,thick,{stealth-stealth}](3,1) -- (4,1);
\draw [red,thick,{stealth-stealth}](3,2) -- (4,2);
\draw [red,thick,{stealth-stealth}](3,3) -- (4,3);
\draw [red,thick,{stealth-stealth}](3,4) -- (4,4);
\draw [-stealth](3,5) -- (4,5);

\draw [red,thick,{stealth-stealth}](4,1) -- (5,1);
\draw [-stealth](4,2) -- (5,2);
\draw [-stealth](4,3) -- (5,3);
\draw [red,thick,{stealth-stealth}](4,4) -- (5,4);
\draw [-stealth](4,5) -- (5,5);

\draw [stealth-](1,1) -- (1,2);
\draw [stealth-](1,2) -- (1,3);
\draw [red,thick,{stealth-stealth}](1,3) -- (1,4);
\draw [red,thick,{stealth-stealth}](1,4) -- (1,5);

\draw [red,thick,{stealth-stealth}](2,1) -- (2,2);
\draw [red,thick,{stealth-stealth}](2,2) -- (2,3);
\draw [red,thick,{stealth-stealth}](2,3) -- (2,4);
\draw [-stealth](2,4) -- (2,5);

\draw [stealth-](3,1) -- (3,2);
\draw [red,thick,{stealth-stealth}](3,2) -- (3,3);
\draw [red,thick,{stealth-stealth}](3,3) -- (3,4);
\draw [red,thick,{stealth-stealth}](3,4) -- (3,5);

\draw [stealth-](4,1) -- (4,2);
\draw [red,thick,{stealth-stealth}](4,2) -- (4,3);
\draw [-stealth](4,3) -- (4,4);
\draw [-stealth](4,4) -- (4,5);

\draw [red,thick,{stealth-stealth}](5,1) -- (5,2);
\draw [red,thick,{stealth-stealth}](5,2) -- (5,3);
\draw [red,thick,{stealth-stealth}](5,3) -- (5,4);
\draw [-stealth](5,4) -- (5,5);
\end{tikzpicture}
 \\ [1ex] 
 \hline
 4-cube & 4 & 16 & 8 & 8$C$ & 30 &
   \begin{tikzpicture} [scale=0.4]
\draw [stealth-](0.5,1) -- (1,1);
\draw [stealth-](0.5,2) -- (1,2);
\draw [stealth-](0.5,3) -- (1,3);
\draw [stealth-](0.5,4) -- (1,4);
\draw [stealth-](0.5,5) -- (1,5);

\draw [-stealth](5,1) -- (5.5,1);
\draw [-stealth](5,2) -- (5.5,2);
\draw [-stealth](5,3) -- (5.5,3);
\draw [-stealth](5,4) -- (5.5,4);
\draw [-stealth](5,5) -- (5.5,5);

\draw [-stealth](1,5) -- (1,5.5);
\draw [-stealth](2,5) -- (2,5.5);
\draw [-stealth](3,5) -- (3,5.5);
\draw [-stealth](4,5) -- (4,5.5);
\draw [-stealth](5,5) -- (5,5.5);

\draw [-stealth](1,1) -- (1,0.5);
\draw [-stealth](2,1) -- (2,0.5);
\draw [-stealth](3,1) -- (3,0.5);
\draw [-stealth](4,1) -- (4,0.5);
\draw [-stealth](5,1) -- (5,0.5);

\draw [red,thick,{stealth-stealth}](1,1) -- (2,1);
\draw [red,thick,{stealth-stealth}](1,2) -- (2,2);
\draw [stealth-](1,3) -- (2,3);
\draw [stealth-](1,4) -- (2,4);
\draw [stealth-](1,5) -- (2,5);

\draw [-stealth](2,1) -- (3,1);
\draw [red,thick,{stealth-stealth}](2,2) -- (3,2);
\draw [red,thick,{stealth-stealth}](2,3) -- (3,3);
\draw [stealth-](2,4) -- (3,4);
\draw [stealth-](2,5) -- (3,5);

\draw [-stealth](3,1) -- (4,1);
\draw [-stealth](3,2) -- (4,2);
\draw [red,thick,{stealth-stealth}](3,3) -- (4,3);
\draw [red,thick,{stealth-stealth}](3,4) -- (4,4);
\draw [stealth-](3,5) -- (4,5);

\draw [-stealth](4,1) -- (5,1);
\draw [-stealth](4,2) -- (5,2);
\draw [-stealth](4,3) -- (5,3);
\draw [red,thick,{stealth-stealth}](4,4) -- (5,4);
\draw [red,thick,{stealth-stealth}](4,5) -- (5,5);

\draw [red,thick,{stealth-stealth}](1,1) -- (1,2);
\draw [-stealth](1,2) -- (1,3);
\draw [-stealth](1,3) -- (1,4);
\draw [-stealth](1,4) -- (1,5);

\draw [red,thick,{stealth-stealth}](2,1) -- (2,2);
\draw [red,thick,{stealth-stealth}](2,2) -- (2,3);
\draw [-stealth](2,3) -- (2,4);
\draw [-stealth](2,4) -- (2,5);

\draw [stealth-](3,1) -- (3,2);
\draw [red,thick,{stealth-stealth}](3,2) -- (3,3);
\draw [red,thick,{stealth-stealth}](3,3) -- (3,4);
\draw [-stealth](3,4) -- (3,5);

\draw [stealth-](4,1) -- (4,2);
\draw [stealth-](4,2) -- (4,3);
\draw [red,thick,{stealth-stealth}](4,3) -- (4,4);
\draw [red,thick,{stealth-stealth}](4,4) -- (4,5);

\draw [stealth-](5,1) -- (5,2);
\draw [stealth-](5,2) -- (5,3);
\draw [stealth-](5,3) -- (5,4);
\draw [red,thick,{stealth-stealth}](5,4) -- (5,5);
\end{tikzpicture}
 \\ [1ex] 
 \hline
\end{tabular}
\end{Large}
\end{table}


 The combinatorial types of 3-dimensional faces of $ASM_n$ include all the combinatorial types of 3-dimensional faces of $B_n$ together with that of the octahedron (See Figure \ref{octahedron}.)  These are all the combinatorial types of 3-dimensional 2-level polytopes \cite{enum2level}.  The combinatorial types of 4-dimensional faces of $ASM_n$ include all the combinatorial types of 4-dimensional faces of $B_n$ other than $B_3$ together with 5 other combinatorial types: $ASM_3$, cubical bipyramid, octahedral prism, polytope $X$ and polytope $Z$. 
  
Paffenholz \cite{Paffenholz1} points out that polytope $Y$ is a wedge over an edge of the base of a square pyramid.  It is also worth noting that polytope $Y$ is the order polytope of the zigzag poset with 4 points.  The paper \cite{Striker4} shows how the order polytope of the zigzag poset with an odd number $2n+1$ of points is a face of $ASM_n$ for which the order polytope of the zigzag poset with $2n$ points is a facet.

The Schlegel diagram of the 4-dimensional Polytope X is in Figure \ref{Schlegel polytope X}.  Note that the outside facet has the combinatorial type of an octahedron, while the other facets are combinatorially equivalent to tetrahedra and pyramids over squares.  Polytope $X$ has a 3-dimensional section that is affinely equivalent to a triangular prism.  In the Schlegel diagram, this section projects to a plane.

\begin{figure}[!h!]
\centering 
\begin{tikzpicture}[scale=0.7]
\draw[fill=black!30](2,0)circle[radius=.1];
\draw[fill=black!30](0,3)circle[radius=.1];
\draw[fill=black!30](4,3)circle[radius=.1];
\draw[fill=black!30](1,6)circle[radius=.1];
\draw[fill=black!30](4,6)circle[radius=.1];
\draw[fill=black!30](2,3.75)circle[radius=.1];
\draw[fill=black!30](2.5,9)circle[radius=.1];
\draw[fill=black!30](2.5,5.25)circle[radius=.1];
\draw [fill=black!20](2,0)--(0,3);
\draw [fill=black!20](2,0)--(4,3);
\draw [fill=black!20](2,0)--(1,6);
\draw [fill=black!20](4,3)--(4,6);
\draw [fill=black!20](0,3)--(4,3);
\draw [fill=black!20](0,3)--(1,6);
\draw [fill=black!20](1,6)--(4,6);
\draw [fill=black!20](4,3)--(2.5,9);
\draw [fill=black!20](1,6)--(2.5,9);
\draw [fill=black!20](0,3)--(2.5,9);
\draw [fill=black!20](4,6)--(2.5,9);
\draw [fill=black!20](2,0)--(4,6);
\draw [fill=black!20](0,3)--(2,3.75);
\draw [fill=black!20](0,3)--(2,3.75);
\draw [fill=black!20](4,3)--(2,3.75);
\draw [fill=black!20](2.5,9)--(2,3.75);
\draw [fill=black!20](4,6)--(2.5,5.25);
\draw [fill=black!20](1,6)--(2.5,5.25);
\draw [fill=black!20](2.5,9)--(2.5,5.25);
\draw [fill=black!20](2,0)--(2.5,5.25);
\draw [fill=black!20](2.5,5.25)--(2,3.75);
\draw [fill=black!20](2,0)--(2,3.75);
\end{tikzpicture}
\caption{\label{Schlegel polytope X}Schlegel diagram of polytope X.}
\end{figure}

The doubly directed graphs of Polytope X, Polytope Z, and the octahedral prism can be obtained from the doubly directed graph of the octahedron in Figure \ref{octahedron} by adding an ear.  For Polytope X, the ear is added connecting two distinct ears.  For Polytope Z, the added ear has distinct endpoints in the same ear. For the octahedral prism, a cycle is added with at most one point in common with one of ears of the octahedron.  These are all the ways to add an ear to the doubly directed graph of the octahedron.  

In the following discussion when we say that a polytope is ``not a product" we mean that it is not affinely equivalent to a product of two polytopes of dimension at least one.  

Polytope $Z$ has $9$ vertices and is not a product and is not centrally symmetric. Polytope $Z$ is a wedge over a facet of an octahedron. It has two octahedron and one triangular prism facets, each of which has $6$ vertices and is opposite to a triangle.  It has three facets that are pyramids over square bases, and each of these is opposite to a tetrahedron facet. 

A 4-dimensional polytope different from polytope $Z$ that is not a product and also has $9$ vertices is the cubic pyramid, so we see that the combinatorial types for the $d$-dimensional faces that are not products and not centrally symmetric and have $2^{d-1}+1$ vertices meeting the second upper bound of Theorem \ref{2connectedub} are not unique.    

There are 19 4-dimensional 2-level polytopes \cite{enum2level}. Four of those, namely $B_3$, a pyramid over an octahedron, the convex hull of the 0/1 vectors in $\mathbb{R}^4$ that have 1 or 2 ones, and the 4-dimensional cross polytope, are not faces of $ASM_n$. These four polytopes can be found in lines 7, 17, 25, and 31, respectively, of Table 1 in \cite{semidefinite4D}. $B_3$ was treated in Section 8, and the same proof shows that the
  pyramid over an octahedron, which has 7 vertices and 9 facets, is not a face of $ASM_n$. The convex hull of the 0/1 vectors in $\mathbb{R}^4$ that have 1 or 2 ones has 5 octahedron faces, each of which is opposite to a tetrahedron face (See line 25 of Table 1 in \cite{semidefinite4D}).  Since it has $10=2^{4-1}+2$ vertices and it is not a product and not centrally symmetric, we see that it violates the second upper bound of $2^{d-1}+1$ vertices from Theorem \ref{2connectedub}.  Therefore, it is not a face of $ASM_n$.    Finally, the 4-dimensional cross-polytope violates the upper bound of Theorem \ref{maxnumfacets} for the number of facets.  

 Finally, Figure \ref{bipyramidover4cube} shows nonuniqueness of ASMs that meet the bound of Theorem \ref{maxnumfacets} and first upper bound of Theorem \ref{2connectedub}.  All three elementary flow grids correspond to centrally symmetric 5-polytopes. The elementary flow grid on the left corresponds to a bipyramid over a 4-cube. It has $16=4(5-1)$ facets and $18=2^{5-1}+2$ vertices, matching the bound of Theorem \ref{maxnumfacets} and the first upper bound of Theorem \ref{2connectedub}. The elementary flow grid in the middle corresponds to a twisted prism (See \cite[Section 3.2]{2levelcombsettings} for a definition) over polytope Y.  It has 16 facets, equal to the upper bound of Theorem \ref{maxnumfacets}, but it has only 16 vertices.  This combinatorial type was discussed in \cite[Example 5.2]{Sanyal1}.  The elementary flow grid on the right has only 14 facets, because there is a set of two ears whose removal disconnects the doubly directed graph.  On the other hand, it has 18 vertices, showing that the bipyramid over a 4-cube is not alone in achieving this maximum over five-dimensional faces that are not products.  







\clearpage

\begin{thebibliography}{99}
\bibitem{Agnarsson1} G. Agnarsson, R. Greenlaw. ``Graph Theory: Modeling, Applications, and Algorithms,"
Pearson Prentice Hall, Upper Saddle River, New Jersey, 2007.
\bibitem{2levelcombsettings} M. Aprile, A. Cevallos, Y. Faenza. ``On 2-level polytopes arising in combinatorial settings."
{\it SIAM Journal on Discrete Mathematics}, 32:1857–1886, 2018.
\bibitem{Balinski1} M. Balinski and A. Russakoff. ``On the assignment polytope," {\it SIAM Review}, 16(4), 1974.
\bibitem{Behrend1} R. Behrend and V. Knight. ``Higher spin alternating sign matrices," {\it Electron. J. Combin,}, 14(1),
2007.
\bibitem{Billera2} L. J. Billera and A. Sarangarajan. ``All 0-1 polytopes are traveling salesman polytopes," {\it Com-
binatorica,} 16:175–188, 1996.
\bibitem{Birkhoff1} G. Birkhoff. ``Tres observaciones sobre el algebra lineal," {\it Univ. Nac. Tucum´an. Revista A.},
pages 147–151, 1946.
\bibitem{enum2level} A. Bohn, Y. Faenza, S. Fiorini, V. Fisikopoulos, M. Macchia, and K. Pashkovich. ``Enumeration
of 2-level polytopes," {\it Mathematical Programming Computation,} 11:173–210, 2019.
\bibitem{borsik2025prefixboundedmatrices} N. A. Borsik, A. Frank, P. Madarasi and T. Tak\'acs, ``Prefix-bounded matrices," arXiv preprint arXiv:2505.10739
\bibitem{Bressoud1} D.M. Bressoud. ``Proofs and Confirmations, The Story of the Alternating Sign Matrix Conjecture," Cambridge University Press, Cambridge UK, 1999.
\bibitem{BruGibI} R. A. Brualdi and P. M. Gibson. ``Convex polyhedra of doubly stochastic matrices I, applica-
tions of the permanent function," {\it J. Combin. Thy. A,} 22:194–230, 1977.
\bibitem{BruGibII} R. A. Brualdi and P. M. Gibson. ``Convex polyhedra of doubly stochastic matrices II, graph
of $\omega_n$," {\it J. Combin. Thy. B,} 22:175–198, 1977.
\bibitem{BruGibIII} R. A. Brualdi and P. M. Gibson. ``Convex polyhedra of doubly stochastic matrices III, affine
and combinatorial properties of $\omega_n$," {\it J. Combin. Thy. A,} 22:338–351, 1977.
\bibitem{Brualdi4} R.A. Brualdi and G. Dahl. ``Alternating sign matrices, extensions and related cones," {\it Advances
in Applied Mathematics,} 86:19–49, 2017.
\bibitem{Diestel} R. Diestel. ``Graph Theory," Springer Verlag, Berlin Germany, 2017.
\bibitem{Dinkelman1} E. Dinkelman. ``Combinatorial Aspects of the Alternating Sign Matrix Polytope," PhD thesis,
George Mason University, 2024.
\bibitem{Edmonds2} J. Edmonds, ``Ray \& me," Fulkerson Centennial Celebration, Cornell University, September 2024.
\bibitem{Fischer1} I. Fischer, ``A new proof of the refined alternating sign matrix theorem," {\it Journal of Combinatorial Theory, Ser. A} 114(2):253-264, 2007.
\bibitem{semidefinite4D} J. Gouveia, K. Pashkovich, R. Z. Robinson, and R.R. Thomas. ``Four-dimensional polytopes
of minimum positive semidefinite rank," {\it Journal of Combinatorial Theory A,} 145:184–226,
2017.
 G. \bibitem{Kuperberg1} Kuperberg. ``Another proof of the alternating sign matrix conjecture," {\it Intern. Math. Not.},
3:129 – 150, 1996.
\bibitem{MacLane1} S. MacLane. ``A combinatorial condition for planar graphs," {\it Fundamentae Mathematicae,}
28:22–32, 1937.
\bibitem{Striker4} K. M\'{e}sz\'{a}ros, A.H. Morales, and J. Striker. ``On flow polytopes, order polytopes, and certain
faces of the alternating sign matrix polytope," {\it Discrete Comput. Geom.}, 62:128–163, 2019.
\bibitem{Mills1} W.H. Mills, D. Robbins, and H. Rumsey. ``Alternating sign matrices and descending plane
partitions," {\it Journal of Combinatorial Theory Series A,} 34:340–359, 1983.
\bibitem{Paffenholz1} A. Paffenholz. ``Faces of Birkhoff polytopes," {\it The Electronic Journal of Combinatorics,} 22(1),
2015.
\bibitem{Sanyal1} R. Sanyal, A. Werner, and G. Ziegler. ``On Kalai's conjectures concerning centrally symmetric polytopes," {\it Discrete and Computational Geometry}, 41:183-198, 2009.
\bibitem{Striker1} J. Striker. ``The alternating sign matrix polytope," {\it Electron. J. Combin.}, 16, 2009.
\bibitem{Whitney1} H. Whitney. ``Non-separable and planar graphs," {\it Transactions of the American Mathematical Society}, 34(2):339-362, 1932.
\bibitem{Zeilberger1} D. Zeilberger. ``Proof of the alternating sign matrix conjecture," {\it Electron. J. Combin.}, 3(2),
1996.
\bibitem{Ziegler1} G. Ziegler. ``Lectures on Polytopes," Springer, New York, 2000



\end{thebibliography}
\end{document}